\documentclass[a4paper, 11pt, twoside]{scrartcl}

\usepackage{etex}
\usepackage{calc}
\usepackage{ifthen}
\usepackage{xargs}
\usepackage[T1]{fontenc}
\usepackage[utf8]{inputenc}
\usepackage[ngerman,british]{babel}
\usepackage{eurosym}
\usepackage{microtype}
\usepackage[intlimits]{amsmath}
\usepackage{amssymb}
\usepackage{mathtools}
\usepackage{bm}
\usepackage[framed, hyperref, thref, amsmath]{ntheorem}
\usepackage{units}
\usepackage{mleftright}
\usepackage[svgnames]{xcolor}
\usepackage[ntheorem, framemethod=tikz]{mdframed}
\allowdisplaybreaks
\usepackage[p]{libertine}
\usepackage[libertine,slantedGreek]{newtxmath} 
\usepackage[bb=px,bbscaled=0.94,cal=dutchcal,calscaled=0.95]{mathalfa}
\useosf
\usepackage{anyfontsize}
\usepackage[singlespacing]{setspace}
\usepackage[headsepline,footsepline,plainfootsepline,headsepline=1pt,footsepline=1pt]{scrlayer-scrpage}
\usepackage[shortlabels]{enumitem}
\setlist{listparindent=\parindent,parsep=\parskip}
\usepackage{graphicx}
\usepackage{placeins}
\usepackage{caption}
\DeclareCaptionFont{gray}{\color{black}}
\usepackage{subfig}
\usepackage{pstricks-add}
\usepackage{pst-plot}
\usepackage{pst-func}
\psset{unit=1cm}
\usepackage{array}
\usepackage{colortbl}
\usepackage{longtable}
\usepackage{ltablex}
\usepackage{varioref}
\usepackage{hyperref}
\hypersetup{hidelinks}
\hypersetup{breaklinks=true}
\hypersetup{bookmarksopen,bookmarksopenlevel=1}
\addto\extrasbritish{%
}
\usepackage{breakurl}
\usepackage{bookmark}
\KOMAoption{BCOR}{5mm}
\KOMAoption{DIV}{classic}
\KOMAoption{headinclude}{false}
\KOMAoption{footinclude}{false}
\KOMAoption{pagesize}{auto}
\recalctypearea
\KOMAoption{headings}{small}
\KOMAoption{numbers}{endperiod}
\newcommand{\AUTHORsave}{}
\newcommand{\AUTHOR}[1]{\renewcommand{\AUTHORsave}{#1}}
\newcommand{\SHORTAUTHORsave}{}
\newcommand{\SHORTAUTHOR}[1]{%
  \ifthenelse{\equal{#1}{}}%
  {\renewcommand{\SHORTAUTHORsave}{\AUTHORsave}}%
  {\renewcommand{\SHORTAUTHORsave}{#1}}%
}
\newcommand{\TITLEsave}{}
\newcommand{\TITLE}[1]{\renewcommand{\TITLEsave}{#1}}
\newcommand{\SHORTTITLEsave}{}
\newcommand{\SHORTTITLE}[1]{%
  \ifthenelse{\equal{#1}{}}%
  {\renewcommand{\SHORTTITLEsave}{\TITLEsave}}%
  {\renewcommand{\SHORTTITLEsave}{#1}}%
}
\newcommand{\INSTITUTEsave}{}
\newcommand{\INSTITUTE}[1]{\renewcommand{\INSTITUTEsave}{#1}}
\newcommand{\CORRESPONDENCEsave}{}
\newcommand{\CORRESPONDENCE}[1]{\renewcommand{\CORRESPONDENCEsave}{#1}}
\newcommand{\ABSTRACTsave}{}
\newcommand{\ABSTRACT}[1]{\renewcommand{\ABSTRACTsave}{#1}}
\newcommand{\KEYWORDSsave}{}
\newcommand{\KEYWORDS}[1]{\renewcommand{\KEYWORDSsave}{#1}}
\newcommand{\AMSCLASSsave}{}
\newcommand{\AMSCLASS}[1]{\renewcommand{\AMSCLASSsave}{#1}}
\newcommand{\Inst}[1]{\textsuperscript{#1}}
\newcommand{\email}[1]{\href{mailto:#1}{\nolinkurl{#1}}}
\newcommand{\TitleHeader}{%
  \noindent%
  \rule{\linewidth}{1pt}\\[10pt]%
  {\sffamily\bfseries\Large\TITLEsave}\\[10pt]%
  {\sffamily\bfseries\small\AUTHORsave}\\%
  {\sffamily\scriptsize\INSTITUTEsave}\\[15pt]%
  \parbox{\linewidth}{\sffamily\scriptsize {\bfseries
      Correspondence}\\\CORRESPONDENCEsave}\\[10pt]%
  \rule{\linewidth}{1pt}\\[15pt]%
  \hspace*{40pt}%
  \parbox{\linewidth-40pt}{%
    \small%
    {\sffamily\bfseries Abstract.} \quad \ABSTRACTsave%
  }\\[5pt]%
  \hspace*{40pt}%
  \parbox{\linewidth-40pt}{%
    \small%
    {\sffamily\bfseries Keywords.} \quad \KEYWORDSsave%
  }\\[5pt]%
  \hspace*{40pt}%
  \parbox{\linewidth-40pt}{%
    \small%
    {\sffamily\bfseries AMS subject classification.} \quad \AMSCLASSsave%
  }\\[5pt]%
  \rule{\linewidth}{1pt}%
}
\clearscrheadfoot
\ohead{\pagemark}
\cehead{\footnotesize\itshape\scshape \SHORTAUTHORsave}
\cohead{\footnotesize\itshape\scshape \SHORTTITLEsave}

\newcommand{\qed}{}
\newlength{\frametopsep}
\colorlet{ThmFrame}{CornflowerBlue!25}
\colorlet{ThmBody}{CornflowerBlue!5}
\mdfdefinestyle{thmframe}{innerleftmargin=5pt, innerrightmargin=5pt,
  innermargin=-5.5pt, outermargin=-5.5pt, skipabove=\frametopsep,
  skipbelow=\frametopsep, 
  roundcorner=5pt, linewidth=0.5pt, linecolor=ThmFrame,
  backgroundcolor=ThmBody, splittopskip=12.5pt, splitbottomskip=5pt}
\makeatletter
\newtheoremstyle{plainheaderbreak}
  {\item[]{\theorem@headerfont ##1\ ##2\theorem@separator}\hskip\labelsep}%
  {\item[]{\theorem@headerfont ##1\ ##2\
    (##3)}\theorem@separator\hskip\labelsep}
\makeatother
\theoremheaderfont{\sffamily\bfseries\scshape}
\theorembodyfont{}
\theoremstyle{plainheaderbreak}
\theoremseparator{.}
\theoremprework{}
\theorempostwork{}
\theoremframepreskip{\frametopsep}
\theoremframepostskip{\frametopsep}
\theoreminframepreskip{0pt}
\theoreminframepostskip{0pt}
\newmdtheoremenv[style=thmframe]{Definition}{Definition}[section]
\newmdtheoremenv[style=thmframe]{Problem}[Definition]{Problem}
\newmdtheoremenv[style=thmframe]{Assumption}[Definition]{Assumption}
\theoremheaderfont{\sffamily\bfseries\scshape\upshape}
\theorembodyfont{\itshape}
\newmdtheoremenv[style=thmframe]{Theorem}[Definition]{Theorem}
\newmdtheoremenv[style=thmframe]{Proposition}[Definition]{Proposition}
\newmdtheoremenv[style=thmframe]{Lemma}[Definition]{Lemma}
\newmdtheoremenv[style=thmframe]{Corollary}[Definition]{Corollary}
\theoremheaderfont{\sffamily\bfseries\scshape\upshape}
\theorembodyfont{\sffamily}
\theoremstyle{break}
\newmdtheoremenv[style=thmframe]{Algorithm}[Definition]{Algorithm}
\newlength{\unframetopsep}
\theoremheaderfont{\itshape}
\theorembodyfont{}
\theoremstyle{plain}
\theoremseparator{.}
\theorempreskip{\unframetopsep}
\theorempostskip{\unframetopsep}
\theoremprework{%
  \renewcommand{\qed}{%
    \hspace*{\fill}\nolinebreak[3]%
    \nopagebreak[3]\hspace*{\fill}{\raisebox{0.25pt}{\scalebox{0.9}{$\medcirc$}}}}}
\theorempostwork{}

\theoremprework{%
  \renewcommand{\qed}{%
    \hspace*{\fill}\nolinebreak[3]%
    \nopagebreak[3]\hspace*{\fill}{\raisebox{0.25pt}{\scalebox{0.9}{$\medcirc$}}}}}
\theorempostwork{}

\theoremprework{%
  \renewcommand{\qed}{%
    \hspace*{\fill}\nolinebreak[3]%
    \nopagebreak[3]\hspace*{\fill}{\raisebox{0.25pt}{\scalebox{0.9}{$\medcirc$}}}}}
\theorempostwork{}
\newtheorem{Remark}[Definition]{Remark}
\makeatletter
\newtheoremstyle{nonumberplainof}%
{\item[\theorem@headerfont\hskip\labelsep ##1\theorem@separator]}%
{\item[\theorem@headerfont\hskip \labelsep ##1\ of\ ##3\theorem@separator]}
\makeatother
\theoremstyle{nonumberplainof}
\theoremprework{%
    \renewcommand{\qed}{%
        \hspace*{\fill}\nolinebreak[3]%
        \nopagebreak[3]\hspace*{\fill}{$\square$}}}
\theorempostwork{}
\newtheorem{Proof}{Proof}

\newcommand{\pers}[1]{%
    {\scshape#1}%
}

\newcommand{\PCholesky}{\pers{Cholesky}}

\newcommand{\PEuclid}{\pers{Euclid}}

\newcommand{\PFenchel}{\pers{Fenchel}}

\newcommand{\PFourier}{\pers{Fourier}}

\newcommand{\PGauss}{\pers{Gauß}}

\newcommand{\PHermit}{\pers{Hermit}}

\newcommand{\PHilbert}{\pers{Hilbert}}

\newcommand{\PKronecker}{\pers{Kronecker}}
\newcommand{\PKrylov}{\pers{Krylov}}

\newcommand{\PLanczos}{\pers{Lanczos}}

\newcommand{\PLegendre}{\pers{Legendre}}
\newcommand{\PLipschitz}{\pers{Lipschitz}}

\newcommand{\PParseval}{\pers{Parseval}}

\newcommand{\PPoisson}{\pers{Poisson}}

\newcommand{\PRademacher}{\pers{Rademacher}}

\newcommand{\PRitz}{\pers{Ritz}}

\newcommand{\PSchatten}{\pers{Schat\-ten}}

\newcommand{\PSchmidt}{\pers{Schmidt}}

\newcommand{\PSobolev}{\pers{Sobolev}}

\newcommand{\PTikhonov}{\pers{Tik\-ho\-nov}}

\newcommand{\eg}{%
    e.g.%
}
\newcommand{\ie}{%
    i.e.%
}

\newcommand{\cf}{%
    cf.%
}

\newcommand{\Vek}[1]{%
    \ensuremath{\bm{#1}}%
}
\newcommand{\Mat}[1]{%
    \ensuremath{\bm{#1}}%
}
\newcommand{\R}{%
    \ensuremath{\mathbb{R}}%
}

\newcommand{\N}{%
    \ensuremath{\mathbb{N}}%
}

\newcommand{\C}{%
    \ensuremath{\mathbb{C}}%
}
\newcommand{\iProd}[2]{%
    \left\langle#1,#2\right\rangle%
}

\newcommand{\iProdn}[2]{%
    \langle#1,#2\rangle%
}
\newcommand{\iProdb}[2]{%
    \bigl\langle#1,#2\bigr\rangle%
}

\newcommand{\absn}[1]{%
    |\hspace{1pt}#1\hspace{1pt}|%
}
\newcommand{\absb}[1]{%
    \bigl|\hspace{1pt}#1\hspace{1pt}\bigr|%
}

\newcommand{\absB}[1]{%
    \Bigl|\hspace{1pt}#1\hspace{1pt}\Bigr|%
}

\newcommand{\pNorm}[1]{%
    \left|\left|\hspace{1pt}#1\hspace{1pt}\right|\right|%
}

\newcommand{\pNormn}[1]{%
    ||\hspace{1pt}#1\hspace{1pt}||%
}

\newcommandx{\lpNorm}[3][3=\empty]{%
    \left\lVert\hspace{1pt}#1\hspace{1pt}\right\rVert_{\ell^{#2}%
        \ifthenelse{\equal{#3}{}}{}{\mleft(#3\mright)}%
    }%
}
\newcommandx{\lpNormm}[3][3=\empty]{%
    \mleft\lVert\hspace{1pt}#1\hspace{1pt}\mright\rVert_{\ell^{#2}%
        \ifthenelse{\equal{#3}{}}{}{\mleft(#3\mright)}%
    }%
}
\newcommandx{\lpNormn}[3][3=\empty]{%
    \lVert\hspace{1pt}#1\hspace{1pt}\rVert_{\ell^{#2}%
        \ifthenelse{\equal{#3}{}}{}{\mleft(#3\mright)}%
    }%
}
\newcommandx{\lpNormb}[3][3=\empty]{%
    \bigl\lVert\hspace{1pt}#1\hspace{1pt}\bigr\rVert_{\ell^{#2}%
        \ifthenelse{\equal{#3}{}}{}{\mleft(#3\mright)}%
    }%
}
\newcommandx{\lpNormbb}[3][3=\empty]{%
    \biggl\lVert\hspace{1pt}#1\hspace{1pt}\biggr\rVert_{\ell^{#2}%
        \ifthenelse{\equal{#3}{}}{}{\mleft(#3\mright)}%
    }%
}
\newcommandx{\lpNormB}[3][3=\empty]{%
    \Bigl\lVert\hspace{1pt}#1\hspace{1pt}\Bigr\rVert_{\ell^{#2}%
        \ifthenelse{\equal{#3}{}}{}{\mleft(#3\mright)}%
    }%
}
\newcommandx{\lpNormBB}[3][3=\empty]{%
    \Biggl\lVert\hspace{1pt}#1\hspace{1pt}\Biggr\rVert_{\ell^{#2}%
        \ifthenelse{\equal{#3}{}}{}{\mleft(#3\mright)}%
    }%
}
\newcommandx{\LpNorm}[3][3=\empty]{%
    \left\lVert\hspace{1pt}#1\hspace{1pt}\right\rVert_{L^{#2}%
        \ifthenelse{\equal{#3}{}}{}{\mleft(#3\mright)}%
    }%
}
\newcommandx{\LpNormm}[3][3=\empty]{%
    \mleft\lVert\hspace{1pt}#1\hspace{1pt}\mright\rVert_{L^{#2}%
        \ifthenelse{\equal{#3}{}}{}{\mleft(#3\mright)}%
    }%
}
\newcommandx{\LpNormn}[3][3=\empty]{%
    \lVert\hspace{1pt}#1\hspace{1pt}\rVert_{L^{#2}%
        \ifthenelse{\equal{#3}{}}{}{\mleft(#3\mright)}%
    }%
}
\newcommandx{\LpNormb}[3][3=\empty]{%
    \bigl\lVert\hspace{1pt}#1\hspace{1pt}\bigr\rVert_{L^{#2}%
        \ifthenelse{\equal{#3}{}}{}{\mleft(#3\mright)}%
    }%
}
\newcommandx{\LpNormbb}[3][3=\empty]{%
    \biggl\lVert\hspace{1pt}#1\hspace{1pt}\biggr\rVert_{L^{#2}%
        \ifthenelse{\equal{#3}{}}{}{\mleft(#3\mright)}%
    }%
}
\newcommandx{\LpNormB}[3][3=\empty]{%
    \Bigl\lVert\hspace{1pt}#1\hspace{1pt}\Bigr\rVert_{L^{#2}%
        \ifthenelse{\equal{#3}{}}{}{\mleft(#3\mright)}%
    }%
}
\newcommandx{\LpNormBB}[3][3=\empty]{%
    \Biggl\lVert\hspace{1pt}#1\hspace{1pt}\Biggr\rVert_{L^{#2}%
        \ifthenelse{\equal{#3}{}}{}{\mleft(#3\mright)}%
    }%
}
\newcommand{\Fourier}{%
    \mathop{\kern0pt\mathcal{F}}\nolimits%
}

\newcommand{\Fresnel}{%
    \mathop{\kern0pt\mathcal{E}}\nolimits%
}

\newcommand{\Laplace}{%
    \mathop{\kern0pt\mathcal{L}}\nolimits%
}
\newcommand{\LCT}{%
    \mathop{\kern0pt\mathcal{C}}\nolimits%
}

\newcommand{\z}{%
    \mathop{\kern0pt\mathcal{Z}}\nolimits%
}
\newcommand{\Landau}{%
    \mathop{\kern0pt\mathcal{O}}\nolimits%
}

\newcommand{\e}{%
    \ensuremath{\mathrm{e}}%
}
\newcommand{\I}{%
    \ensuremath{\mathrm{i}}%
}

\newcommand{\argmin}{%
    \ensuremath{\operatorname*{argmin}}%
}

\newcommand{\diag}{%
    \ensuremath{\operatorname{diag}}%
}

\newcommand{\maximize}{%
    \ensuremath{\operatorname*{maximize}}%
  }
  \newcommand{\minimize}{%
    \ensuremath{\operatorname*{minimize}}%
}

\newcommand{\prox}{%
    \ensuremath{\operatorname{prox}}%
}

\newcommand{\rank}{%
    \ensuremath{\operatorname{rank}}%
}

\newcommand{\sgn}{%
    \ensuremath{\operatorname{sgn}}%
}

\newcommand{\tr}{%
    \ensuremath{\operatorname{tr}}%
}

\newcommand{\T}{%
    \ensuremath{\mathrm{T}}%
}

\newcommand{\ind}[2][]{%
    \ifthenelse{\equal{#1}{}}{%
        \index{#2}\margin{#2}%
    }{%
        \index{#2}\margin{#1}%
    }%
}
\newcommand{\addmathskip}[1][5pt]{%
    \vspace*{#1}%
}
\newcommand{\submathskip}[1][-5pt]{%
    \vspace*{#1}%
}
\newlength{\subalignskip}
\newlength{\subalignaboveskip}
\newlength{\subalignbelowskip}
\newlength{\fsmallskip}
\newlength{\fskip}
\newlength{\subintertextskip}
\newlength{\addintertextskip}
\newlength{\OverlineExtension}

\newlength{\optmarginbox}
\newcommand{\margin}[1]{%
    \marginpar[\raggedleft\raisebox{0cm}{%
        \parbox{\marginparwidth}{\raggedleft\sffamily\tiny\color{Gray}#1}}%
    ]{\raggedright\raisebox{0cm}{%
            \parbox{\marginparwidth}{\raggedright\sffamily\tiny\color{Gray}#1}}}%
}
\newcommand{\marginBox}[1]{%
    \settowidth{\optmarginbox}{#1}
    \addtolength{\optmarginbox}{8pt}
    \ifthenelse{\lengthtest{\optmarginbox < \marginparwidth}}{
        \marginpar[\raggedleft\raisebox{0cm}{%
                \psframebox[linecolor=LightGray, linearc=5pt, cornersize=absolute,%
                framesep=3pt, fillcolor=GhostWhite, fillstyle=solid]%
                {\sffamily\tiny\color{Gray}#1}}%
        ]{\raggedright\raisebox{0cm}{%
                \psframebox[linecolor=LightGray, linearc=5pt, cornersize=absolute,%
                framesep=3pt, fillcolor=GhostWhite, fillstyle=solid]%
                {\sffamily\tiny\color{Gray}#1}}}%
    }{%
        \marginpar[\raggedleft\raisebox{0cm}{%
                \psframebox[linecolor=LightGray, linearc=5pt, cornersize=absolute,%
                framesep=3pt, fillcolor=GhostWhite, fillstyle=solid]%
                {\parbox{\marginparwidth-8pt}{\raggedleft\sffamily\tiny\color{Gray}#1}}}%
        ]{\raggedright\raisebox{0cm}{%
                \psframebox[linecolor=LightGray, linearc=5pt, cornersize=absolute,%
                framesep=3pt, fillcolor=GhostWhite, fillstyle=solid]%
                {\parbox{\marginparwidth-8pt}{\raggedright\sffamily\tiny\color{Gray}#1}}}}%
    }%
}
\newcommand{\marginBoxA}[1]{%
    \settowidth{\optmarginbox}{#1}
    \addtolength{\optmarginbox}{8pt}
    \ifthenelse{\lengthtest{\optmarginbox < \marginparwidth}}{
        \marginpar[\raggedleft\raisebox{0cm}{%
                \psframebox[linecolor=Maroon, linearc=5pt, cornersize=absolute,%
                framesep=3pt, fillcolor=MistyRose, fillstyle=solid]%
                {\sffamily\tiny\color{Maroon}#1}}%
        ]{\raggedright\raisebox{0cm}{%
                \psframebox[linecolor=Maroon, linearc=5pt, cornersize=absolute,%
                framesep=3pt, fillcolor=MistyRose, fillstyle=solid]%
                {\sffamily\tiny\color{Maroon}#1}}}%
    }{%
        \marginpar[\raggedleft\raisebox{0cm}{%
                \psframebox[linecolor=Maroon, linearc=5pt, cornersize=absolute,%
                framesep=3pt, fillcolor=MistyRose, fillstyle=solid]%
                {\parbox{\marginparwidth-8pt}{\raggedleft\sffamily\tiny\color{Maroon}#1}}}%
        ]{\raggedright\raisebox{0cm}{%
                \psframebox[linecolor=Maroon, linearc=5pt, cornersize=absolute,%
                framesep=3pt, fillcolor=MistyRose, fillstyle=solid]%
                {\parbox{\marginparwidth-8pt}{\raggedright\sffamily\tiny\color{Maroon}#1}}}}%
    }%
}

\newcommand{\bq}[1]{%
    `#1'%
}

\AUTHOR{Robert \pers{Beinert}\Inst{1} and Kristian
  \pers{Bredies}\Inst{1}}

\SHORTAUTHOR{R.~Beinert and K.~Bredies}

\TITLE{Tensor-Free Proximal Methods for Lifted Bilinear/Quadratic
  Inverse Problems with Applications to Phase Retrieval}

\SHORTTITLE{Tensor-free methods for bilinear/quadratic inverse problems}

\INSTITUTE{\Inst{1}\parbox[t]{0.95\linewidth}{Institut für Mathematik
    und Wissenschaftliches Rechnen\\
    \pers{Karl}-\pers{Franzens}-Universität Graz\\
    Heinrichstraße 36\\
    8010 Graz, Austria
  }
}

\CORRESPONDENCE{R. \pers{Beinert}: \email{robert.beinert@uni-graz.at}\\
  K. \pers{Bredies}: \email{kristian.bredies@uni-graz.at}}

\ABSTRACT{ We propose and study a class of novel algorithms that aim
  at solving bilinear and quadratic inverse problems.  Using a convex
  relaxation based on tensorial lifting, and applying first-order
  proximal algorithms, these problems could be solved numerically by
  singular value thresholding methods.  However, a direct realization
  of these algorithms for, e.g., image recovery problems is often
  impracticable, since computations have to be performed on the
  tensor-product space, whose dimension is usually tremendous.  To
  overcome this limitation, we derive tensor-free versions of common
  singular value thresholding methods by exploiting low-rank
  representations and incorporating an augmented \PLanczos\ process.
  Using a novel reweighting technique, we further improve the
  convergence behavior and rank evolution of the iterative algorithms.
  Applying the method to the two-dimensional masked \PFourier\ phase
  retrieval problem, we obtain an efficient recovery method.
  Moreover, the tensor-free algorithms are flexible enough to
  incorporate a-priori smoothness constraints that greatly improve the
  recovery results.}


\AMSCLASS{45Q05, 65F10, 65R32, 65K10}

\begin{document}
\thispagestyle{empty}
\TitleHeader

\section{Introduction}
\label{sec:introduction}

The theory of inverse problems is nowadays one of the main tools to
deal with recovery problems in medicine, engineering, and life
sciences.  The real-world applications of this theory embrace for
instance computed tomography, magnetic resonance imaging, and
deconvolution problems in microscopy, see \cite{BB98, MS12, Ram05,
  SW13, Uhl03, Uhl13}.  Besides these recent monographs, which are
only a small selection, there exist many further publications about
applications, regularization, and numerical solvers.

In this paper, we consider the subclass of bilinear and quadratic
inverse problems \cite{BB18}.  Problem formulations of these kinds
originate form real-world applications in imaging and physics
\cite{SGG+09} like blind deconvolution \cite{BS01,JR06},
deautoconvolution \cite{GH94, GHB+14, ABHS16}, phase retrieval
\cite{DF87, Mil90, SSD+06}, parallel imaging in MRI \cite{BBM+04}, and
parameter identification in EIT \cite{MS12}.  Although these problems
may be seen as non-linear inverse problems and can be solved by
appropriate non-linear solvers, the question arises if we can exploit
the specific bilinear or quadratic structure of these problem.

One of the most famous approaches of this kind is PhaseLift
\cite{CSV13}, where the generic phase retrieval problem is lifted to a
linear inverse problem with rank-one constraint using the universal
property of the tensor product.  After a relaxation the lifted problem
is solved by semi-definite programming.  From the theoretical side,
one has proved that this methodology yields the true solution with
very high probability.  Based on PhaseLift, we will derive convex
lifting methods for general bilinear and quadratic inverse problems.

Similarly to PhaseLift, the relaxations here require the minimization
of a nuclear norm functional on the tensor product, which can be done
by applying proximal minimization methods, see for instance
\cite{CP16} and references therein.  Since the dimension of the lifted
problem can become huge, the question arises how one can
implement these methods numerically.  Therefore, we develop
tensor-free, equivalent versions that can be performed in an efficient
and memory-saving manner.  The main benefit is here that the
tensor-free reformulations of the required tensorial operations are
performed exactly without additional error.  To improve the
convergence and solution behaviour, we additionally propose a novel
reweighting technique

The paper is organized as follows: In
\autoref{sec:bilin-inverse-probl}, we introduce the considered
bilinear and quadratic inverse problems in more detail.  The focus
here lies on the bilinear setting since quadratic formulations are
based on underlying bilinear structures.  Based on the universal
property of bilinear mappings and the nuclear norm heuristic, we then
derive a relaxed convex minimization problem with linear lifted
forward operator.  To stabilize the lifted problem regarding noise and
measurement errors, we additionally consider a \PTikhonov\ approach.

In \autoref{sec:sing-value-thresh}, we develop a proximal solver based
on the first-order primal-dual method of \pers{Chambolle} and
\pers{Pock} \cite{CP11} to solve the lifted problem numerically.  The
primal-dual iteration is here only one explicit example and can be
replaced by any other proximal method.  In so doing, we obtain a
singular value thresholding depending on the actual \PHilbert\ spaces
building the domain of the original problem.  Although the tensorial
lifting allows us to apply linear methods, the dimension of the
relaxed minimization problem becomes tremendous.

To overcome this issue, we derive a tensor-free representation of the
suggested algorithm.  The efficient computation of the required
singular value thresholding is here ensured by exploiting an
orthogonal power iteration or, alternatively, an augmented \PLanczos\
process, see \autoref{sec:circ-tens-prod}.  Moreover, in
\autoref{sec:reduc-rank-adapt}, we introduce a novel \PHilbert\ space
reweighting to promote low-rank iterations and solutions.  We complete
the paper with a numerical study, where we consider the masked
\PFourier\ phase retrieval, see \autoref{sec:mask-phase-retr}.

The contribution of the paper is twofold: Firstly, we derive a novel
tensor-free proximal algorithm based on a convex lifting and
relaxation.  Secondly, we introduce a novel phase retrieval technique,
which solves high-dimensional instances of the masked phase retrieval
problem.  Moreover, our approach allows us to incorporate smoothness
constraints or relation between different pixels to improve the
convergence of the algorithm.

\section{Convex liftings of bilinear and quadratic inverse problems}
\label{sec:bilin-inverse-probl}

Bilinear and quadratic problem formulations arise in a wide range of
applications in imaging and physics \cite{SGG+09} like blind
deconvolution \cite{BS01,JR06}, deautoconvolution \cite{GH94, GHB+14,
  ABHS16}, phase retrieval \cite{DF87, Mil90, SSD+06}, parallel
imaging in MRI \cite{BBM+04}, and parameter identification in EIT
\cite{MS12}.  Since we are mainly interested in computing a numerical
solution, we restrict ourselves to \emph{finite-dimensional bilinear
problems} of the form
\begin{equation}
  \label{eq:gen-bilin-prob}
  \mathcal B(\Vek u, \Vek v) = \Vek g^\dagger,
  \tag{$\mathfrak B$}
  \addmathskip
\end{equation}
where $\mathcal B$ is a bilinear operator from
$\mathcal H_1 \times \mathcal H_2$ into $\mathcal K$, and where
$\mathcal H_1 \simeq \R^{N_1}$, $\mathcal H_2 \simeq \R^{N_2}$, and
$\mathcal K \simeq \R^M$ are real finite-dimensional \PHilbert\ spaces
equipped with some inner product.  Simultaneously, we study
\emph{finite-dimensional quadratic problems} of the form
\begin{equation}
  \label{eq:gen-quad-prob}
  \mathcal Q(\Vek u) = \Vek g^\dagger.
  \tag{$\mathfrak Q$}
\end{equation}
A quadratic operator $\mathcal Q \colon \mathcal H \to \mathcal K$ is
here the restriction of an associate bilinear operator
$\mathcal B_{\mathcal Q} : \mathcal H \times \mathcal H \to \mathcal
K$
to its diagonal and is thus given by
$\mathcal Q(\Vek u) \coloneqq \mathcal B_{\mathcal Q}(\Vek u, \Vek
u)$.
Again, the \PHilbert\ spaces $\mathcal H$ and $\mathcal K$ are finite
dimensional.  Without loss of generality, one may assume that the
associated bilinear operator is symmetric, \ie\
$\mathcal B_{\mathcal Q} (\Vek u, \Vek v) = \mathcal B_{\mathcal Q}
(\Vek v, \Vek u)$
for all $\Vek u$, $\Vek v \in \mathcal H$, which can be enforced by
considering
$\nicefrac12 \, \mathcal B_{\mathcal Q} (\Vek u, \Vek v) + \nicefrac12
\, \mathcal B_{\mathcal Q} (\Vek v, \Vek u)$.
However, we will make no use of this restriction.

For the sake of simplicity, we use the following notation and write an
arbitrary vector $\Vek u \in \mathcal H_1$ in the form
$\Vek u \coloneqq (u_n)_{n=0}^{N_1-1} \in \R^{N_1}$, which may be
interpreted as the coefficient vector with respect to some finite
basis.  The inner product of $\mathcal H_1$ can now be stated in the
form
$\iProdn{\cdot}{\cdot}_{\mathcal H_1} \coloneqq \iProdn{\Mat H_1 \,
  \cdot}{\cdot} = \iProdn{\cdot}{\Mat H_1 \, \cdot}$,
where $\Mat H_1$ is some symmetric, positive definite matrix in
$\R^{N_1 \times N_1}$, and where the inner products on the right-hand
side denote the usual \PEuclid{ian} inner product
$\iProd{\Vek u}{\Vek v} \coloneqq \Vek v^* \Vek u$.  Here $\cdot^*$
labels the transposition of a vector or matrix.  For the remaining
spaces $\mathcal H_2$ and $\mathcal K$, we proceed in the same manner
with the associate matrices $\Mat H_2 \in \R^{N_2 \times N_2}$ and
$\Mat K \in \R^{M \times M}$, respectively.

Although we restrict ourselves to the real-valued setting, all
following algorithms and statements remain valid for the
complex-valued setting, where one considers sesquilinear mappings
$\mathcal B \colon \mathcal H_1 \times \mathcal H_2 \to \mathcal K$
with $\mathcal H_1 \simeq \C^{N_1}$, $\mathcal H_2 \simeq \C^{N_2}$,
and $\mathcal K \simeq \C^{M}$.
Changing the notation to
$\Vek u \coloneqq (u_n)_{n=0}^{N_1-1} \in \C^{N_1}$, replacing the
property \bq{symmetric} by \bq{\PHermit{ian},} and using the real part
of the inner products, \eg\ $\iProd{\Vek u}{\Vek v} = \Re[ \Vek v^*
\Vek u]$, where $\cdot^*$ is the transposition and conjugation, all
considerations translate one to one.  In the complex case, the
associate matrices $\Mat H_1$, $\Mat H_2$, and $\Mat K$ may also be
complex-valued, \PHermit{ian}, and positive definite.

Inspired by PhaseLift \cite{CESV13,CSV13}, which exploits the solution
strategy developed for matrix completion problems \cite{CCS10,MGC11},
we suggest to tackle the general bilinear and quadratic inverse
problem \eqref{eq:gen-bilin-prob} and \eqref{eq:gen-quad-prob} by
convex liftings and relaxations.  Our approach is here based on the
so-called universal property of the tensor product with respect to
bilinear mappings, see for instance \cite[Theorem~2.9]{Rya02}.  In the
finite-dimensional setting, the lifting may be stated as follows.

\begin{Theorem}[Bilinear lifting, \cite{Rya02}]
  \label{the:bilin-lifting}
  For every bilinear mapping
  $\mathcal B \colon \mathcal H_1 \times \mathcal H_2 \rightarrow
  \mathcal K$,
  there exists a unique linear mapping
  $\breve{\mathcal B} \colon \mathcal H_1 \otimes \mathcal H_2
  \rightarrow \mathcal K$
  such that
  $\breve{\mathcal B}(\Vek u \otimes \Vek v) = \mathcal B(\Vek u, \Vek
  v)$.
\end{Theorem}

In the finite-dimensional setting considered by us, the tensor product
$\mathcal H_1 \otimes \mathcal H_2$ can be simply identified with the
matrix space $\R^{N_2 \times N_1}$, where the rank-one tensor
$\Vek u \otimes \Vek v$ corresponds to the matrix
$\Vek v \Vek u^* = (v_{n_2} u_{n_1})_{n_2=0,n_1=0}^{N_2-1,N_1-1}$.
The universal property is also applicable to quadratic mappings, where
the tensor product $\mathcal H \otimes \mathcal H$ is restricted to
the subspace of symmetric tensors
$\mathcal H \otimes_{\mathrm{sym}} \mathcal H$ which is isomorphic to
the space of symmetric matrices.

\begin{Corollary}[Quadratic lifting]
  \label{cor:quad-lifting}
  For every bounded quadratic mapping
  $\mathcal Q \colon \mathcal H \times \rightarrow \mathcal K$, there
  exists a unique linear mapping
  $\breve{\mathcal Q} \colon \mathcal H \otimes_{\mathrm{sym}}
  \mathcal H \rightarrow \mathcal K$
  such that
  $\breve{\mathcal Q}(\Vek u \otimes \Vek u) = \mathcal Q(\Vek u)$.
\end{Corollary}

\begin{Proof}
  Since the lifting $\breve{\mathcal B}$ of the associate bilinear
  mapping $\mathcal B$ is unique, the restriction
  $\breve{\mathcal Q} \coloneqq \breve{\mathcal B} \vert_{\mathcal H
    \otimes_{\mathrm{sym}} \mathcal H}$
  yields a unique linear mapping with the asserted properties.  \qed
\end{Proof}

Due to the uniqueness of the lifting, the bilinear inverse problem
\eqref{eq:gen-bilin-prob} and the quadratic inverse problem
\eqref{eq:gen-quad-prob} are equivalent to the linear
inverse problems
\begin{align}
  \label{eq:lift-bilin-prob}
  \breve{\mathcal B} (\Vek w) = \Vek g^\dagger
  &\qquad\text{subject to}\qquad
    \rank(\Vek w) \leq 1
    \tag{$\breve{\mathfrak B}$}
  \\  
   \shortintertext{and}
  \label{eq:lift-quad-prob}
  \breve{\mathcal Q} (\Vek w) = \Vek g^\dagger
  &\qquad\text{subject to}\qquad
    \rank(\Vek w) \leq 1,
    \quad\Vek w \succeq \Vek 0,
    \tag{$\breve{\mathfrak Q}$}
\end{align}
where the additional constraint $\Vek w \succeq \Vek 0$ means that
$\Vek w$ is positive semi-definite.  Although the positive
semi-definiteness of the lifted quadratic problem is not mandatory, it
strongly reduces the space of possible solutions.  The central benefit
of these reformulations is the shift of the non-linearity of the
forward operator into the rank-one constraint.  Although the problem
is now linear, we have to deal with an additional non-convex side
condition.

In order to eliminate the non-linear constraint for bilinear
operators, we first rewrite \eqref{eq:lift-bilin-prob} into the rank
minimization problem
\begin{equation*}
  \minimize
  \quad
  \rank(\Vek w)
  \quad\text{with}\quad
  \breve{\mathcal B}(\Vek w) = \Vek g^\dagger
\end{equation*}
and then relax the non-convex objective function by replacing it with
the \emph{nuclear} or \emph{projective norm}
$\pNormn{\cdot}_{\mathcal H_1 \otimes_\uppi \mathcal H_2}$ of the
tensor $\Vek w$.  Depending on the \PHilbert\ spaces $\mathcal H_1$
and $\mathcal H_2$, this norm is defined by
\begin{equation*}
  \pNormn{\Vek w}_{\mathcal H_1 \otimes_\uppi \mathcal H_2}
  \coloneqq
  \inf \biggl\{ \sum_{n=1}^N 
  \pNormn{\Vek u_n}_{\mathcal H_1} \, \pNormn{\Vek v_n}_{\mathcal H_2} :
  \Vek w = \sum_{n=1}^N \Vek u_n \otimes \Vek v_n, N \in \N \biggr\},
\end{equation*}
where the infimum is taken over all finite representations of the
tensor $\Vek w$.  In so doing, we finally obtain the convex
minimization problem
\begin{equation}
  \label{eq:nucl-min-bilin}
  \minimize
  \quad
  \pNormn{\Vek w}_{\mathcal H_1 \otimes_\uppi \mathcal H_2}
  \quad\text{subject to}\quad
  \breve{\mathcal B}(\Vek w) = \Vek g^\dagger,
  \tag{$\mathfrak B_0$}
  \addmathskip
\end{equation}
i.e., with linear constraints.

Since $\mathcal H_1$ as well as $\mathcal H_2$ is a \PHilbert\ space,
the nuclear norm here coincides with the trace class norm or with the
\PSchatten\ one-norm; so the nuclear norm is simply the sum of the
corresponding singular values of the matrix
$\Vek w \in \R^{N_2 \times N_1}$ with respect to the \PHilbert\ spaces
$\mathcal H_1$ and $\mathcal H_2$, \cf\ for instance
\cite[Satz~VI.5.5]{Wer02}.

\begin{Lemma}[Projective norm, \cite{Wer02}]
  For $\Vek w \in \mathcal H_1 \otimes \mathcal H_2$, the projective
  norm is given by
  $\pNormn{\Vek w}_{\mathcal H_1 \otimes_\uppi \mathcal H_2} =
  \sum_{n=0}^{R-1} \sigma_n$,
  where $\sigma_n$ denotes the $n$-th singular value and $R$ the rank of
  $\Vek w$.
\end{Lemma}

The main idea behind the nuclear norm heuristic is that the projective
norm is the convex envelope of the rank on the unit ball with respect
to the spectral norm.  Usually, the nuclear norm heuristic empirically
yields low-rank solutions of the matrix equation
$\breve{\mathcal B}(\Vek w) = \Vek g^\dagger$, see for example
\cite{CCS10,MGC11,RFP10} and references therein.  If the lifted
operator $\breve{\mathcal B}$ fulfills an appropriate restricted
isometry property, one can rigorously prove that the bilinear inverse
problem \eqref{eq:gen-bilin-prob} and the relaxed nuclear norm
minimization problem \eqref{eq:nucl-min-bilin} are equivalent, see
\cite[Theorem~3.3]{RFP10}.

For quadratic inverse problems, we follow a similar approach.  This
means that we first rewrite \eqref{eq:gen-quad-prob} into a rank
minimization problem over the positive semi-definite (symmetric)
matrix cone.  In order to obtain a convex minimization problem, we
again replace the objective function by the projective norm.  Since
the singular values $\sigma_n$ of a symmetric tensor $\Vek w$ coincide
with the absolute value of its eigenvalues $\lambda_n$, the positive
semi-definiteness may be incorporated into the projective norm by
\begin{equation*}
  \pNormn{\Vek w}_{\mathcal H \otimes_\uppi \mathcal H}^+
  \coloneqq
  \sum_{n=0}^{R-1} \lambda_n + \chi_{[0,\infty)}(\lambda_n),
\end{equation*}
where the indicator function $\chi_{[0, \infty)}$ is equal to $0$ for
arguments in $[0, \infty)$ and $+\infty$ otherwise; so the modified
projective ``norm''
$\pNormn{\cdot}_{\mathcal H \otimes_\uppi \mathcal H}^+$ sums up the
non-negative eigenvalues for positive semi-definite tensors and is
infinity otherwise.  To solve the quadratic inverse problem
\eqref{eq:gen-quad-prob} numerically, we thus consider the convex
minimization problem
\begin{equation}
  \label{eq:nucl-min-quad}
  \minimize \quad
  \pNormn{\Vek w}_{\mathcal H \otimes_\uppi \mathcal H}^+
  \quad\text{subject to}\quad
  \breve{\mathcal Q}(\Vek w) = \Vek g^\dagger.
  \tag{$\mathfrak Q_0$}
\end{equation}

Up to this point, the given data $\Vek g^\dagger$ have been known
exactly.  A first approach to incorporate noisy measurements into the
inverse problems \eqref{eq:gen-bilin-prob} and
\eqref{eq:gen-quad-prob} could be to extend the subspace of exact
solutions with respect to a supposed error level.  More precisely, one
may consider the minimization problems
\begin{alignat}{2}
  \label{eq:inexact-data-bilin}
  &\minimize \quad
  \pNormn{\Vek w}_{\mathcal H_1 \otimes_\uppi \mathcal H_2}
  &&\quad\text{subject to}\quad
  \pNormn{\breve{\mathcal B}(\Vek w)
  - \Vek g^\epsilon}_{\mathcal K} \le \epsilon
  \tag{$\mathfrak B_{\epsilon}$}
  \\
\shortintertext{and}
  \label{eq:inexact-data-quad}
  &\minimize \quad
  \pNormn{\Vek w}_{\mathcal H \otimes_\uppi \mathcal H}^+
  &&\quad\text{subject to}\quad
  \pNormn{\breve{\mathcal Q}(\Vek w)
  - \Vek g^\epsilon}_{\mathcal K} \le \epsilon.
  \tag{$\mathfrak Q_{\epsilon}$}
\end{alignat}
In other words, we minimize over all solutions that approximate the
given noisy data $\Vek g^\epsilon$ with
$\pNormn{\Vek g - \Vek g^\epsilon}_{\mathcal K} \le \epsilon$ up to
the error level $\epsilon$.

Another approach is to incorporate the data fidelity of the possible
solutions directly into the objective function.  Following this
approach, we may minimize a \PTikhonov\ functional with projective
norm regularization to solve \eqref{eq:gen-bilin-prob} and
\eqref{eq:gen-quad-prob}.  In more detail, we consider the problems 
\begin{align}
  \label{eq:tikh-func-bilin}
  &\minimize \quad
  \tfrac12 \, \pNormn{\breve{\mathcal B}(\Vek w) - \Vek g^\epsilon}_{\mathcal K}^2
  + \alpha \, \pNormn{\Vek w}_{\mathcal H_1 \otimes_\uppi \mathcal
    H_2}
  \tag{$\mathfrak B_\alpha$}
    \shortintertext{and}
  \label{eq:tikh-func-quad}
  &\minimize \quad
  \tfrac12 \, \pNormn{\breve{\mathcal Q}(\Vek w) - \Vek g^\epsilon}_{\mathcal K}^2
  + \alpha \, \pNormn{\Vek w}^+_{\mathcal H \otimes_\uppi \mathcal
    H}.
  \tag{$\mathfrak Q_\alpha$}
\end{align}

All proposed convex relaxations of the lifted problems
\eqref{eq:lift-bilin-prob} and \eqref{eq:lift-quad-prob} have in
common that the minimization of the projective norm usually promotes
low-rank tensors and, in the best case, yields a rank-one solutions.
The later case forms our basic premise to derive numerical solutions.

\begin{Assumption}[Basic premise]
  \label{ass:basic-pre}
  Suppose that the solutions $\Vek w$ of the relaxation
  \eqref{eq:nucl-min-bilin}, \eqref{eq:inexact-data-bilin}, or
  \eqref{eq:tikh-func-bilin} are at most rank-one tensors
  $\Vek u \otimes \Vek v$ such that $(\Vek u, \Vek v)$ (approximately)
  solves the bilinear inverse problem \eqref{eq:gen-bilin-prob}.
  Likewise, suppose that the solutions $\Vek w$ of the relaxations
  \eqref{eq:nucl-min-quad}, \eqref{eq:inexact-data-quad}, or
  \eqref{eq:tikh-func-quad} have at most rank one so that they
  (approximately) solve the quadratic inverse problem
  \eqref{eq:gen-quad-prob}.
\end{Assumption}

\section{Proximal algorithms for the lifted problem}
\label{sec:sing-value-thresh}

To exploit the nuclear norm heuristic, we have to solve the
minimization problem \eqref{eq:nucl-min-bilin},
\eqref{eq:inexact-data-bilin}, \eqref{eq:tikh-func-bilin},
\eqref{eq:nucl-min-quad}, \eqref{eq:inexact-data-quad}, and
\eqref{eq:tikh-func-quad} in an efficient manner.  Looking back at the
comprehensive literature about matrix completion
\cite{CR09,CCS10,MGC11}, low-rank solutions of matrix equations
\cite{RFP10}, and PhaseLift \cite{CESV13,CSV13}, there exists several
numerical methods like interior-point methods for semi-definite
programming, fixed point iterations, singular value thresholding
algorithm, projected subgradient methods, and low-rank parametrization
approaches.

In order to solve the lifted bilinear and quadratic inverse problems,
we follow another approach.  Let us first consider the actual
structure of the six derived minimization problems in
\autoref{sec:bilin-inverse-probl}, which is given by
\begin{equation}
  \label{eq:gen-prob-struc}
  \minimize \quad
  F(\mathcal A (\Vek w)) + G(\Vek w),
\end{equation}
where
$\mathcal A \colon \mathcal H_1 \otimes \mathcal H_2 \to \mathcal K$
denotes the lifted bilinear or quadratic forward operator,
$F \colon \mathcal K \to \overline \R$ with
$\overline \R \coloneqq \R \cup \{-\infty, +\infty\}$ describes the
data fidelity, and
$G \colon \mathcal H_1 \times \mathcal H_2 \to \overline \R$ is the
(modified) projective norm.  Since the regularization mapping $G$ and,
in some circumstances, the data fidelity mapping $F$ are non-smooth
but convex functions, we may apply proximal first-order methods like
the forward-backward splitting, the primal-dual method by
\pers{Chambolle--Pock}, the alternating directions method of
multipliers (ADMM), the \pers{Douglas--Rachford} splitting, and
several variants of these and other algorithms, see for instance
\cite{CP16}.

Although we can apply any of these algorithms, we exemplarily consider
the forward-backward splitting \cite{LM79,CW05}
\begin{equation}
  \label{eq:forward-backward}
  \Vek w^{(n+1)} 
  \coloneqq
  \prox_{\tau G} \bigl(\Vek w^{(n)} 
  - \tau \mathcal A^* \nabla F(\mathcal A\Vek w^{(n)}) \bigr)
\end{equation}
 and the primal-dual method \cite[Algorithm~1]{CP11}  
\begin{equation}
  \label{eq:primal-dual}
  \begin{aligned}
    \Vek y^{(n+1)} &\coloneqq \prox_{\sigma F^*}
    \bigl(\Vek y^{(n)} + \sigma \, \mathcal A(\breve{\Vek w}^{(n)})\bigr)
    \\[\fskip]
    \Vek w^{(n+1)} &\coloneqq \prox_{\tau G} \bigl(\Vek
    w^{(n)} - \tau \, \mathcal A^* (\Vek y^{(n+1)})\bigr)
    \\[\fskip]
    \breve{\Vek w}^{(n+1)} &\coloneqq \Vek w^{(n+1)} + \theta \, (\Vek
    w^{(n+1)} - \Vek w^{(n)})
  \end{aligned}
\end{equation}
with fixed parameters $\tau, \sigma > 0$ and $\theta \in [0,1]$.  The
details of these methods are given below.

First, the forward-backward splitting and the primal-dual method are
originally defined for linear forward operators between
finite-dimensional \PHilbert\ spaces; so we have to equip the tensor
product $\mathcal H_1 \otimes \mathcal H_2$ with a corresponding
structure.  In the following, we always assume that this structure
arises from the inner product defined by
\begin{equation*}
  \iProdn{\Vek u_1 \otimes \Vek v_1}{\Vek u_2 \otimes \Vek
    v_2}_{\mathcal H_1 \otimes \mathcal H_2}
  \coloneqq 
  \iProdn{\Vek u_1}{\Vek u_2}_{\mathcal H_1} 
  \iProdn{\Vek v_1}{\Vek v_2}_{\mathcal H_2},
\end{equation*}
see for instance \cite[Section~2.6]{KR83}.  Using the matrices
$\Mat H_1$ and $\Mat H_2$ defining the inner products of
$\mathcal H_1$ and $\mathcal H_2$, we can write the resulting inner
product of the \PHilbert{ian} tensor product
$\mathcal H_1 \otimes \mathcal H_2$ in the form
\begin{equation}
  \label{eq:new-inner-prod}
  \iProdn{\Vek w_1}{\Vek w_2}_{\mathcal H_1 \otimes \mathcal
    H_2} 
  \coloneqq
  \iProdn{\Mat H_2 \Vek w_1 \Mat H_1}{\Vek w_2}
  =
  \iProdn{\Vek w_1}{\Mat H_2 \Vek w_2 \Mat H_1}
  =
  \tr(\Vek w_2^* \Mat H_2 \Vek w_1 \Mat H_1),
\end{equation}
where the inner products on the right-hand side denote the
\PHilbert--\PSchmidt\ inner product for matrices.  For the quadratic
setting, the symmetric \PHilbert{ian} tensor product $\mathcal H
\otimes_{\mathrm{sym}} \mathcal H$ is the related subspace
embracing the symmetric tensors.

Next, the above stated methods are mainly based on concepts from
convex analysis, which is reflected in the presuppositions; so the
data fidelity mapping $F$ and, similarly, the regularization mapping
$G$ have to be \emph{convex} and \emph{lower semi-continuous}.
Commonly, a function $f \colon \mathcal X \to \overline \R$ on the
real \PHilbert\ space $\mathcal X$ 
is called \emph{convex} when
\begin{equation*}
  f(t \Vek x_1 + (1-t) \, \Vek x_2) 
  \le
  t \, f(\Vek x_1) + (1-t) \, f(\Vek x_2)
  \addmathskip
\end{equation*}
for all $\Vek x_1, \Vek x_2 \in \mathcal X$ and all $t \in [0,1]$, and
\emph{lower semi-continuous} when
\begin{equation*}
  f(\Vek x) 
  \le
  \liminf_{n\to\infty} f(\Vek x_n)
\end{equation*}
for all sequences $(\Vek x_n)$ in $\mathcal X$ with
$\Vek x_n \to \Vek x$.  Since $F$ and $G$ in the relaxed minimization
problems of \autoref{sec:bilin-inverse-probl} represent norms or
indicator functions on closed convex sets, here the assumptions for
the primal-dual method are
always fulfilled.  The forward-backward splitting additionally
requires a differentiable data fidelity $F$ with
\PLipschitz-continuous derivative; so this method can only be applied
to the \PTikhonov\ relaxations.

For the primal-dual iteration \eqref{eq:primal-dual}, the first
proximity operator $\prox_{\sigma F^*}$ is computed with respect to the
\PLegendre--\PFenchel\ conjugate $F^*$.  For any function $f \colon
\mathcal X \to \overline \R$, where $\mathcal X$ denotes an arbitrary
\PHilbert\ space, the \emph{\PLegendre--\PFenchel\ conjugate} $f^*
\colon \mathcal X \to \overline \R$ is defined by
\begin{equation*}
  f^*(\Vek x') \coloneqq
  \sup_{\Vek x \in \mathcal X} \iProdn{\Vek x'}{\Vek x}_{\mathcal X} - f(\Vek x)
\end{equation*}
and is always convex and lower semi-continuous, see \cite{Roc70}.  If
the function $f \colon \mathcal X \to \overline \R$ is lower
semi-continuous and convex, the \emph{subdifferential} $\uppartial f$
at a certain point $\Vek x \in \mathcal X$ is given by
\begin{equation*}
  \uppartial f (\Vek x)
  \coloneqq
  \bigl\{ \Vek x' \in \mathcal X
  :
  f(\Vek y) \ge f(\Vek x) + \iProdn{\Vek x'}{ \Vek y - \Vek
    x}_{\mathcal X}
  \;\text{for all}\;
  \Vek y \in \mathcal X
  \bigr\}
\end{equation*}
and figuratively consists of all linear minorants, see \cite{Roc70}.
Finally, the \emph{proximation} $\prox_{f}$ of a lower
semi-continuous, convex function
$f \colon \mathcal X \to \overline \R$ is defined as the unique
minimizer 
\begin{equation*}
  \prox_{f}(\Vek x)
  \coloneqq
  \argmin_{\Vek y \in \mathcal X} f(\Vek x) + \tfrac{1}{2} \,
  \pNormn{\Vek y - \Vek x}_{\mathcal X}^2.
  \addmathskip
\end{equation*}
Using the subdifferential calculus, one can show that the proximation
coincides with the \emph{resolvent}, \ie\
\begin{equation*}
  \prox_{f}(\Vek x) = (I + \uppartial f)^{-1}(\Vek x),
  \addmathskip
\end{equation*}
see for instance \cite{Roc70, CP16}.

The most crucial step in the forward-backward splitting
\eqref{eq:forward-backward} and the primal-dual iteration
\eqref{eq:primal-dual} is the application of the proximal projective
norm
$\prox_{\tau \pNormn{\cdot}_{\mathcal H_1 \otimes_\uppi \mathcal
    H_2}}$,
whereas the computation of proximal conjugated data fidelity
$\prox_{\sigma F^*}$ is usually much simpler.  To determine the
proximal projective norm explicitly, we exploit the singular value
decomposition of the argument with respect to the underlying
\PHilbert\ spaces $\mathcal H_1$ and $\mathcal H_2$, which can be
derived by an adaption of the classical singular value decomposition
for matrices with respect to the \PEuclid{ian} inner product.

\begin{Lemma}[Singular value decomposition]
  \label{lem:svd}
  Let $\Vek w$ be a tensor in $\mathcal H_1 \otimes \mathcal H_2$.
  The singular value decomposition of $\Vek w$ with respect to
  $\mathcal H_1$ and $\mathcal H_2$ is given by
  \begin{equation*}
    \Vek w 
    = 
    \sum_{n=0}^{R-1} \sigma_n \,
    (\widetilde{\Vek u}_n \otimes \widetilde{\Vek v}_n)
    \qquad\text{with}\qquad
    \widetilde{\Vek u}_n 
    \coloneqq \Mat H_1^{-\nicefrac12} \Vek u_n
    \quad\text{and}\quad
    \widetilde{\Vek v}_n 
    \coloneqq \Mat H_2^{-\nicefrac12} \Vek v_n,
  \end{equation*}
  where $\sum_{n=0}^{R-1} \sigma_n \, (\Vek u_n \otimes \Vek v_n)$ is
  the classical singular value decomposition of
  $\Mat H_2^{\nicefrac12} \Vek w \, (\Mat H_1^{\nicefrac12})^*$ with
  respect to the \PEuclid{ian} inner product. 
\end{Lemma}

\begin{Remark}
  \label{rem:svd:1}
  Unless stated otherwise, the square roots
  $\Mat H_1^{\nicefrac12} \in \R^{N_1 \times N_1}$ and
  $\Mat H_2^{\nicefrac12} \in \R^{N_2 \times N_2}$ are taken with
  respect to the factorizations
  \begin{equation*}
    (\Mat H_1^{\nicefrac12})^* \, \Mat H_1^{\nicefrac12} = \Mat H_1
    \qquad\text{and}\qquad
    (\Mat H_2^{\nicefrac12})^* \, \Mat H_2^{\nicefrac12} = \Mat H_2.
  \end{equation*}
  Allowing also non-symmetric but invertible factorizations, the root
  $\Mat H_1^{\nicefrac12}$ and $\Mat H_2^{\nicefrac12}$ are here
  non-unique.  Possible candidates are the symmetric positive definite
  square root or the \PCholesky\ decomposition of $\Mat H_1$ and
  $\Mat H_2$.  In the following, the roots are solely required to
  derive the proximal projective norm mathematically.  Their actual
  computation is not necessary in the final tensor-free algorithm.
  \qed
\end{Remark}

\begin{Proof}[\thref{lem:svd}]
  By assumption the possibly non-symmetric square roots
  $\Mat H_1^{\nicefrac12}$ and $\Mat H_2^{\nicefrac12}$ are
  invertible.  Considering the classical \PEuclid{ian} singular value
  decomposition of the matrix
  $\Mat H_2^{\nicefrac12} \Vek w \, (\Mat H_1^{\nicefrac12})^*$, we
  immediately obtain
  \begin{equation*}
    \Vek w
    =
    \Mat H_2^{-\nicefrac12} \, 
    \biggl( \sum_{n=0}^{R-1} \sigma_n \, (\Vek u_n \otimes \Vek v_n)
    \biggr) \, (\Mat H_1^{-\nicefrac12})^*
    =
    \sum_{n=0}^{R-1} \sigma_n \, \bigl((\Mat H_1^{-\nicefrac12} \Vek u_n )
    \otimes (\Mat H_2^{-\nicefrac12} \Vek v_n)\bigr).
  \end{equation*}
  The last arrangement may be easily validated by using the matrix
  notation $\Vek v_n \Vek u_n^*$ of the rank-one tensor
  $\Vek u_n \otimes \Vek v_n$.  Due to the identity 
  \begin{equation*}
    \iProdn{\widetilde{\Vek u}_n}{\widetilde{\Vek u}_m}_{\mathcal H_1}
    =
    \iProdn{\Vek u_n}{(\Mat H_1^{-\nicefrac12})^* \Mat H_1 \Mat
      H_1^{-\nicefrac12} \Vek u_m}
    =
    \iProdn{\Vek u_n}{\Vek u_m}
  \end{equation*}
  for all $n,m \in \{0, \dots, R-1\}$, the singular vectors
  $\{\widetilde{\Vek u}_n : n=0,\dots,R-1\}$ form an orthonormal
  system with respect to $\mathcal H_1$ as well as their counterparts
  $\{\widetilde{\Vek v}_n : n=0,\dots, R-1\}$ with respect to
  $\mathcal H_2$.  \qed
\end{Proof}

\begin{Remark}
  \label{rem:svd:2}
  The singular value decomposition can also be interpreted as a matrix
  factorization of the tensor $\Vek w$.  In this case, we have the
  factorization
  \begin{equation*}
    \Vek w 
    =
    \widetilde{\Mat V} \Mat \Sigma \widetilde{\Mat U}^*
    \qquad\text{with}\qquad
    \widetilde{\Mat U} 
    \coloneqq \Mat H_1^{-\nicefrac12} \Mat U
    \quad\text{and}\quad
    \widetilde{\Mat V}
    \coloneqq \Mat H_2^{-\nicefrac12} \Mat V,
  \end{equation*}
  where $\Mat V \Mat \Sigma \Mat U^*$ is the classical \PEuclid{ian}
  singular value decomposition of
  $\Mat H_2^{\nicefrac12} \Vek w \, (\Mat H_1^{\nicefrac12})^*$ with
  the left singular vectors
  $\Mat V \coloneqq [\Vek v_0, \dots, \Vek v_{R-1}]$, the right
  singular vectors $\Mat U \coloneqq [\Vek u_0, \dots, \Vek u_{R-1}]$,
  and the singular values
  $\Mat \Sigma \coloneqq \diag(\sigma_0, \dots, \sigma_{R-1})$.  \qed
\end{Remark}

For the quadratic setting, we will rely on the eigenvalue
decomposition instead of the singular value decomposition.

\begin{Corollary}[Eigenvalue decomposition]
  \label{cor:evd}
  Let $\Vek w$ be a tensor in
  $\mathcal H \otimes_{\mathrm{sym}} \mathcal H$.  The eigenvalue
  decomposition of $\Vek w$ with respect to $\mathcal H$ is given by
  \begin{equation*}
    \Vek w 
    = 
    \sum_{n=0}^{R-1} \lambda_n \,
    (\widetilde{\Vek u}_n \otimes \widetilde{\Vek u}_n)
    \qquad\text{with}\qquad
    \widetilde{\Vek u}_n 
    \coloneqq \Mat H^{-\nicefrac12} \Vek u_n,
  \end{equation*}
  where $\sum_{n=0}^{R-1} \lambda_n \, (\Vek u_n \otimes \Vek v_n)$ is
  the classical eigenvalue decomposition of
  $\Mat H^{\nicefrac12} \Vek w \, (\Mat H^{\nicefrac12})^*$ with
  respect to the \PEuclid{ian} inner product. 
\end{Corollary}

\begin{Proof}
  Due to the close relation between singular value decomposition and
  eigenvalue decomposition for symmetric matrices, the assertion
  follows from \thref{lem:svd} with
  $\lambda_n \coloneqq \sigma_n \, \iProdn{\widetilde{\Vek
      u}_n}{\widetilde{\Vek v}_n}_{\mathcal H}$.  \qed
\end{Proof}

With the adaption in \thref{lem:svd}, we can apply any numerical
singular value method to compute the singular value decomposition of a
given tensor.  The next ingredient is the subdifferential of the
nuclear norm based on $\mathcal H_1$ and $\mathcal H_2$ with respect
to $\mathcal H_1 \otimes \mathcal H_2$.  In the
following, the set-valued signum function $\sgn$ is defined by
\begin{equation*}
  \sgn(t) \coloneqq 
  \begin{cases}
    \{ 1 \} & \text{if} \; t > 0, \\
    [-1,1] & \text{if} \; t = 0, \\
    \{ -1 \} & \text{if} \; t < 0.
  \end{cases}
\end{equation*}

\begin{Lemma}[Subdifferential]
  \label{lem:subdiff}
  Let $\Vek w$ be a tensor in
  $\mathcal H_1 \otimes \mathcal H_2$.
  Then the subdifferential of the projective norm
  $\pNormn{\cdot}_{\mathcal H_1 \otimes_\uppi \mathcal H_2}$ at
  $\Vek w$ with respect to
  $\mathcal H_1 \otimes \mathcal H_2$ is given by
  \begin{equation*}
    \uppartial
    \pNormn{\cdot}_{\mathcal H_1 \otimes_\uppi
      \mathcal H_2} (\Vek w)
    =
    \biggl\{
    \sum_{n=0}^{R-1} \mu_n \, (\widetilde{\Vek u}_n \otimes
    \widetilde{\Vek v}_n) : 
    \mu_n \in \sgn( \sigma_n),
    \Vek w = \sum_{n=0}^{R-1} \sigma_n \, 
    (\widetilde{\Vek u}_n \otimes \widetilde{\Vek v}_n)
    \biggr\}, 
  \end{equation*}
  where
  $\Vek w = \sum_{n=0}^{R-1} \sigma_n \, (\widetilde{\Vek u}_n \otimes
  \widetilde{\Vek v}_n)$
  is a valid singular value decomposition of $\Vek w$ with respect to
  $\mathcal H_1$ and $\mathcal H_2$.   
\end{Lemma}

\begin{Proof}
  The central idea to compute the subdifferential is to rely on the
  corresponding statement for the \PEuclid{ian} setting in
  \cite{Lew95}.  More precisely, if $\mathcal H_1 = \R^{N_1}$ and
  $\mathcal H_2 = \R^{N_2}$ are equipped with the \PEuclid{ian} inner
  product, then the subdifferential
  $\uppartial_{\mathcal{H\!S}}$ with respect to the
  \PHilbert--\PSchmidt\ inner product on $\R^{N_1} \otimes \R^{N_2}$
  is given by
  \begin{equation*}
    \uppartial_{\mathcal{H\!S}} \pNormn{\cdot}_{\R^{N_1} \otimes_\uppi
      \R^{N_2}} (\Vek w)
    =
    \biggl\{
    \sum_{n=0}^{R-1} \mu_n \, (\Vek u_n \otimes
    \Vek v_n) : 
    \mu_n \in \sgn( \sigma_n),
    \Vek w = \sum_{n=0}^{R-1} \sigma_n \, 
    (\Vek u_n \otimes \Vek v_n)
    \biggr\}, 
  \end{equation*}
  where
  $\Vek w = \sum_{n=0}^{R-1} \sigma_n \, (\Vek u_n \otimes \Vek v_n)$
  is an \PEuclid{ian} singular value decomposition of $\Vek w$, see
  \cite[Corollary~2.5]{Lew95}.  The upper bound $R$ is here some
  number less than or equal to $\min \{N_1, N_2\}$, and the singular
  value decomposition of $\Vek w$ may contain zero as singular value.

  Next, we adapt this result to our specific case.  Therefore, we
  exploit that the projective norm is the sum of the singular
  values.  Using \thref{lem:svd}, we notice that the generalized
  projective norm of a tensor $\Vek w$ is given by
  \begin{equation}
    \label{eq:rescal-nukl-norm}
    \pNormn{\Vek w}_{\mathcal H_1 \otimes_\uppi \mathcal H_2}
    = \pNormn{\Mat H_2^{\nicefrac12} \Vek w
      \, (\Mat H_1^{\nicefrac12})^*}_{\R^{N_1} \otimes_\uppi \R^{N_2}},
  \end{equation}
  where the norm on the right-hand side is the usual projective norm
  with respect to the \PEuclid{ian} inner product.  In order to
  consider the inner product of
  $\mathcal H_1 \otimes \mathcal H_2$ in the
  subdifferential, we exploit that
  \begin{equation*}
    \pNormn{\breve{\Vek w}}_{\mathcal H_1 \otimes_\uppi \mathcal  H_2} 
    \ge 
    \pNormn{\Vek w}_{\mathcal H_1 \otimes_\uppi \mathcal  H_2}
    + \iProd{\Vek \xi}{\breve{\Vek w} - \Vek w}_{\mathcal H_1
      \otimes \mathcal H_2}
    \qquad
    \text{for all}
    \qquad
    \breve{\Vek w} \in \mathcal H_1
    \otimes \mathcal H_2
  \end{equation*}
  if and only if
  \begin{equation*}
    \pNormn{\breve{\Vek w}}_{\mathcal H_1 \otimes_\uppi \mathcal  H_2} 
    \ge 
    \pNormn{\Vek w}_{\mathcal H_1 \otimes_\uppi \mathcal  H_2}
    + \iProd{\Mat H_2 \Vek \xi \Mat H_1}{\breve{\Vek w} - \Vek w}
    \qquad
    \text{for all}
    \qquad
    \breve{\Vek w} \in \R^{N_2 \times N_1},
  \end{equation*}
  where the inner product on the right-hand side is the usual
  \PHilbert--\PSchmidt\ scalar product for matrices, see
  \eqref{eq:new-inner-prod}.  Thus, the subdifferential with respect
  to the $\mathcal H_1 \otimes \mathcal H_2$ scalar product can be
  expressed in terms of $\uppartial_{\mathcal{H\!S}}$ by
  \begin{equation}
    \label{eq:rescal-subdiff}
    \uppartial
    \pNormn{\cdot}_{\mathcal H_1 \otimes_\uppi \mathcal H_2} (\Vek w)
    = \Mat H_2^{-1} \, \uppartial_{\mathcal{H\!S}}
    \pNormn{\cdot}_{\mathcal H_1 \otimes_\uppi \mathcal H_2} (\Vek w) \,
    \Mat H_1^{-1}.
  \end{equation}
  Plugging \eqref{eq:rescal-nukl-norm} into \eqref{eq:rescal-subdiff},
  and using the chain rule, we obtain the assertion.  \qed
\end{Proof}

With the characterization of the subdifferential, we are ready to
determine the proximity operator for the projective norm with respect
to the underlying \PHilbert\ spaces $\mathcal H_1$ and $\mathcal H_2$.
In the following, the soft-thresholding operator with respect to the
level $\tau$ is defined by
\begin{equation*}
  S_\tau(t)
  \coloneqq
  \begin{cases}
    t - \tau
    & \text{if} \; t > \tau,
    \\
    t + \tau
    & \text{if} \; t < - \tau,
    \\
    0
    & \text{otherwise.}
  \end{cases}
\end{equation*}

\begin{Theorem}[Proximal projective norm]
  \label{the:resolv-op}
  Let $\Vek w$ be a tensor  in $\mathcal H_1 \otimes
  \mathcal H_2$.  The proximation of the projective norm is given by
  \begin{equation*}
    \prox_{\tau \pNorm{\cdot}_{\mathcal H_1 \otimes_\uppi \mathcal
      H_2}}(\Vek w)
    =
    \sum_{n=0}^{R-1} S_\tau (\sigma_n) \, (\widetilde{\Vek u}_n \otimes
    \widetilde{\Vek v}_n), 
  \end{equation*}
  where
  $\sum_{n=0}^{R-1} \sigma_n \, (\widetilde{\Vek u}_n \otimes
  \widetilde{\Vek v}_n)$
  is a singular value decomposition of $\Vek w$ with respect to the
  \PHilbert\ spaces $\mathcal H_1$ and $\mathcal H_2$.
\end{Theorem}

  \begin{Proof}
    In order to establish the statement, we only have to convince
    ourselves that
    $\breve{\Vek w} \coloneqq \sum_{n=0}^{R-1} S_\tau (\sigma_n) \,
    (\widetilde{\Vek u}_n \otimes \widetilde{\Vek v}_n)$ is the
    resolvent for the given $\Vek w$, \ie\
    $\Vek w \in (I + \tau \, \uppartial \pNormn{\cdot}_{\mathcal H_1
      \otimes_\uppi \mathcal H_2})(\breve{\Vek w})$.  Since
    $\breve{\Vek w}$ is already represented by its singular value
    decomposition, \thref{lem:subdiff} implies
    \begin{equation*}
      \sum_{n=0}^{R-1} \bigl[S_\tau(\sigma_n) + \tau
      \sgn(S_\tau(\sigma_n))\bigr] \, (\widetilde{\Vek u}_n \otimes
      \widetilde{\Vek v}_n)
      \subset
      (I + \tau \, \uppartial \pNormn{\cdot}_{\mathcal H_1
        \otimes_\uppi \mathcal H_2})(\breve{\Vek w}).
    \end{equation*}
    If $\sigma_n > \tau$, the related summand becomes
    $\sigma_n \, (\widetilde{\Vek u}_n \otimes \widetilde{\Vek v}_n)$.
    Otherwise, the summand is
    $\mu_n \, (\widetilde{\Vek u}_n \otimes \widetilde{\Vek v}_n)$
    with $\mu \in [-\tau, \tau]$.  Since the singular value
    decomposition of $\Vek w$ obviously has this form, the proof is
    completed.  \qed
  \end{Proof}

\begin{Remark}[Singular value thresholding]
  \label{rem:resolv-op}
  The proximation of the projective norm with respect to
  $\mathcal H_1$ and $\mathcal H_2$ is simply a soft thresholding of
  the corresponding singular values.  In the following, we denote the
  matrix-valued operation
  \begin{equation*}
    \Vek w 
    \mapsto
    \sum_{n=0}^{R-1} S_\tau(\sigma_n) \, (\widetilde{\Vek u}_n \otimes
    \widetilde{\Vek v}_n)
  \end{equation*}
  as \emph{(soft) singular value thresholding} $\mathcal S_\tau$.  \qed 
\end{Remark}

The proximation
$\prox_{\tau \pNormn{\cdot}_{\mathcal H \otimes_\uppi \mathcal H}^+}$
associated to the modified projective norm used in the quadratic
setting may be computed analogously.  As preliminary step, we
determine the subdifferential of the modified norm with respect to the
symmetric \PHilbert{ian} tensor product.  In a nutshell, we have
simply to replace the set-valued signum function by the modified
signum $\sgn^+$ defined by
\begin{equation*}
  \sgn^+(t) \coloneqq
  \begin{cases}
    \{ 1\} & \text{if } t > 0,\\
    (-\infty,1] & \text{if } t = 0,\\
    \emptyset & \text{if } t < 0.
  \end{cases}
\end{equation*}
Since the truncated absolute value $\absn{\cdot} + \chi_{[0,\infty)}$
is not an absolutely symmetric function, one cannot apply the
subdifferential characterization in \cite{Lew96}; there we compute the
subdifferential directly.

\begin{Lemma}[Subdifferential]
  \label{lem:subdiff-plus}
  Let $\Vek w$ be a positive semi-definite tensor in
  $\mathcal H \otimes_{\mathrm{sym}} \mathcal H$.  The subdifferential
  of the modified projective norm
  $\pNormn{\cdot}^+_{\mathcal H \otimes_\uppi \mathcal H}$ at $\Vek w$
  with respect to
  $\mathcal H \otimes_{\mathrm{sym}} \mathcal H$ is given
  by
  \begin{equation*}
    \uppartial
    \pNormn{\cdot}^+_{\mathcal H \otimes_\uppi \mathcal H}
    (\Vek w)
    =
    \biggl\{ \sum_{n=0}^{R_1} \mu_n \, (\widetilde{\Vek u}_n \otimes
    \widetilde{\Vek u}_n) :
    \mu_n \in \sgn^+(\lambda_n), 
    \Vek w = \sum_{n=0}^{R-1} \lambda_n \, (\widetilde{\Vek u}_n \otimes
    \widetilde{\Vek u}_n) \biggr\},
  \end{equation*}
  where
  $\Vek w = \sum_{n=0}^{R-1} \lambda_n \, (\widetilde{\Vek u}_n
  \otimes \widetilde{\Vek u}_n)$
  is an eigenvalue decomposition of $\Vek w$ with respect to
  $\mathcal H$.  If $\Vek w$ is not positive semi-definite, then the
  subdifferential is empty.
\end{Lemma}

\begin{Proof}
  Similarly to \thref{lem:subdiff}, we rely on the corresponding
  statement for the \PEuclid{ian} setting in \cite{Lew99}.  If we
  endow $\mathcal H \simeq \R^N$ with the \PEuclid{ian} inner product,
  then the subdifferential of the modified projective norm is
  determined by
  \begin{equation*}
    \uppartial
    \pNormn{\cdot}^+_{\R^N \otimes_\uppi \R^N}
    (\Vek w)
    =
    \biggl\{ \sum_{n=0}^{R_1} \mu_n \, ({\Vek u}_n \otimes
    {\Vek u}_n) :
    \mu_n \in \sgn^+(\lambda_n), 
    \Vek w = \sum_{n=0}^{R-1} \lambda_n \, ({\Vek u}_n \otimes
    {\Vek u}_n) \biggr\},
  \end{equation*}
  where
  $\Vek w = \sum_{n=0}^{R-1} \lambda_n \, ({\Vek u}_n \otimes {\Vek
    u}_n)$
  is a eigenvalue decomposition of $\Vek w$, see
  \cite[Theorem~6]{Lew99}.  The transformation to an arbitrary
  \PHilbert\ space $\mathcal H$ works exactly as in the proof of
  \thref{lem:subdiff}.  \qed
\end{Proof}

Together with the positive soft-thresholding operator $S^+_\tau$
defined by
\begin{equation*}
  S_\tau^+(t) \coloneqq
  \begin{cases}
    t - \tau & \text{if } t > \tau, \\
    0 & \text{otherwise},
  \end{cases}
\end{equation*}
the computed subdifferential leads us to the following proximity
operator.

\begin{Theorem}[Proximal projective norm]
  \label{the:mod-resolv-oper}
  Let $\Vek w$ be a tensor in $\mathcal H \otimes \mathcal H$.  Then
  the proximation of the modified projective norm is given by
  \begin{equation*}
    \prox_{\tau \pNormn{\cdot}^+_{\mathcal H \otimes_\uppi \mathcal H}}
    (\Vek w)
    =
    \sum_{n=0}^{R-1} S_\tau^+ (\lambda_n) \, (\widetilde{\Vek u}_n
    \otimes \widetilde{\Vek u}_n),
  \end{equation*}
  where $\sum_{n=0}^{R-1} \lambda_n \, (\widetilde{\Vek u}_n \otimes
  \widetilde{\Vek u}_n)$ is an eigenvalue decomposition of $\Vek w$
  with respect to $\mathcal H$.
\end{Theorem}

\begin{Proof}
  The assertion follows similarly to \thref{the:resolv-op} by
  replacing the subdifferential in \thref{lem:subdiff} by
  \thref{lem:subdiff-plus} and the singular value decomposition by an
  eigenvalue decomposition.  \qed
\end{Proof}

\begin{Remark}
  \label{rem:mod-resolv-oper:1}
  Analogously to the tensor-valued singular value thresholding operator
  $\mathcal S_\tau$, we define the positive eigenvalue thresholding operator
  $\mathcal S_\tau^+$, which additionally projects the argument to the
  positive semi-definite cone, by
  \begin{equation*}
    \Vek w \mapsto \sum_{n=0}^{R-1} S_\tau^+(\lambda_n) \,
    (\widetilde{\Vek u}_n \otimes \widetilde{\Vek u}_n).
    \tag*{\qed}
  \end{equation*} 
\end{Remark}

Knowing the proximation of the (modified) projective norm, we are now
able to perform proximal algorithms to solve the
minimization problem in \autoref{sec:bilin-inverse-probl}.  Although
we can use any of the mentioned method, here we restrict ourselves the
primal-dual iteration \eqref{eq:primal-dual}.  First, we consider the
bilinear minimization problems \eqref{eq:nucl-min-bilin},
\eqref{eq:inexact-data-bilin}, and \eqref{eq:tikh-func-bilin}.

For exactly given data $\Vek g^\dagger$ corresponding to the
minimization problem \eqref{eq:nucl-min-bilin}, the data fidelity
functional corresponds to $F \colon \mathcal K \to \overline\R$ with
$F(\Vek y) \coloneqq \chi_{\{0\}}(\Vek y - \Vek g^\dagger)$.  A simple
computation shows that the conjugate $F^*$ is given by $F^*(\Vek y) \coloneqq
\iProdn{\Vek y}{\Vek g^\dagger}_{\mathcal K}$ and the associate
proximal mapping by
\begin{equation*}
  \prox_{\sigma F^*}(\Vek y)
  =(I + \sigma \, \uppartial F^*)^{-1}(\Vek y)
  = \Vek y - \sigma \, \Vek g^\dagger.
\end{equation*}
Thus, we obtain the following algorithm.

\begin{Algorithm}[Primal-dual for exact data]
  \label{alg:exact-data}
  \begin{enumerate}[(i)]
  \item Initiation: Fix the parameters $\tau, \sigma > 0$ and
    $\theta \in [0,1]$.  Choose an arbitrary start value
    $(\Vek w^{(0)}, \Vek y^{(0)})$ in
    $(\mathcal H_1 \otimes \mathcal H_2) \times \mathcal
    K$, and set $\breve{\Vek w}^{(0)}$ to $\Vek w^{(0)}$.
  \item Iteration: For $n>0$, update $\Vek w^{(n)}$, $\breve{\Vek
      w}^{(n)}$, and $\Vek y^{(n)}$ by 
    \begin{align*}
      \Vek y^{(n+1)}
      &\coloneqq
        \Vek y^{(n)} + \sigma \, (\breve{\mathcal B}(\breve{\Vek
        w}^{(n)}) - \Vek g^\dagger) 
      \\[\fskip]
      \Vek w^{(n+1)}
      &\coloneqq
        \mathcal S_{\tau}\bigl(\Vek w^{(n)} - \tau \,
        \breve{\mathcal B}^* (\Vek 
        y^{(n+1)})\bigr)
      \\[\fskip]
      \breve{\Vek w}^{(n+1)}
      &\coloneqq
        \Vek w^{(n+1)} + \theta \, ( \Vek w^{(n+1)} - \Vek w^{(n)}).
    \end{align*}
  \end{enumerate}
\end{Algorithm}

\begin{Remark}
  \label{rem:exact-data}
  If the projective norm in \eqref{eq:nucl-min-bilin} is weighed with
  a parameter $\alpha > 0$ in order to control the influence of the
  data fidelity and the regularization, \cf\
  \eqref{eq:tikh-func-bilin}, then the iteration in
  \thref{alg:exact-data} changes slightly.  More precisely, one has to
  replace $\mathcal S_\tau$ by $\mathcal S_{\alpha \tau}$.  \qed
\end{Remark}

For inexact data $\Vek g^\epsilon$, we first consider the \PTikhonov\
minimization \eqref{eq:tikh-func-bilin}, whose data fidelity
corresponds to
$F(\Vek y) \coloneqq \nicefrac12 \, \pNormn{\Vek y - \Vek
  g^\epsilon}_{\mathcal K}^2$.
Here the conjugate $F^*$ is given by
$F^*(\Vek y) = \nicefrac12 \, \pNormn{\Vek y}_{\mathcal K}^2 +
\iProd{\Vek y}{\Vek g^\epsilon}_{\mathcal K}$
with subdifferential
$\uppartial F^*(\Vek y) = \{\Vek y + \Vek g^\epsilon\}$.  Again a
simple computation leads to the proximation
\begin{equation*}
  \prox_{\sigma F^*} (\Vek y)
  =
  (I + \sigma \, \uppartial F^*)^{-1}(\Vek y)
  = 
  \tfrac{1}{\sigma + 1} \, (\Vek y - \sigma \, \Vek g^\epsilon),
\end{equation*}
which yields the following algorithm.

\begin{Algorithm}[\PTikhonov\ regularization]
  \label{alg:tikh-reg}
  \begin{enumerate}[(i)]
  \item Initiation: Fix the parameters $\tau, \sigma > 0$ and
    $\theta \in [0,1]$.  Choose an arbitrary start value
    $(\Vek w^{(0)}, \Vek y^{(0)})$ in
    $(\mathcal H_1 \otimes \mathcal H_2) \times \mathcal
    K$, and set $\breve{\Vek w}^{(0)}$ to $\Vek w^{(0)}$. 
  \item Iteration: For $n>0$, update $\Vek w^{(n)}$, $\breve{\Vek
      w}^{(n)}$, and $\Vek y^{(n)}$ by 
    \begin{align*}
      \Vek y^{(n+1)}
      &\coloneqq
        \tfrac{1}{\sigma + 1} \bigl(
        \Vek y^{(n)} + \sigma \, (\breve{\mathcal B}(\breve{\Vek w}^{(n)}) -
        \Vek g^\epsilon) \bigr)
      \\[\fskip]
      \Vek w^{(n+1)}
      &\coloneqq
        \mathcal S_{\tau\alpha}\bigl(\Vek w^{(n)} - \tau \,
        \breve{\mathcal B}^* (\Vek 
        y^{(n+1)})\bigr)
      \\[\fskip]
      \breve{\Vek w}^{(n+1)}
      &\coloneqq
        \Vek w^{(n+1)} + \theta \, ( \Vek w^{(n+1)} - \Vek w^{(n)}).
    \end{align*}
  \end{enumerate}  
\end{Algorithm}

\begin{Remark}  
  \label{rem:tikh-reg}
  Since the data fidelity $F$ for the \PTikhonov\ functional is
  differentiable, one may here apply the forward-backward splitting as
  an alternative for the primal-dual iteration.  In so doing, the
  whole iteration in \thref{alg:tikh-reg}.ii becomes
  \begin{equation*}
      \Vek w^{(n+1)}
      \coloneqq
      \mathcal S_{\tau\alpha}\bigl(\Vek w^{(n)} - \tau \,
      \breve{\mathcal B}^* ( \breve{\mathcal B} \Vek 
      w^{(n)} - \Vek g^\epsilon)\bigr).
  \end{equation*}
  Analogously, one can here apply FISTA \cite{BT09} in order to improve the
  convergence.  \qed
\end{Remark}

If we incorporate the measurement errors by extending the solution
space as in \eqref{eq:inexact-data-bilin}, then the data fidelity is
chosen by
$F(\Vek y) \coloneqq \chi_{\epsilon \mathbb B_{\mathcal K}}(\Vek y -
\Vek g^\epsilon)$, where $\mathbb B_{\mathcal K}$ denotes the closed
unit ball in the \PHilbert\ space $\mathcal K$, and
$\chi_{\epsilon \mathbb B_{\mathcal K}}$ is the indicator functional
of the closed $\epsilon$-ball in $\mathcal K$, i.e.,
$\chi_{\epsilon \mathbb B_{\mathcal K}}(\Vek y) = 0$ if
$\pNormn{\Vek y}_{\mathcal K} \leq \epsilon$ and $\infty$ otherwise.
Since the conjugation of the unit ball yields the corresponding norm,
we obtain
$F^*(\Vek y) = \epsilon \, \pNormn{\Vek y}_{\mathcal K} + \iProd{\Vek
  y}{\Vek g^\epsilon}$ with subdifferential
\begin{equation}
  \label{eq:subdiff-hilb-norm}
  \uppartial F^*(\Vek y) =
  \begin{cases}
    \bigl\{\tfrac{\epsilon \Vek y}{\pNormn{\Vek y}_{\mathcal K}} + \Vek
    g^\epsilon \bigr\}
    & \text{if} \; \Vek y \ne \Vek 0,
    \\[\fskip]
    \epsilon \, \mathbb B_{\mathcal K} + \Vek g^\epsilon
    & \text{if} \; \Vek y = \Vek 0,
  \end{cases}
\end{equation}
\cf\ \cite[Exercise~8.27]{RW09}.  Since the proximation is not as
simple as in the previous cases, we give a more detailed computation.

\begin{Lemma}[Proximity operator]
  \label{lem:resol-hilb-norm}
  Let the functional $F \colon \mathcal K \to \overline \R$ be defined
  by
  $F(\Vek y) \coloneqq \chi_{\epsilon \mathbb B_{\mathcal K}}(\Vek y -
  \Vek g^\epsilon)$.  The proximation of $F^*$ is then given by
  \begin{equation*}
    \prox_{\sigma F^*} (\Vek y)
        = 
    \begin{cases}
      \Vek 0
      & \text{if} \; \pNormn{\Vek y - \sigma \Vek
        g^\epsilon}_{\mathcal K} \le \sigma \epsilon,
      \\[\fskip]
      \bigl(1 - \tfrac{\sigma \epsilon}{\pNormn{\Vek y - \sigma \Vek
          g^\epsilon}_{\mathcal K}} \bigr) \, (\Vek y - \sigma \Vek g^\epsilon)
      & \text{otherwise.}
    \end{cases}
  \end{equation*}
\end{Lemma}

\begin{Proof}
  The vector $\breve{\Vek y}$ is the resolvent
  $(I + \sigma \, \uppartial F^*)^{-1}(\Vek y)$ if and only if
  \begin{equation*}
    \Vek y \in
    \begin{cases}
      \bigl\{ \breve{\Vek y} + \sigma \,
      \bigl(\tfrac{\epsilon \breve{\Vek y}}{\pNormn{\Vek y}_{\mathcal
          K}} + \Vek g^\epsilon \bigr) \bigr\}  
      & \text{if} \; \breve{\Vek y} \ne \Vek 0,
      \\[\fskip]
      \sigma \, ( \epsilon \mathbb B_{\mathcal K} + \Vek g^\epsilon )
      & \text{if} \; \breve{\Vek y} = \Vek 0,
    \end{cases}
    \addmathskip
  \end{equation*}
  which is an immediate consequence of~\eqref{eq:subdiff-hilb-norm}.
  Bringing $\sigma \Vek g^\epsilon$ to the left-hand side, we are
  looking for a $\breve{\Vek y}$ such that
  \begin{equation*}
    \Vek y - \sigma \Vek g^\epsilon \in
    \begin{cases}
      \bigl\{ \bigl( 1 + \tfrac{\sigma \epsilon}{\pNormn{\breve{\Vek
            y}}_{\mathcal K}} \bigr) \, \breve{\Vek y} \bigr\}
      & \text{if} \; \breve{\Vek y} \ne \Vek 0,
      \\[\fskip]
      \sigma \epsilon \mathbb B_{\mathcal K}
      & \text{if} \; \breve{\Vek y} = \Vek 0.
    \end{cases}
    \addmathskip
  \end{equation*}
  For
  $\pNormn{\Vek y - \sigma \Vek g^\epsilon}_{\mathcal K} \le \sigma
  \epsilon$, the last condition is fulfilled for
  $\breve{\Vek y} = \Vek 0$.  Otherwise, it follows that
  $\breve{\Vek y} = \gamma \, (\Vek y - \sigma \Vek g^\epsilon)$ for
  some $\gamma > 0$.  With the notation
  $\Vek z \coloneqq \Vek y - \sigma \Vek g^\epsilon$, the first
  condition becomes
  \begin{equation*}
    \Vek z = \bigl( 1 + \tfrac{\sigma \epsilon}{\gamma \, \pNormn{\Vek
        z}_{\mathcal K}} \bigr) \, \gamma \Vek z
    \qquad\text{or}\qquad
    \bigl(1 - \tfrac{\sigma \epsilon}{\pNormn{\Vek z}_{\mathcal K}}
    \bigr) \, \Vek z = \gamma \Vek z.
  \end{equation*}
  Since $\pNormn{\Vek z}_{\mathcal K} > \sigma \epsilon$, we obtain
  $\gamma = 1 - \nicefrac{(\sigma \epsilon)}{\pNormn{\Vek z}_{\mathcal
      K}}$, and consequently, the assertion. \qed
\end{Proof}

\begin{Remark}
  \label{rem:resol-hilb-norm}
  The central part of the resolvent in \thref{lem:resol-hilb-norm} is
  given by the operator $\mathcal P_\gamma \colon \mathcal K \to
  \mathcal K$ with
  \begin{equation*}
    \mathcal P_\gamma (\Vek z) 
    \coloneqq
    \begin{cases}
      \Vek 0 
      & \text{if} \; \pNormn{\Vek z}_{\mathcal K} \le \gamma,
      \\[\fskip]
      \bigl( 1 - \tfrac{\gamma}{\pNormn{\Vek z}_{\mathcal K}} \bigr)
      \, \Vek z
      & \text{otherwise.}
    \end{cases}
    \addmathskip
  \end{equation*}
  Pictorially, this operator may be interpreted as shrinkage or
  contraction around the origin.  \qed
\end{Remark}

After this small digression to compute the proximation of the
conjugated data fidelity, the minimization problem
\eqref{eq:inexact-data-bilin} may be solved by the following
primal-dual iteration.

\begin{Algorithm}[Primal-dual for inexact data]
  \label{alg:inexact-data}
  \begin{enumerate}[(i)]
  \item Initiation: Fix the parameters $\tau, \sigma, \epsilon > 0$ and
    $\theta \in [0,1]$.  Choose an arbitrary start value
    $(\Vek w^{(0)}, \Vek y^{(0)})$ in
    $(\mathcal H_1 \otimes \mathcal H_2) \times \mathcal
    K$, and set $\breve{\Vek w}^{(0)}$ to $\Vek w^{(0)}$.
  \item Iteration: For $n>0$, update $\Vek w^{(n)}$, $\breve{\Vek
      w}^{(n)}$, and $\Vek y^{(n)}$ by 
    \begin{align*}
      \Vek y^{(n+1)}
      &\coloneqq
        \mathcal P_{\sigma \epsilon}\bigl(
        \Vek y^{(n)} + \sigma \, ( \breve{\mathcal B}(\breve{\Vek w}^{(n)}) -
        \Vek g^\epsilon )\bigr)
      \\[\fskip]
      \Vek w^{(n+1)}
      &\coloneqq
        \mathcal S_{\tau}\bigl(\Vek w^{(n)} - \tau \, \breve{\mathcal B}^* (\Vek
        y^{(n+1)})\bigr)
      \\[\fskip]
      \breve{\Vek w}^{(n+1)}
      &\coloneqq
        \Vek w^{(n+1)} + \theta \, ( \Vek w^{(n+1)} - \Vek w^{(n)}).
    \end{align*}
  \end{enumerate}  
\end{Algorithm}

The weighing between data fidelity and regularization in
\thref{rem:exact-data} analogously holds for \thref{alg:inexact-data}.

The central differences between the primal-dual iterations in
\thref{alg:exact-data}, \ref{alg:tikh-reg}, and \ref{alg:inexact-data}
for the minimization problems \eqref{eq:nucl-min-bilin},
\eqref{eq:tikh-func-bilin}, and \eqref{eq:inexact-data-bilin}
are contained in the dual update of $\Vek y^{(n+1)}$.  If the
parameter $\sigma$ is chosen close to zero, the three iterations nearly
coincide.  Thus, all three iterations should yield similar results; so
\thref{alg:exact-data} should also be able to deal with noisy
measurements.

\begin{Remark}
  \label{rem:quad-alg}
  For the relaxations \eqref{eq:nucl-min-quad},
  \eqref{eq:tikh-func-quad}, and \eqref{eq:inexact-data-quad} of the
  quadratic inverse problem \eqref{eq:gen-quad-prob}, we can derive
  analogous algorithms.  We do not state these variants more closely
  since the only differences are the application of the eigenvalue
  thresholding $\mathcal S_\tau^+$ instead of the singular value
  thresholding $\mathcal S_\tau$ and the usage of the quadratic
  lifting $\breve{\mathcal Q}$ instead of the bilinear lifting
  $\breve{\mathcal B}$.  Up to these two small modifications the
  methods completely coincide with the derived algorithms for bilinear
  inverse problems.  \qed
\end{Remark}

\section{Tensor-free singular value thresholding}
\label{sec:circ-tens-prod}

Each of the proposed methods solving bilinear or quadratic inverse
problems is based on a singular value thresholding on the tensor
product $\mathcal H_1 \otimes \mathcal H_2$ or
$\mathcal H \otimes_{\mathrm{sym}} \mathcal H$
respectively.  If the dimension of the original space
$\mathcal H_1 \times \mathcal H_2$ or $\mathcal H$ is already
enormous, then the dimension of the tensor product literally explodes,
which makes the computation of the required singular value
decomposition impracticable.  This difficulty occurs for nearly all
bilinear and quadratic image recovery problems.  However, since the
tensor $\Vek w^{(n)}$ is generated by a singular value thresholding,
the iterations $\Vek w^{(n)}$ usually possesses a very low rank.
Hence, the involved tensors can be stored in an efficient and
storage-saving manner.  In order to determine this low-rank
representation, we only compute a partial singular value decomposition
of the argument $\Vek w$ of $\mathcal S_\tau$ by deriving iterative
algorithms only requiring the left- and right-hand actions of $\Vek
w$.

Our first algorithm is based on the orthogonal iteration with \PRitz\
acceleration, see \cite{Ste69, GV13}.  In order to compute the leading
$\ell$ singular values, the main idea is here a joint power iteration
over two $\ell$-dimensional subspaces
$\widetilde{\mathcal U}_n \subset \mathcal H_1$ and
$\widetilde{\mathcal V}_n \subset \mathcal H_2$ alternately generated
by
$\widetilde{\mathcal U}_n \coloneqq \Vek w^* \Mat H_2
\widetilde{\mathcal V}_{n-1}$ and
$\widetilde{\mathcal V}_n \coloneqq \Vek w \Mat H_1
\widetilde{\mathcal U}_{n}$.  These subspaces are represented by
orthonormal bases
$\widetilde{\Mat U}_n \coloneqq [\widetilde{\Vek u}_0^{(n)}, \dots,
\widetilde{\Vek u}_{\ell-1}^{(n)}]$ in $\mathcal H_1$ and
$\widetilde{\Mat V}_n \coloneqq [\widetilde{\Vek v}_0^{(n)}, \dots,
\widetilde{\Vek v}_{\ell-1}^{(n)}]$ in $\mathcal H_2$.

\begin{Algorithm}[Subspace iteration]
  \label{alg:subspace-iteration}
   {\scshape Input}:
    $\Vek w \in \mathcal H_1 \otimes \mathcal H_2$, $\ell > 0$,
    $\delta > 0$.
  \begin{enumerate}[(i)]
  \item Choose
    $\widetilde{\Vek V}_0 \in \R^{N_2 \times \ell}$, whose columns are
    orthonormal with respect to $\mathcal H_2$.
  \item For $n>0$, 
  repeat:
    \begin{enumerate}[(a)]
    \item Compute $\widetilde{\Mat E}_n \coloneqq \Mat w^* \Mat H_2
      \widetilde{\Mat V}_{n-1}$, and reorthonormalize the columns in
      $\mathcal H_1$.
    \item Compute $\widetilde{\Mat F}_n \coloneqq \Mat w \Mat H_1
      \widetilde{\Mat E}_{n}$, and reorthonormalize the columns in
      $\mathcal H_2$. 
    \item Determine the \PEuclid{ian} singular value decomposition
      \begin{equation*}
        \widetilde{\Mat F}_n^* \Mat H_2 \Vek w \Mat H_1 \widetilde{\Mat
          E}_n = \Mat Y_n \Mat \Sigma_n \Mat Z_n^*,
      \end{equation*}
      and set
      $\widetilde{\Mat U}_n \coloneqq \widetilde{\Mat E}_n \Mat Z_n$
      and $\widetilde{\Mat V}_n \coloneqq \widetilde{\Mat F}_n \Mat Y_n$.
    \end{enumerate}
    until $\ell$ singular vectors have converged, which
    means
    \begin{equation*}
      \pNorm{\Vek w^* \Mat H_2
        \widetilde{\Vek v}_m^{(n)} - \sigma_m^{(n)} \widetilde{\Vek
          u}_m^{(n)}}_{\mathcal H_1} \le \delta \, \pNormn{\Vek
        w}_{\mathcal L(\mathcal H_1, \mathcal H_2)}
      \qquad\text{for}\qquad
      0 \le m < \ell,
    \end{equation*}
    where $\pNormn{\cdot}_{\mathcal L(\mathcal H_1, \mathcal H_2)}$
    denotes the operator norm estimated by $\sigma_0^{(n)}$.
  \end{enumerate}
  {\scshape Output}:
  $\widetilde{\Mat U}_n \in \R^{N_1 \times \ell}$,
  $\widetilde{\Mat V}_n \in \R^{N_2 \times \ell}$,
  $\Mat \Sigma_n \in \R^{\ell \times \ell}$ with
  $\widetilde{\Mat V}_n^* \Mat H_2 \Vek w \Mat H_1 \widetilde{\Mat
    U}_n = \Mat \Sigma_n$.
\end{Algorithm}

Here, reorthonormalization means that for each applicable $m$, the
span of the first $m$ columns of the matrix and its
reorthonormalization coincide, and that the reorthonormalized matrix
has orthonormal columns. This can, for instance, be achieved by the
well-known \pers{Gram--Schmidt} procedure.

Under mild conditions on the subspace associated to
$\widetilde{\Mat V}_0$, the matrices $\widetilde{\Vek U}_n$,
$\widetilde{\Vek V}_n$, and
\raisebox{0pt}[0pt][0pt]{$\Mat \Sigma_n \coloneqq \diag
  (\sigma_0^{(n)}, \dots, \sigma_{\ell-1}^{(n)})$}
converge to leading singular vectors as well as to the leading
singular values of a singular value decomposition
\raisebox{0pt}[0pt][0pt]{$\Vek w = \sum_{n=0}^{R-1} \sigma_n \,
  (\widetilde{\Vek u}_n \otimes \widetilde{\Vek v}_n)$}.

\begin{Theorem}[Subspace iteration]
  \label{the:subspace-iteration}
  If none of the basis vectors in $\widetilde{\Mat V}_0$ is orthogonal
  to the $\ell$ leading singular vectors
  $\widetilde{\Vek v}_0, \dots, \widetilde{\Vek v}_{\ell-1}$, and if
  $\sigma_{\ell-1} > \sigma_\ell$, then the singular values
  $\sigma_0^{(n)} \ge \cdots \ge \sigma_{\ell-1}^{(n)}$ in
  \thref{alg:subspace-iteration} converge to
  $\sigma_0 \ge \cdots \ge \sigma_{\ell-1}$ with a rate of
  \begin{equation*}
    \absb{\bigl[\sigma_m^{(n)} \bigr]^2 - \sigma_m^2}
    =
    \Landau \Bigl( \absB{\frac{\sigma_\ell}{\sigma_m}}^{2k} \Bigr)
    \qquad\text{and}\qquad
    \sigma^{(n)}_m \le \sigma_m.
  \end{equation*}
\end{Theorem}

\begin{Proof}
  By the construction in steps (a) and (b), the columns in
  $\widetilde{\Mat E}_n$ and $\widetilde{\Mat F}_n$ form orthonormal
  systems in $\mathcal H_1$ and $\mathcal H_2$.  In this proof, we
  denote the corresponding subspaces by $\widetilde{\mathcal E}_n$ and
  $\widetilde{\mathcal F}_n$, which are related by
  $\widetilde{\mathcal E}_n = \Vek w^* \Mat H_2 \widetilde{\mathcal
    F}_{n-1}$ and
  $\widetilde{\mathcal F}_n = \Vek w \Mat H_1 \widetilde{\mathcal
    E}_n$.  Due to the basis transformation in (c), the columns of
  $\widetilde{\Mat U}_n$ and $\widetilde{\Mat V}_n$ also form
  orthonormal bases of $\widetilde{\mathcal E}_n$ and
  $\widetilde{\mathcal F}_n$.  Next, we exploit that the projection
  $\Mat P_n \coloneqq \widetilde{\Mat V}_n \widetilde{\Mat V}_n^* \Mat
  H_2$ onto $\widetilde{\mathcal F}_n$ acts as identity on
  $\Vek w \Mat H_1 \widetilde{\mathcal E}_n$ by construction.  Since
  $\widetilde{\Mat U}_n$ is a basis of $\widetilde{\mathcal E}_n$, and
  since
  $\widetilde{\Mat V}_n^* \Mat H_2 \Vek w \Mat H_1 \widetilde{\Mat
    U}_n = \Mat \Sigma_n$ by the singular value decomposition in step
  (c), we have
  \begin{equation}
    \label{subspace-svd}
    \begin{aligned}
      \widetilde{\Mat U}_n^* \Mat H_1 \Vek w^* \Mat H_2 \Vek w \Mat
      H_1 \widetilde{\Mat U}_n
      &=
      \widetilde{\Mat U}_n^* \Mat H_1 \Vek
      w^* \Mat H_2 \Mat P_n \Vek w \Mat H_1 \widetilde{\Mat U}_n
      \\[\fskip]
      &=
      \widetilde{\Mat U}_n^* \Mat H_1 \Vek w^* \Mat H_2
      \widetilde{\Mat V}_n \widetilde{\Mat V}_n^* \Mat H_2 \Vek w \Mat
      H_1 \widetilde{\Mat U}_n = \Mat \Sigma_n^2,
    \end{aligned}
  \end{equation}
  and $\widetilde{\Mat U}_n$ diagonalizes
  $\Mat H_1 \Vek w^* \Mat H_2 \Vek w \Mat H_1$ on the subspace
  $\widetilde{\mathcal E}_n$.

  Using the substitutions
  \begin{equation*}
    \Mat E_n \coloneqq \Mat H_1^{\nicefrac12} \, \widetilde{\Mat E}_n,
    \quad
    \Mat U_n \coloneqq \Mat H_1^{\nicefrac12} \, \widetilde{\Mat U}_n,
    \quad
    \Mat F_n \coloneqq \Mat H_2^{\nicefrac12} \, \widetilde{\Mat F}_n,
    \quad
    \Mat V_n \coloneqq \Mat H_2^{\nicefrac12} \, \widetilde{\Mat V}_n
  \end{equation*}
  as well as
  \begin{equation*}
    \mathcal E_n = \Mat H_1^{\nicefrac12} \, \widetilde{\mathcal E}_n
    \quad\text{and}\quad
    \mathcal F_n = \Mat H_2^{\nicefrac12} \, \widetilde{\mathcal F}_n,
    \addmathskip
  \end{equation*}
  we notice that the iteration in \thref{alg:subspace-iteration} is
  composed of two main steps.  First, in (a) and (b), we compute an
  orthonormal basis $\Mat E_n$ of
  \begin{equation*}
    \mathcal E_n = (\Mat H_1^{\nicefrac12} \Vek w^*  (\Mat
    H_2^{\nicefrac12})^*) (\Mat H_2^{\nicefrac12} \Vek w \, (\Mat
    H_1^{\nicefrac12})^*) \, \mathcal E_{n-1}.
  \end{equation*}
  Secondly, \eqref{subspace-svd} implies that we determine an
  \PEuclid{ian} eigenvalue decomposition on the subspace
  $\mathcal E_n$ by
  \begin{equation*}
    \Mat E_n^* (\Mat H_1^{\nicefrac12} \Vek w^* (\Mat H_2^{\nicefrac12})^*
    )
    (\Mat H_2^{\nicefrac12} \Vek w \, (\Mat H_1^{\nicefrac12} )^*) \Mat E_n 
    =
    \Mat Z_n \Mat \Sigma_n^2 \Mat Z_n^*
  \end{equation*}
  and $\Mat U_n \coloneqq \Mat E_n \Mat Z_n$.

  This two-step iteration exactly coincides with the orthogonal
  iteration with \PRitz\ acceleration for the matrix
  $ (\Mat H_1^{\nicefrac12} \Vek w^* (\Mat H_2^{\nicefrac12})^* ) (\Mat
  H_2^{\nicefrac12} \Vek w \, (\Mat H_1^{\nicefrac12} )^*)$,
  see \cite{Ste69, GV13}.  Under the given assumptions, this iteration
  converges to the $\ell$ leading eigenvalues and eigenvectors with
  the asserted rates.  In view of \thref{lem:svd}, the columns in
  $\widetilde{\Mat U}_n$ and $\widetilde{\Mat V}_n$ together with
  $\Mat \Sigma_n$ converge to the leading components of the singular
  value decomposition of $\Mat w$ with respect to $\mathcal H_1$ and
  $\mathcal H_2$.  \qed
\end{Proof}

Considering the subspace iteration (\thref{alg:subspace-iteration}),
notice that the algorithm does not need an explicit representation of
its argument $\Vek w$ but the left- and right-hand actions of $\Vek w$
as a matrix-vector multiplication.  We may thus use the subspace
iteration to compute the singular value thresholding
$\mathcal S_\tau(\Vek w)$ without a tensor representation of $\Vek w$.

\begin{Algorithm}[Tensor-free singular value thresholding]
  \label{alg:tf-threshold:subspace}
  {\scshape Input}: $\Vek w \in \mathcal H_1 \otimes \mathcal H_2$,
  $\tau > 0$, $\ell > 0$, $\delta > 0$.
  \begin{enumerate}[(i)]
  \item Apply \thref{alg:subspace-iteration} with the following
    modifications:
    \begin{itemize}
    \item If $\sigma_m^{(n)} > \tau$ for all $0 \le m < \ell$,
      increase $\ell$ and extend $\widetilde{\Mat V}_n$ by further
      orthonormal columns, unless $\ell = \rank \Vek w$, \ie{}, when
      the columns of $\widetilde{\Mat E}_n$ would become linearly
      dependent.
    \item Additionally, stop the subspace iterations when the first
      $\ell' + 1$ singular values with $\ell' < \ell$ have converged
      and \raisebox{0pt}[0pt][0pt]{$\sigma_{\ell'+1}^{(n)} < \tau$}.
      Otherwise, continue the iteration until all non-zero singular values
      converge and set $\ell' = \ell$.
    \end{itemize}
  \item Set
    $\widetilde{\Mat U}' \coloneqq [\widetilde{\Vek u}_0, \dots,
    \widetilde{\Vek u}_{\ell'-1}]$,
    $\widetilde{\Mat V}' \coloneqq [\widetilde{\Vek v}_0, \dots,
    \widetilde{\Vek v}_{\ell'-1}]$, and
    \begin{equation*}
      \Mat \Sigma' \coloneqq
      \diag\bigl(S_\tau \bigl(\sigma_0^{(n)}\bigr), \dots,
      S_\tau\bigl(\sigma_{\ell'}^{(n)}\bigr)\bigr).
  \end{equation*}
  \end{enumerate}
  {\scshape Output}: $\widetilde{\Mat U}' \in \R^{N_1 \times \ell'}$,
  $\widetilde{\Mat V}' \in \R^{N_2 \times \ell'}$,
  $\Mat \Sigma' \in \R^{\ell' \times \ell'}$ with
  $\widetilde{\Mat V}' \Mat \Sigma' (\widetilde{\Mat U}')^* = \mathcal
  S_\tau(\Vek w)$.
\end{Algorithm}

\begin{Corollary}[Exact singular value thresholding] 
  \label{cor:tf-threshold:subspace}
  If the non-zero singular values of $\Vek w$ are distinct, and if
  none of the columns in $\widetilde{\Mat V}_n$ is orthogonal to the
  singular vectors with $\sigma_n > \tau$, then
  \thref{alg:tf-threshold:subspace} computes the low-rank
  representation of $\mathcal S_\tau(\Vek w)$.
\end{Corollary}

Although \thref{alg:tf-threshold:subspace} for generic start values
always yields the singular value thresh\-old\-ing, the convergence of the
subspace iteration is rather slow.  Therefore, we now derive an
algorithm that is based on the \PLanczos-based bidiagonalization
method proposed by \pers{Golub} and \pers{Kahan} in \cite{GK65} and
the \pers{Ritz} approximation in \cite{GV13}.  This method again only
require the left-hand and right-hand action of $\Vek w$ with respect
to a given vector.  For simplifying the following considerations, we
initially present the employed \PLanczos\ process with respect to the
\PEuclid{ian} singular value decomposition.

\sloppy The central idea is here to construct, for fixed $k$, orthonormal
matrices
$\Mat F_k = [\Vek f_0, \dots, \Vek f_{k-1}] \in \R^{N_2 \times k}$ and
$\Mat E_k = [\Vek e_0, \dots, \Vek e_{k-1}] \in \R^{N_1 \times k}$
such that the transformed matrix
\begin{equation}
  \label{eq:bidiag-Bk}
  \Mat F_k^* \Vek w \Mat E_k =
  \Mat B_k =
  \begin{bmatrix}
    \beta_0 & \gamma_0 & & & \\
    & \beta_1 & \gamma_1 & & \\
    & & \ddots & \ddots & \\
    & & & \beta_{k-2} & \gamma_{k-2} \\
    & & & & \beta_{k-1}
  \end{bmatrix}
\end{equation}
is bidiagonal, and then to compute the singular value decomposition of
$\Mat B_k$ by determining orthogonal matrices $\Mat Y_k$, $\Mat Z_k$,
and $\Mat \Sigma_k$ in $\R^{k \times k}$ such that
\begin{equation*}
  \Mat Y_k^* \Mat B_k \Mat Z_k 
  =
  \Mat \Sigma_k
  =
  \diag(\sigma_0, \dots, \sigma_{k-1}).
\end{equation*}
Defining $\Vek U_k \in \R^{N_1 \times k}$ and
$\Vek V_k \in \R^{N_2 \times k}$ as
\begin{equation*}
  \Mat U_k \coloneqq \Mat E_k \Mat Z_k
  \qquad\text{and}\qquad
  \Mat V_k \coloneqq \Mat F_k \Mat Y_k,
\end{equation*}
we finally obtain a set of approximate right-hand and left-hand
singular vectors, see \cite{GK65, BR05, GV13}.

The values $\beta_n$ and $\gamma_n$ of the bidiagonal matrix
$\Mat B_k$ and the related vectors $\Vek e_n$ and $\Vek f_n$ can be
determined by the following iterative procedure \cite{GK65}: Choose
an arbitrary unit vector $\Vek p_{-1} \in \R^{N_1}$ with respect to the
\PEuclid{ian} norm, and compute
\begin{equation*}
  \left.
    \begin{aligned}
      \Vek e_{m+1} &\coloneqq \gamma_m^{-1} \, \Vek p_m, 
      \\
      \Vek q_{m+1} &\coloneqq \Vek w \Vek e_{m+1} - \gamma_m \Vek f_m,
      \\
      \beta_{m+1} &\coloneqq \pNormn{\Vek q_{m+1}},
    \end{aligned}
    \qquad\middle|\qquad
    \begin{aligned}
      \Vek f_{m+1} &\coloneqq \beta_{m+1}^{-1} \Vek q_{m+1},
      \\
      \Vek p_{m+1} &\coloneqq \Vek w^* \Vek f_{m+1} - \beta_{m+1} \Vek
      e_{m+1},
      \\
      \gamma_{m+1} &\coloneqq \pNormn{\Vek p_{m+1}}.
    \end{aligned}
  \right.
\end{equation*}
For the first iteration, we set $\gamma_{-1} \coloneqq 1$ and
$\Vek f_{-1} \coloneqq \Vek 0$. If $\gamma_{m+1}$ vanishes, then we
stop the \PLanczos\ process since we have found an invariant \PKrylov\
subspace such that the computed singular values become exact.

In order to compute an approximate singular value decomposition with
respect to the \PHilbert\ spaces $\mathcal H_1$ and $\mathcal H_2$, we
exploit \thref{lem:svd} and perform the \PLanczos\ bidiagonalization
regarding the transformed matrix
$\Mat H_2^{\nicefrac12} \Vek w \, (\Mat H_1^{\nicefrac12})^*$.
Moreover, we incorporate the back transformation in \thref{lem:svd}
with the aid of the substitutions 
\begin{equation}
  \label{eq:sub-lanc-proc}
  \widetilde{\Vek e}_m 
  \coloneqq 
  \Mat H_1^{-\nicefrac12} \Vek e_m,
  \quad
  \widetilde{\Vek p}_m 
  \coloneqq 
  \Mat H_1^{-\nicefrac12} \Vek p_m
  \quad\text{and}\quad
  \widetilde{\Vek f}_m 
  \coloneqq 
  \Mat H_2^{-\nicefrac12} \Vek f_m,
  \quad
  \widetilde{\Vek q}_m 
  \coloneqq 
  \Mat H_2^{-\nicefrac12} \Vek q_m.
\end{equation}
In this manner, the square roots $\Mat H_1^{\nicefrac12}$ and $\Mat
H_2^{\nicefrac12}$ and their inverses cancel out, and we obtain the
following algorithm, which only relies on the original matrices $\Mat
H_1$ and $\Mat H_2$.

\begin{Algorithm}[\PLanczos\ bidiagonalization]
  \label{alg:lanc-bidiag}
  {\scshape Input}: $\Vek w \in \mathcal H_1 \otimes \mathcal H_2$, $k>0$.
  \begin{enumerate}[(i)]
  \item Initiation: Set $\gamma_{-1} \coloneqq 1$ and
    $\widetilde{\Vek f}_{-1} \coloneqq \Vek 0$.  Choose a unit vector
    $\widetilde{\Vek p}_{-1} \in \mathcal H_1$.
  \item \PLanczos\ bidiagonalization: For $m = -1, \dots, k-2$ while
    $\gamma_m \ne 0$, repeat:
    \begin{enumerate}[(a)]
    \item Compute
      $\widetilde{\Vek e}_{m+1} \coloneqq \gamma_m^{-1} \,
      \widetilde{\Vek p}_m$,
      and reorthogonalize with $\widetilde{\Vek e}_0, \dots,
      \widetilde{\Vek e}_m$ in $\mathcal H_1$.
    \item Determine
      $\widetilde{\Vek q}_{m+1} \coloneqq \Vek w \Mat H_1
      \widetilde{\Vek e}_{m+1} - \gamma_m \widetilde{\Vek f}_m$,
      and set
      $\beta_{m+1} \coloneqq \pNormn{\widetilde{\Vek q}_{m+1}}_{\mathcal H_2}$.
      \newline Compute
      $\widetilde{\Vek f}_{m+1} \coloneqq \beta_{m+1}^{-1}
      \widetilde{\Vek q}_{m+1}$
      and reorthogonalize with
      $\widetilde{\Vek f}_0, \dots, \widetilde{\Vek f}_m$ in
      $\mathcal H_2$.
    \item Determine
      $\widetilde{\Vek p}_{m+1} \coloneqq \Vek w^* \Mat H_2
      \widetilde{\Vek f}_{m+1} - \beta_{m+1} \widetilde{\Vek
        e}_{m+1}$,
      and set
      $\gamma_{m+1} \coloneqq \pNormn{\widetilde{\Vek p}_{m+1}}$.
    \end{enumerate}
  \item Compute the \PEuclid{ian} singular value decomposition of
    $\Mat B_k$ according to \eqref{eq:bidiag-Bk}, \ie\
    $\Mat B_k = \Mat Y_k \Mat \Sigma_k \Mat Z_k^*$, and set
    $\widetilde{\Mat U}_k \coloneqq \widetilde{\Mat E}_k \Mat Z_k$ and
    $\widetilde{\Vek V}_k \coloneqq \widetilde{\Mat F}_k \Mat Y_k$.
  \end{enumerate}
  {\scshape Output}: $\widetilde{\Mat U}_k \in \R^{N_1 \times k}$,
  $\widetilde{\Mat V}_k \in \R^{N_2 \times k}$, $\Mat \Sigma_k \in \R^{k \times k}$ with
  $\widetilde{\Mat V}_k^* \Mat H_2 \Vek w \Mat H_1 \widetilde{\Mat U}_k
  = \Mat \Sigma_k$.
\end{Algorithm}

\begin{Remark}
  \label{rem:lanc-bidiag:1}
  The bidiagonalization by \pers{Golub} and \pers{Kahan} is based on a
  \PLanczos-type process, which is numerically unstable in the
  computation of $\widetilde{\Vek e}_n$ and
  \raisebox{0pt}[0pt][0pt]{$\widetilde{\Vek f}_n$}.  For this reason,
  we have to reorthogonalize all newly generated vectors
  $\widetilde{\Vek e}_n$ and
  \raisebox{0pt}[0pt][0pt]{$\widetilde{\Vek f}_n$} with the previously
  generated vectors, see \cite{GK65}.  This amounts to projecting
  $\widetilde{\Vek e}_{m+1}$ to the orthogonal complement of the span
  of $\{\widetilde{\Vek e}_{0}, \ldots, \widetilde{\Vek e}_{m} \}$ and
  the analog for $\widetilde{\Vek f}_{m+1}$, for instance, via the
  \pers{Gram--Schmidt} procedure.  \qed
\end{Remark}

\begin{Remark}
  \label{rem:lanc-bidiag:2}
  The computation of the last $\widetilde{\Vek p}_{k-1}$ seems to be superfluous
  since it is not needed for the determination of the matrix $\Mat
  B_k$.  On the other side, this vector represents the residuals of
  the approximate singular value decomposition.  More precisely, we
  have 
  \begin{equation}
    \label{eq:svd-approx-error}
    \Mat w \Mat H_1 \widetilde{\Vek u}_m = \sigma_m \widetilde{\Vek v}_m
    \qquad\text{and}\qquad
    \Mat w^* \Mat H_2 \widetilde{\Vek v}_m = \sigma_m \widetilde{\Vek
      u}_m + \widetilde{\Vek p}_{k-1} \Vek 
    \eta_{k-1}^* \Vek y_m
  \end{equation}
  for $m = 0, \dots, k-1$, see \cite{BR05}.  Here the vectors
  $\widetilde{\Vek u}_m$, $\widetilde{\Vek v}_m$, and $\Vek y_m$
  denote the columns of the matrices
  $\widetilde{\Mat U}_k = [\widetilde{\Vek u}_0, \dots,
  \widetilde{\Vek u}_{k-1}]$,
  $\widetilde{\Mat V}_k = [\widetilde{\Vek v}_0, \dots,
  \widetilde{\Vek v}_{k-1}]$, and
  $\Mat Y_k = [\Vek y_0, \dots, \Vek y_{k-1}]$ respectively; the
  singular values $\sigma_m$ of $\Mat B_k$ are given by
  $\Mat \Sigma_k = \diag(\sigma_0, \dots, \sigma_{k-1})$; the vector
  $\Vek \eta_{k-1} \in \R^k$ represents the last unit vector
  $(0, \dots, 0, 1)^*$. \qed
\end{Remark}

Since the bidiagonalization method by \pers{Golub} and \pers{Kahan} is
based on the \PLanczos\ process for symmetric matrices, one can apply
the related convergence theory to show that the approximate singular
values and singular vectors -- for increasing $k$ -- converge to the
wanted singular value decomposition of $\Vek w$, see \cite{GV13}.
Since we are only interested in the leading singular values and
singular vectors, and since we want to choose the matrix $\Mat B_k$ as
small as possible, this convergence theory does not apply to our
setting.

In order to improve the quality of the approximate singular value
decomposition computed by \thref{alg:lanc-bidiag}, we here use a
restarting technique proposed by \pers{Baglama} and \pers{Reichel}
\cite{BR05}.  The central idea is to adapt the \PLanczos\
bidiagonalization such that the method can be restarted by a set of
$\ell$ previously computed \PRitz\ vectors.  For this purpose,
\pers{Baglama} and \pers{Reichel} suggest a modified bidiagonalization
of the form
\begin{equation}
  \label{eq:aug-bidiag-Bk}
  \Mat F_{k,n}^* \Vek w \Mat E_{k,n}
  = \Mat B_{k,n} =
  \begin{bmatrix}
    \sigma_0^{(n-1)} & & & \rho_0^{(n)} \\
    & \ddots & & \vdots \\
    & & \sigma_{\ell-1}^{(n-1)} & \rho_{\ell - 1}^{(n)} \\
    & & & \beta_\ell^{(n)} & \gamma_\ell^{(n)} \\
    & & & & \ddots & \ddots \\
    & & & & & \beta_{k-2}^{(n)} & \gamma_{k-2}^{(n)} \\
    & & & & & & \beta_{k-1}^{(n)} \\
  \end{bmatrix},
\end{equation}
where the first $\ell$ columns of the orthonormal matrices
\begin{equation*}
  \Mat E_{k,n} = [\Vek u_0^{(n-1)}, \dots, \Vek u_{\ell - 1}^{(n-1)},
  \dots]
  \qquad\text{and}\qquad
  \Mat F_{k,n} = [\Vek v_0^{(n-1)}, \dots, \Vek v_{\ell - 1}^{(n-1)}, \dots ]
\end{equation*}
are predefined by the \PRitz\ vectors of the previous iteration.  For
the computation of the first $\ell < k$ leading singular values and
singular vectors, we employ the following algorithm \cite{BR05}, which
has been adapted to our setting by incorporating \thref{lem:svd} and
the substitution \eqref{eq:svd-approx-error}.

\begin{Algorithm}[Augmented \PLanczos\ Bidiagonalization]
  \label{alg:aug-lanc-bidiag}
  {\scshape Input}: $\Vek w \in \mathcal H_1 \otimes \mathcal H_2$,
  $\ell > 0$ $k> \ell$, $\delta>0$.
  \begin{enumerate}[(i)]
  \item Apply \thref{alg:lanc-bidiag} to compute an approximate
    singular value decomposition $\widetilde{\Mat V}_{k,0}^* \Mat H_2
    \Vek w \Mat H_1 \widetilde{\Mat U}_{k,0} = \Mat \Sigma_{k,0}$.
  \item For $n > 0$, until $\ell$ singular vectors have
    converged, which means
    \begin{equation*}
      \beta_{k-1}^{(n-1)} \absn{\Vek \eta_{k-1}^* \Vek y_m^{(n-1)}}
      \le \delta \pNormn{\Vek w}_{\mathcal L(\mathcal H_1,
        \mathcal H_2)}
      \qquad\text{for}\qquad
      0 \le m < \ell,
    \end{equation*}
    where
    $\pNormn{\cdot}_{\mathcal L(\mathcal H_1, \mathcal H_2)}$ denotes
    the operator norm estimated by $\sigma_0^{(n-1)}$, repeat:
    \begin{enumerate}[(a)]
    \item Initialize the new iteration by setting
      $\widetilde{\Vek e}_m^{(n)} \coloneqq \widetilde{\Vek
        u}_m^{(n-1)}$ and
      $\widetilde{\Vek f}\kern0pt_m^{(n)} \coloneqq \widetilde{\Vek
        v}_m^{(n-1)}$ for $m = 0, \dots, \ell - 1$.  Further, set
      $\widetilde{\Vek p}_{\ell - 1}^{(n)} \coloneqq \widetilde{\Vek
        p}_{k-1}^{(n-1)}$ and
      $\gamma_{\ell - 1}^{(n)} \coloneqq \pNormn{\widetilde{\Vek
          p}_{\ell - 1}^{(n)}}_{\mathcal H_1}$.
    \item Compute
      $\widetilde{\Vek e}_\ell^{(n)} \coloneqq (\gamma_{\ell -
        1}^{(n)})^{-1} \,
      \widetilde{\Vek p}_{\ell - 1}^{(n)}$,
      and reorthogonalize with
      $\widetilde{\Vek e}_0^{(n)}, \dots, \widetilde{\Vek e}_{\ell - 1}^{(n)}$ in
      $\mathcal H_1$.
    \item Determine
      $\widetilde{\Vek q}_{\ell}^{(n)} \coloneqq \Vek w \Mat H_1
      \widetilde{\Vek e}_\ell^{(n)}$, compute the inner products
      $\rho_m^{(n)} \coloneqq \iProdn{\widetilde{\Vek
          f}\kern0pt_m^{(n)}}{\widetilde{\Vek q}_\ell^{(n)}}_{\mathcal
        H_2}$ for $m = 0, \dots, \ell - 1$, and reorthogonalize
      $\widetilde{\Vek q}_\ell^{(n)}$ in $\mathcal H_2$ by
      \begin{equation*}
        \widetilde{\Vek q}_\ell^{(n)} \coloneqq \widetilde{\Vek
          q}_\ell^{(n)} - \sum_{m=0}^{\ell - 1} \rho_m^{(n)}
        \widetilde{\Vek f}\kern0pt_m^{(n)} .
      \end{equation*}
    \item Set
      $\beta_\ell^{(n)} \coloneqq \pNormn{\widetilde{\Vek
          q}_\ell^{(n)}}_{\mathcal H_2}$ and
      $\widetilde{\Vek f}\kern0pt_\ell^{(n)} \coloneqq
      (\beta_\ell^{(n)})^{-1} \, \widetilde{\Vek q}_\ell^{(n)}$.
    \item Determine
      $\widetilde{\Vek p}_\ell^{(n)} \coloneqq \Vek w^* \Mat H_2
      \widetilde{\Vek f}\kern0pt_\ell^{(n)} - \beta_\ell^{(n)}
      \widetilde{\Vek e}_\ell^{(n)}$, and set
      $\gamma_\ell^{(n)} \coloneqq \pNormn{\widetilde{\Vek
          p}_\ell^{(n)}}_{\mathcal H_1}$.
    \item Calculate the remaining values of $\Mat B_{k,n}$ by applying step (ii) of
      \thref{alg:lanc-bidiag} with $m = \ell, \dots, k - 2$.
    \item Compute the \PEuclid{ian} singular value decomposition of
      $\Mat B_{k,n}$ in~\eqref{eq:aug-bidiag-Bk}, \ie\
      $\Mat B_{k,n} = \Mat Y_{k,n} \Mat \Sigma_{k,n} \Mat Z_{k,n}^*$,
      and set
      $\widetilde{\Mat U}_{k,n} \coloneqq \widetilde{\Mat E}_{k,n}
      \Mat Z_{k,n}$ and
      $\widetilde{\Mat V}_{k,n} \coloneqq \widetilde{\Mat F}_{k,n}
      \Mat Y_{k,n}$.
    \end{enumerate}
  \item Set
    $\widetilde{\Mat U} \coloneqq [\widetilde{\Vek u}_0^{(n)}, \dots,
    \widetilde{\Vek u}_{\ell-1}^{(n)}]$,
    $\widetilde{\Mat V} \coloneqq[\widetilde{\Vek v}_0^{(n)}, \dots,
    \widetilde{\Vek v}_{\ell-1}^{(n)}]$, and
    $\Mat \Sigma \coloneqq \diag(\sigma_0^{(n)}, \dots,
    \sigma_{\ell-1}^{(n)})$.
  \end{enumerate}
  {\scshape Output}:
  $\widetilde{\Mat U} \in \R^{N_1 \times \ell}$,
  $\widetilde{\Mat V} \in \R^{N_2 \times \ell}$,
  $\Mat \Sigma \in \R^{\ell \times \ell}$ with
  $\widetilde{\Mat V}^* \Mat H_2 \Vek w \Mat H_1 \widetilde{\Mat
    U} = \Mat \Sigma$.
\end{Algorithm}

\begin{Remark}
  \label{rem:aug-lanc-bidiag:1}
  The stopping criterion in step (ii) originates from the error
  representation in \eqref{eq:svd-approx-error}.  For the operator
  norm
  $\pNormn{\widetilde{\Vek w}}_{\mathcal L(\mathcal H_1, \mathcal
    H_2)}$, one may use the maximal leading singular values of the
  former iterations, which usually gives a sufficiently good
  approximation, see \cite{BR05}.  \qed
\end{Remark}

Although the numerical effort of the restarted augmented \PLanczos\
process is enormously reduced compared with the subspace iteration, we
are unfortunately not aware of a convergence and error analysis for
this specific variant of \PLanczos-type method.  Nevertheless, we can
employ the obtained partial singular value decomposition to determine
the singular value thresholding.

\begin{Algorithm}[Tensor-free singular value thresholding]
  \label{alg:tf-threshold:lanczos}
  {\scshape Input}: $\Vek w \in \mathcal H_1 \otimes \mathcal H_2$,
  $\tau > 0$, $\ell > 0$, $k>\ell$, $\delta > 0$.
  \begin{enumerate}[(i)]
  \item Apply \thref{alg:aug-lanc-bidiag} with the following
    modifications:
    \begin{itemize}
    \item If $\sigma_m^{(n)} > \tau$ for all $0 \le m < \ell$,
      increase $\ell$ and $k$ with $\ell < k$, unless
      $k = \rank \Vek w$, \ie{}, when $\gamma_{k}^{(n)}$ in
      \thref{alg:lanc-bidiag} vanishes.
    \item Additionally, stop the augmented \PLanczos\ method when the
      first $\ell' + 1$ singular values with $\ell' < \ell$ have
      converged and
      \raisebox{0pt}[0pt][0pt]{$\sigma_{\ell'+1}^{(n)} < \tau$}.
      Otherwise, continue the iteration until all non-zero singular
      values converge and set $\ell' = \ell$.
    \end{itemize}
  \item Set
    $\widetilde{\Mat U}' \coloneqq [\widetilde{\Vek u}_0^{(n)}, \dots,
    \widetilde{\Vek u}_{\ell'-1}^{(n)}]$,
    $\widetilde{\Mat V}' \coloneqq [\widetilde{\Vek v}_0^{(n)}, \dots,
    \widetilde{\Vek v}_{\ell'-1}^{(n)}]$, and
    \begin{equation*}
      \Mat \Sigma' \coloneqq
      \diag\bigl(S_\tau \bigl(\sigma_0^{(n)}\bigr), \dots,
      S_\tau\bigl(\sigma_{\ell'}^{(n)}\bigr)\bigr).
  \end{equation*}
  \end{enumerate}
  {\scshape Output}: $\widetilde{\Mat U}' \in \R^{N_1 \times \ell'}$,
  $\widetilde{\Mat V}' \in \R^{N_2 \times \ell'}$,
  $\Mat \Sigma' \in \R^{\ell' \times \ell'}$ with
  $\widetilde{\Mat V}' \Mat \Sigma' (\widetilde{\Mat U}')^* = \mathcal
  S_\tau(\Vek w)$.  
\end{Algorithm}

\begin{Remark}[Tensor-free eigenvalue thresholding]
  \label{rem:eigen-threshold}
  Using the relation between eigenvalues and singular values, we
  apply \thref{alg:subspace-iteration} and \ref{alg:aug-lanc-bidiag}
  in the same manner to compute the positive eigenvalue thresholding.
  More precisely, with
  $\lambda_m \coloneqq \sigma_m \iProd{\widetilde{\Vek
      u}_m}{\widetilde{\Vek v}_m}$
  and
  $\Mat \Lambda_n \coloneqq \diag(\lambda_0, \dots,
  \lambda_{\ell-1})$, we obtain the eigenvalue decomposition from the
  singular value decomposition, \ie{},
  \begin{equation*}
    \widetilde{\Mat U}_n^* \Mat H \Vek w \Mat H \widetilde{\Mat U}_n
    = \Mat \Lambda_n
    \qquad\text{from}\qquad
    \widetilde{\Mat V}_n^* \Mat H \Vek w \Mat H \widetilde{\Mat U}_n
    = \Mat \Sigma_n.
  \end{equation*}
  Analogously, we can transfer \thref{alg:tf-threshold:subspace} and
  \ref{alg:tf-threshold:lanczos} to the quadratic setting.  \qed
\end{Remark}

Besides the singular value thresholding, the proximal methods in
\autoref{sec:sing-value-thresh} to solve the lifted and relaxed
bilinear and quadratic problems in \autoref{sec:bilin-inverse-probl}
require the application of the lifted operators $\breve{\mathcal B}$
and $\breve{\mathcal Q}$ as well as their adjoints
$\breve{\mathcal B}^*$ and $\breve{\mathcal Q}^*$.  Both operations
can be computed in a tensor-free manner.  Assuming that $\Vek w$
has a low rank, one may compute the lifted bilinear forward operator with the
aid of the universal property in \thref{the:bilin-lifting}.

\begin{Corollary}[Tensor-free bilinear lifting]
  \label{cor:tf-for-oper:bilin}
  Let $\mathcal B \colon \mathcal H_1 \times \mathcal H_2 \to \mathcal
  K$ be a bilinear mapping.
  If $\Vek w \in \mathcal H_1 \otimes \mathcal H_2$ has the
  representation $\Vek w = \widetilde{\Mat V} \Mat \Sigma
  \widetilde{\Mat U}^*$ with $\widetilde{\Mat U} \coloneqq
  [\widetilde{\Vek u}_0, \dots, \widetilde{\Vek u}_{\ell -1}]$, $\Mat
  \Sigma \coloneqq \diag(\sigma_0, \dots, \sigma_{\ell-1})$, and $\widetilde{\Mat V} \coloneqq
  [\widetilde{\Vek v}_0, \dots, \widetilde{\Vek v}_{\ell -1}]$, then
  the lifted forward operator $\breve{\mathcal B}$ acts by
  \begin{equation*}
    \breve{\mathcal B}(\Vek w) = \sum_{n=0}^{\ell-1} \sigma_n \,
    \mathcal B(\widetilde{\Vek u}_n, \widetilde{\Vek v}_n).
  \end{equation*}
\end{Corollary}

Likewise, one may apply the lifted quadratic forward operator by using
\thref{cor:quad-lifting}.

\begin{Corollary}[Tensor-free quadratic lifting]
  \label{cor:tf-for-oper:quad}
  Let $\mathcal Q \colon \mathcal H \to \mathcal K$ be a quadratic
  mapping.  If $\Vek w \in \mathcal H \otimes_{\mathrm{sym}} \mathcal H$ has the
  representation $\Vek w = \widetilde{\Mat U} \Mat \Lambda
  \widetilde{\Mat U}^*$ with $\widetilde{\Mat U} \coloneqq
  [\widetilde{\Vek u}_0, \dots, \widetilde{\Vek u}_{\ell -1}]$, $\Mat
  \Lambda \coloneqq \diag(\lambda_0, \dots, \lambda_{\ell-1})$, then
  the lifted forward operator $\breve{\mathcal Q}$ acts by
  \begin{equation*}
    \breve{\mathcal Q}(\Vek w) = \sum_{n=0}^{\ell-1} \lambda_n \,
    \mathcal Q(\widetilde{\Vek u}_n).
  \end{equation*}
\end{Corollary}

Considering the proximal methods, we see that the adjoint lifting only
occurs in the argument of the singular value thresholding.  If one
applies the subspace iteration or the augmented \PLanczos\ process, it
is hence enough to study the left-hand and right-hand actions of the
adjoint liftings.  In the bilinear setting, these actions can be
expressed by the \emph{left-hand} or \emph{right-hand adjoint} of the
original bilinear mapping $\mathcal B$.

\begin{Lemma}[Tensor-free adjoint bilinear lifting]
  \label{lem:tens-free-adj-lift:bilin}
  Let
  $\mathcal B \colon \mathcal H_1 \times \mathcal H_2 \to \mathcal K$
  be a bilinear mapping.  The left-hand and right-hand actions of the
  adjoint lifting
  $\breve{\mathcal B}^*(\Vek y) \in \mathcal H_1 \otimes
  \mathcal H_2$ with $\Vek y \in \mathcal K$ are given by
  \begin{equation*}
    \breve{\mathcal B}^*(\Vek y) \, \Mat H_1 \Vek e
    =
    [\mathcal B(\Vek e, \cdot)]^*(\Vek y)
    \qquad\text{and}\qquad
    [\breve{\mathcal B}^*(\Vek y)]^* \, \Mat H_2 \Vek f
    =
    [\mathcal B(\cdot, \Vek f)]^*(\Vek y)
  \end{equation*}
  for $\Vek e \in \mathcal H_1$ and $\Vek f \in \mathcal H_2$. 
\end{Lemma}

\begin{Proof}
  Testing the right-hand action of the image
  $\breve{\mathcal B}^*(\Vek y)$ on $\Vek e \in \mathcal H_1$ with an
  arbitrary vector $\Vek f \in \mathcal H_2$, we obtain
  \begin{align*}
    \iProdb{\breve{\mathcal B}^*(\Vek y) \, \Mat H_1 \Vek
    e}{\Vek f}_{\mathcal H_2}
    &= 
      \tr \bigl(\Vek f^*  \Mat H_2\, \breve{\mathcal B}^*(\Vek y)
      \, \Mat H_1 \Vek e \bigr)
      =
      \tr \bigl( \Vek e \Vek f^* \Mat H_2 \, \breve{\mathcal B}^*(\Vek
      y) \, \Mat H_1 \bigr)
    \\[\fskip]
    &=
      \iProdb{\breve{\mathcal B}^*(\Vek y)}{\Vek e
      \otimes \Vek f}_{\mathcal H_1 \otimes \mathcal H_2}
      =
      \iProdb{\Vek y}{\mathcal B(\Vek e, \Vek f)}_{\mathcal K}
      =
      \iProdb{[\mathcal B(\Vek e, \cdot)]^*(\Vek y)}{\Vek f}_{\mathcal
      H_2}.
  \end{align*}
  The left-hand action follows analogously.  \qed
\end{Proof}

An similar observation holds for the quadratic setting, where the
adjoint is taken with respect to the symmetric subspace $\mathcal H
\otimes_{\mathrm{sym}} \mathcal H$.  The associate
bilinear mapping to $\mathcal Q$ is again denoted by $\mathcal
B_{\mathcal Q}$.

\begin{Lemma}[Tensor-free adjoint quadratic lifting]
  \label{lem:tens-free-adj-lift:quad}
  Let $\mathcal Q \colon \mathcal H \to \mathcal K$ denote a quadratic
  mapping.  The action of the (symmetric) adjoint
  lifting
  $\breve{\mathcal Q}^*(\Vek y) \in \mathcal H \otimes_{\mathrm{sym}}
  \mathcal H$ with $\Vek y \in \mathcal K$ is given by
  \begin{equation*}
    \breve{\mathcal Q}^*(\Vek y) \, \Mat H \Vek e
    =
    \tfrac12 \,
    [\mathcal B_{\mathcal Q}(\Vek e, \cdot)]^*(\Vek y)
    + \tfrac12 \,
    [\mathcal B_{\mathcal Q}(\cdot, \Vek e)]^*(\Vek y)
    \submathskip
  \end{equation*}
  for $\Vek e \in \mathcal H$. 
\end{Lemma}

\begin{Proof}
  Similarly to before, we test the action of the image
  $\breve{\mathcal Q}^*(\Vek y)$ on $\Vek u \in \mathcal H$ with an
  arbitrary vector $\Vek v \in \mathcal H$.  Exploiting the symmetry,
  we obtain
  \begin{align*}
    \iProdb{\breve{\mathcal Q}^*(\Vek y) \, \Mat H \Vek
    e}{\Vek f}_{\mathcal H}
    &= 
      \tfrac12  \tr \bigl(\Vek f^*  \Mat H\, \breve{\mathcal Q}^*(\Vek y)
      \, \Mat H \Vek e \bigr)
      + \tfrac12  \tr \bigl(\Vek e^*  \Mat H\, \breve{\mathcal Q}^*(\Vek y)
      \, \Mat H \Vek f \bigr)
    \\[\fskip]
    &=
      \iProdb{\breve{\mathcal Q}^*(\Vek y)}{\tfrac12 \, (\Vek e
      \otimes \Vek f) + \tfrac12 \, (\Vek f \otimes \Vek e)}_{\mathcal
      H \otimes_{\mathrm{sym}} \mathcal H}
    \\[\fskip]
    &=
      \iProdb{\Vek y}{\tfrac12 \, [\mathcal Q(\Vek e +
      \Vek f) - \mathcal Q(\Vek e) - \mathcal Q(\Vek
      f )]}_{\mathcal K}
    \\[\fskip]
    &=
      \iProdb{\Vek y}{\tfrac12 \,[\mathcal B_{\mathcal Q}(\Vek e, \Vek
      f) + \mathcal B_{\mathcal Q}(\Vek f, \Vek e)]}_{\mathcal K}
    \\[\fskip]
    & =
      \tfrac12 \,
      \iProdb{[\mathcal B_{\mathcal Q}(\Vek e, \cdot)]^*(\Vek y)}{\Vek f}_{\mathcal
      H}
      +\tfrac12 \,
      \iProdb{[\mathcal B_{\mathcal Q}(\cdot, \Vek e)]^*(\Vek y)}{\Vek f}_{\mathcal
      H}. \tag*{\qed}
  \end{align*}
\end{Proof}

\begin{Remark}[Composed tensor-free adjoint lifting]
  \label{rem:com-tf-adj-lift}
  Since the left-hand and right-hand actions of the tensor
  $\Vek w^{(n)} = \sum_{k=0}^{R-1} \sigma_k^{(n)} \, (\widetilde{\Vek
    u}_k^{(n)} \otimes \widetilde{\Vek v}_k^{(n)})$
  are simply given by
  \begin{align}
    \label{eq:tensor-right-act}
    \Vek w^{(n)} \Mat H_1 \Vek e 
    &=
    \sum_{k=0}^{R-1} \sigma_k^{(n)} \, \iProdn{\Vek e}{\widetilde{\Vek
        u}_k^{(n)}}_{\mathcal H_1}
    \, \widetilde{\Vek v}_k^{(n)}
    \shortintertext{and}
    \label{eq:tensor-left-act}
    (\Vek w^{(n)})^* \Mat H_2 \Vek f
    &=
    \sum_{k=0}^{R-1} \sigma_k^{(n)} \, \iProdn{\Vek f}{\widetilde{\Vek
        v}_k^{(n)}}_{\mathcal H_2}
    \, \widetilde{\Vek u}_k^{(n)},
  \end{align}
  the right-hand action of the singular value thresholding argument
  $\Vek w = \Vek w^{(n)} - \tau \, \breve{\mathcal B}^*(\Vek
  y^{(n+1)})$ within the proximal methods in
  Section~\ref{sec:sing-value-thresh} is given by
  \begin{equation}
    \label{eq:sing-thres:right-act}
    \Vek w \Mat H_1 \Vek e 
    = 
    - \tau \, [\mathcal B(\Vek e, \cdot)]^*\bigl(\Vek y^{(n+1)}\bigr)
    + \sum_{k=0}^{R-1} \sigma_k^{(n)} \, \iProdb{\Vek e}{\widetilde{\Vek
        u}_k^{(n)}}_{\mathcal H_1}
    \, \widetilde{\Vek v}_k^{(n)}
  \end{equation}
  and the left-hand action by
  \begin{equation}
    \label{eq:sing-thres:left-act}
    \Vek w^* \Mat H_2 \Vek f 
    = 
    - \tau \, [\mathcal B( \cdot, \Vek f)]^*\bigl(\Vek y^{(n+1)}\bigr)
    + \sum_{k=0}^{R-1} \sigma_k^{(n)} \, \iProdb{\Vek f}{\widetilde{\Vek
        v}_k^{(n)}}_{\mathcal H_2}
    \, \widetilde{\Vek u}_k^{(n)},
  \end{equation}
  where $\Vek w^{(n)} = \sum_{k=0}^{R-1} \sigma_k^{(n)} \,
  (\widetilde{\Vek u}_k^{(n)} \otimes \widetilde{\Vek v}_k^{(n)})$.
  Analogously, the actions in the quadratic setting are given by
  \begin{equation}
    \label{eq:sing-thres:quad-act}
    \Vek w \Mat H \Vek e 
    = 
    - \tau \, [\mathcal B_{\mathcal Q}(\Vek e, \cdot)]^*\bigl(\Vek y^{(n+1)}\bigr)
    + \sum_{k=0}^{R-1} \lambda_k^{(n)} \, \iProdb{\Vek e}{\widetilde{\Vek
        u}_k^{(n)}}_{\mathcal H_1}
    \, \widetilde{\Vek u}_k^{(n)}
  \end{equation}
  where $\Vek w^{(n)} = \sum_{k=0}^{R-1} \lambda_k^{(n)} \,
  (\widetilde{\Vek u}_k^{(n)} \otimes \widetilde{\Vek u}_k^{(n)})$.  \qed
\end{Remark}

Now we are ready to rewrite the proximal methods in
\autoref{sec:sing-value-thresh} into tensor-free variants.
Exemplarily, we consider the primal-dual method for bilinear operators
and exact data, see \thref{alg:exact-data}.

\begin{Algorithm}[Tensor-free primal-dual for exact data]
  \label{alg:tensor-free}
  \begin{enumerate}[(i)]
  \item Initiation: Fix the parameters $\tau, \sigma > 0$
    and $\theta \in [0,1]$.  Choose the start value $(\Vek w^{(0)},
    \Vek y^{(0)}) = (\Vek 0 \otimes \Vek 0, \Vek 0)$ in $(\mathcal H_1
    \otimes \mathcal H_2) \times \mathcal K$, and set
    $\Vek w^{(-1)}$ to $\Vek w^{(0)}$.
  \item Iteration: For $n \ge 0$, update $\Vek w^{(n)}$ and $\Vek
    y^{(n)}$:
    \begin{enumerate}[(a)]
    \item Using the tensor-free computations in
      \thref{cor:tf-for-oper:bilin}, determine
      \begin{align*}
        \Vek y^{(n+1)}
        \coloneqq \Vek y^{(n)} + \sigma \, \bigl( (1 + \theta) \; 
        \breve{\mathcal B}(\Vek w^{(n)})
          -  \theta \; \breve{\mathcal B}(\Vek w^{(n-1)}) - \Vek g^\dagger \bigr).
      \end{align*}
    \item Compute a low-rank representation
      $\Vek w^{(n+1)} = \widetilde{\Mat V}^{(n+1)} \Mat \Sigma^{(n+1)}
      \widetilde{\Mat U}^{(n+1)}$
      of the singular value threshold
      \begin{equation*}
        \mathcal S_\tau ( \Vek w^{(n)} - \tau \, \breve{\mathcal
          B}^*(\Vek y^{(n+1)}))
      \end{equation*}
      with \thref{alg:tf-threshold:lanczos} (or
      \ref{alg:tf-threshold:subspace}).  The required actions are
      given in \eqref{eq:sing-thres:right-act} and
      \eqref{eq:sing-thres:left-act}.
    \end{enumerate}
  \end{enumerate}
\end{Algorithm}

\begin{Remark}
  \label{rem:tensor-free:1}
  As starting value for the augmented \PLanczos\ bidiagonalization
  according to Algorithm~\ref{alg:tf-threshold:lanczos} required for step
  (ii.b) of Algorithm~\ref{alg:tensor-free}, we suggest a linear
  combination of the right-hand singular vectors of the previous
  iteration $\Vek w^{(n)}$ in the hope that they are good
  approximations of the new singular vectors.  \qed
\end{Remark}

Adapting the computation of $\Vek y^{(n+1)}$, one may analogously
apply \thref{alg:tikh-reg} and \ref{alg:inexact-data} in a complete
tensor-free manner.  Using \thref{cor:tf-for-oper:quad},
\thref{rem:eigen-threshold}, and
Equation~\eqref{eq:sing-thres:quad-act}, we obtain the corresponding
tensor-free proximal methods for quadratic inverse problems.

Because the singular value thresholding can be computed with arbitrary
high accuracy, the convergence results for the primal-dual algorithm
translates to our setting.  The convergence analysis
\cite[Theorem~1]{CP11} yields the following convergence guarantee, where
the norm of the bilinear operator $\mathcal B$ is defined by
\begin{equation*}
  \pNormn{\mathcal B} 
  \coloneqq 
  \sup_{\Vek u \in \mathcal H_1 \setminus \{ \Vek 0\}} \,
  \sup_{\Vek v \in \mathcal H_2 \setminus \{ \Vek 0\}} \,
  \frac{\pNormn{\mathcal B(\Vek u, \Vek v)}_{\mathcal K}}{\pNormn{\Vek
    u}_{\mathcal H_1} \, \pNormn{\Vek v}_{\mathcal H_2}}.
\end{equation*}

\begin{Theorem}[Convergence]
  \label{the:convergence}
  Under \thref{ass:basic-pre} and the parameter choice rule
  $\theta = 1$ as well as $\tau \sigma \pNormn{\mathcal B}^2 < 1$, the
  iteration $(\Vek w^{(n)}, \Vek y^{(n)})$ in \thref{alg:tensor-free}
  converges to a point $(\Vek w^\dagger, \Vek y^\dagger)$, where
  $\Vek w^\dagger = \Vek u^\dagger \otimes \Vek v^\dagger$ corresponds
  to a solution $(\Vek u^\dagger, \Vek v^\dagger)$ of the bilinear
  inverse problem \eqref{eq:gen-bilin-prob}.
\end{Theorem}

\begin{Proof}
  For the general minimization problem \eqref{eq:gen-prob-struc}, the
  related saddle-point problem is given by
  \begin{equation}
    \label{eq:gen-primal-dual}
    \minimize_{\Vek w \in \mathcal H_1 \otimes \mathcal
      H_2}
    \quad
    \maximize_{\Vek y \in \mathcal K}
    \quad
    \iProdn{\mathcal A(\Vek w)}{\Vek y} + G(\Vek w) - F^*(\Vek y),
    \addmathskip
  \end{equation}
  \cf\ \cite{CP11}.  Hence, the bilinear relaxation with exact data
  \eqref{eq:nucl-min-bilin} corresponds to the primal-dual formulation
  \begin{equation}
    \label{eq:exact-primal-dual}
    \minimize_{\Vek w \in \mathcal H_1 \otimes \mathcal
      H_2} 
    \quad 
    \maximize_{\Vek y \in \mathcal K}
    \quad
    \iProdn{\breve{\mathcal B}(\Vek w) - \Vek g^\dagger}{\Vek y} +
    \pNormn{\Vek w}_{\mathcal H_1 \otimes_\uppi \mathcal H_2}.
  \end{equation}
  Due to \cite[Theorem~28.3]{Roc70}, the first components
  $\widetilde{\Vek w}$ of the saddle-points
  $(\widetilde{\Vek w}, \widetilde{\Vek y})$ of
  \eqref{eq:exact-primal-dual} are solutions of
  \eqref{eq:nucl-min-bilin}.  Vice versa,
  \cite[Corollary~28.2.2]{Roc70} implies that the solutions
  $\widetilde{\Vek w}$ of \eqref{eq:nucl-min-bilin} are saddle-points
  of \eqref{eq:exact-primal-dual}.  In particular, the saddle-point
  problem \eqref{eq:exact-primal-dual} has at least one solution since
  the given data are exact.

  Now, \cite[Theorem~1]{CP11} yields the convergence
  $(\Vek w^{(n)}, \Vek y^{(n)}) \to (\Vek w^\dagger, \Vek y^\dagger)$
  of the primal-dual iteration in \thref{alg:tensor-free}, where the
  limit $(\Vek w^\dagger, \Vek y^\dagger)$ denotes a saddle-point of
  \eqref{eq:exact-primal-dual}, and $\Vek w^\dagger$ thus a solution
  of \eqref{eq:nucl-min-bilin}.  Under \thref{ass:basic-pre}, this
  solution has at most rank one; so every
  $(\Vek u^\dagger, \Vek v^\dagger)$ with
  $\Vek u^\dagger \otimes \Vek v^\dagger = \Vek w^\dagger$ is a
  solution of the bilinear problem \eqref{eq:gen-bilin-prob}.  \qed
\end{Proof}

Similar convergence guarantees can be obtained for the bilinear
relaxations \eqref{eq:inexact-data-bilin} and
\eqref{eq:tikh-func-bilin} as well as for the quadratic variants
\eqref{eq:nucl-min-quad}, \eqref{eq:inexact-data-quad}, and
\eqref{eq:tikh-func-quad}.  Depending on the considered problem -- the
bilinear or quadratic forward operator -- and on the applied proximal
algorithm, one may even obtain explicit convergence rates.

\section{Reducing rank by Hilbert space reweighting}
\label{sec:reduc-rank-adapt}

The projective norm heuristic usually ensures that the solutions of
the relaxed lifted minimization problems in
\autoref{sec:bilin-inverse-probl} have low rank, but how may we
decrease the rank of the iteration $\Vek w^{(n)}$ even further to
speed up the convergence, and how can we obtain meaningful solutions
if \thref{ass:basic-pre} is not fulfilled?

If we look back at the employed methods in compressed sensing, where
one like to recover a sparse vector from a set of linear measurements,
one idea to improve the sparsity of the solution is to reweight the
$\ell_1$-norm.  Since the nuclear norm coincides with the
$\ell_1$-norm of the singular values, one can try to incorporate this
reweighting approach, which however is not possible.  The main reason
here is that the singular values do not have a special ordering.  If a
tensor, for instance, possesses a multiple singular value, then the
reweighting will be ambiguous.

Instead of reweighting the projective norm itself, we propose to
modify the \PHilbert\ norms of the underlying spaces $\mathcal H_1$
and $\mathcal H_2$, where we initially consider the bilinear setting.
From a heuristic point of view, we may interpret the leading
singular value $\sigma_0$ together with corresponding singular vectors
$\widetilde{\Vek u}_0 \otimes \widetilde{\Vek v}_0$ of the variable
$\Vek w^{(n)}$ as an approximate solution of the inverse problem
\eqref{eq:gen-bilin-prob}; so we may adapt the inner products of
$\mathcal H_1$ and $\mathcal H_2$ in a way that promote vectors in
this direction.

More generally, we want to reweight the norms of $\mathcal H_1$ and
$\mathcal H_2$ with respect to some orthonormal bases, which, for
instance, result from the singular value decomposition of the primal
variable $\Vek w^{(n)}$.  In the following, we only consider the
reweighting of $\mathcal H_1$.  The reweighting of $\mathcal H_2$
can be done completely analogously.  If
$\Mat \Phi \coloneqq [\Vek \phi_0, \dots, \Vek \phi_{N_1-1}]$ denotes an
arbitrary orthonormal basis of $\mathcal H_1$, then \PParseval's
identity states
\begin{equation}
  \label{eq:parseval-H1}
  \pNormn{\Vek u}_{\mathcal H_1}^2 
  = 
  \sum_{n=0}^{N_1-1} \iProd{\Vek u}{\Vek \phi_n}_{\mathcal H_1}^2 
  =
  (\Vek u^* \Mat H_1 \Mat \Phi) \, (\Mat \Phi^* \Mat H_1 \Vek u).
\end{equation}
In other words, the matrix $\Mat H_1$ corresponding to the inner
product on $\mathcal H_1$ can be written in the form
$\Mat H_1 = \Mat H_1 \Mat \Phi \Mat \Phi^* \Mat H_1$, which incidentally
shows $\Mat \Phi \Mat \Phi^* = \Mat H_1^{-1}$.  To reweight the
$\mathcal H_1$-norm \eqref{eq:parseval-H1} with respect to the basis
$\Mat \Phi$, we introduce the weights
$\Mat \Xi \coloneqq \diag(\xi_0, \dots, \xi_{N_1-1})$ and the adapted
norm $\pNormn{\cdot}_{\mathcal H_1(\Mat \Xi)}$ defined by
\begin{equation}
  \label{eq:parselval-H1-rew}
  \pNormn{\Vek u}_{\mathcal H_1(\Mat \Xi)}^2
  =
  \sum_{n=0}^{N_1-1} \xi_n \, \iProdn{\Vek u}{\Vek \phi_n}_{\mathcal
    H_1}^2
  =
  (\Vek u^* \Mat H_1 \Mat \Phi) \, \Mat \Xi \, (\Mat \Phi^* \Mat H_1 \Vek u).
\end{equation}
In so doing, we obtain the updated inner product matrix
$\Mat H_1(\Mat \Xi) = \Mat H_1 \Mat \Phi \Mat \Xi \Mat \Phi^* \Mat H_1$ with
inverse $ \Mat H_1^{-1}(\Mat \Xi) = \Mat \Phi \Mat \Xi^{-1} \Mat \Phi^*$.
Depending on the weights, some directions are more promoted or
penalized than others.

Unfortunately, for large-scale bilinear inverse problems, the proposed
approach is impractical since we have to store a complete orthonormal
basis.  To overcome this disadvantages, we make a slightly refinement
of our approach, where we only reweight a small set of basis vectors
that we want to promote; so we choose the weights $\Mat \Xi$ as
$\xi_n = 1 - \lambda_n$ with $\lambda_n \in (0,1)$ for
$n=0, \dots, S-1$ and $\lambda_n = 0$ otherwise, and thus make the
approach
\begin{equation}
  \label{eq:re-weight-H1}
  \Mat H_1 (\Vek \Xi)
  \coloneqq
  \Mat H_1 - \sum_{n=0}^{S-1} \lambda_n \, \Mat H_1 \Vek \phi_n \Vek \phi_n^*
  \Mat H_1
  =
  \Mat H_1 - \sum_{n=0}^{S-1} \lambda_n \widetilde{\Vek \phi}_n \widetilde{\Vek \phi}_n^*
    \addmathskip
\end{equation}
with $\widetilde{\Vek \phi}_n \coloneqq \Mat H_1 \Vek \phi_n$.  The
inverse is here given by
\begin{equation}
  \label{eq:re-weight-H1-inv}
  \Mat H_1^{-1}(\Mat \Xi)
  =
  \sum_{n=0}^{N_1 - 1} \frac{1}{1 - \lambda_n} \, \Vek \phi_n \Vek \phi_n^*
  =
  \Mat H_1^{-1} - \sum_{n=0}^{S-1} \Bigl( 1 - \frac{1}{1 - \lambda_n}
  \Bigr) \, \Vek \phi_n \Vek \phi_n^*.
\end{equation}
Hence, to update the inner product matrices, we only require the
original matrices $\Mat H_1$ and $\Mat H_1^{-1}$, the (transformed)
promoted vectors $\Vek \phi_n$ and $\widetilde{\Vek \phi}_n$, and the
weights $\lambda_n$.  For the second \PHilbert\ space $\mathcal H_2$,
we proceed completely analogously. 

The reweighting of the \pers{Hilbert} spaces has consequences for the
proposed algorithms.  On the one hand, notice that the adjoint of the
lifted operator $\breve{\mathcal B}$ directly depends on the actual
\PHilbert\ spaces $\mathcal H_1$ and $\mathcal H_2$.  In order to
update the adjoint, we first determine the standardized adjoint
$\breve{\mathcal B}^*_{\mathcal{H\!S, E}}$ with respect to the
\PHilbert-\PSchmidt\ and \PEuclid{ian} inner product.  Afterwards, we
transform this adjoint to the actual spaces by the following lemma.

\begin{Lemma}[Adjoint operator]
  \label{lem:adj-oper}
  The adjoint operator
  $\breve{\mathcal B}^*_{\mathcal H_1 \otimes \mathcal H_2,
    \mathcal K}$
  with respect to the \PHilbert\ spaces
  $\mathcal H_1 \otimes \mathcal H_2$ and $\mathcal K$ is
  given by
  \begin{equation*}
    \breve{\mathcal B}^*_{\mathcal H_1 \otimes \mathcal H_2,
      \mathcal K} (\Vek y)
    = \Mat H_2^{-1} \,
    \breve{\mathcal B}^*_{\mathcal{H\!S}, \mathcal E}(\Mat K \Vek y) \,
    \Mat H_1^{-1},
  \end{equation*}
  where $\breve{\mathcal B}^*_{\mathcal{H\!S}, \mathcal E}$ denotes the adjoint
  with respect to \PHilbert--\PSchmidt\ and \PEuclid{ian} inner
  product.
\end{Lemma}

\begin{Proof}
  The assertion immediately follows from
  \begin{equation*}
    \iProd{\breve{\mathcal B}(\Vek w)}{\Vek y}_{\mathcal K}
    = \iProd{\breve{\mathcal B}(\Vek w)}{\Mat K \Vek y}
    = \iProd{\Vek w}{\breve{\mathcal B}^*_{\mathcal{H\!S}, \mathcal E}(\Mat K \Vek
      y)}
    = \iProd{\Vek w}{\Mat H_2^{-1} \,
    \breve{\mathcal B}^*_{\mathcal{H\!S}, \mathcal E}(\Mat K \Vek y) \,
    \Mat H_1^{-1}}_{\mathcal H_1 \otimes \mathcal H_2}
  \end{equation*}
  for all $\Vek w \in \mathcal H_1 \otimes \mathcal H_2$
  and all $\Vek y \in \mathcal K$. \qed
\end{Proof}

On the one side, \thref{lem:adj-oper} allows us to transform the
adjoint $\breve{\mathcal B}^*_{\mathcal{H\!S}, \mathcal E}$, which can
usually be determined more easily, to the \PHilbert\ spaces
$\mathcal H_1$, $\mathcal H_2$, and $\mathcal K$.  On the other side
and more generally, we may transform the adjoint lifted operator
between arbitrary \PHilbert\ space structures.  For our specific
setting in \eqref{eq:re-weight-H1} and \eqref{eq:re-weight-H1-inv},
for instance, we obtain the following transformation rule, where
$\Mat T_1(\Mat \Xi)$ denotes the transformation matrix
\begin{equation}
  \label{eq:weight-trans}
  \Mat T_1(\Mat \Xi) 
  \coloneqq
  \Mat H_1 \Mat H_1^{-1}(\Mat \Xi)
  =
  \Mat I - \sum_{n=0}^{S-1} \Bigl( 1 - \frac{1}{1 - \lambda_n}
  \Bigr) \, \widetilde{\Vek \phi}_n \Vek \phi_n^*
\end{equation}
with the identity $\Mat I$.  The transformation $\Mat T_2(\Mat \Xi)$
is defined analogously.

\begin{Corollary}[Adjoint operator]
  \label{cor:adj-oper}
  The adjoint operator
  $\breve{\mathcal B}^*_{\mathcal H_1(\Mat \Xi) \otimes
    \mathcal H_2 (\Mat \Xi),
    \mathcal K}$
  with respect to the \PHilbert\ spaces
  $\mathcal H_1(\Mat \Xi) \otimes \mathcal H_2(\Mat \Xi)$ and $\mathcal K$ is
  given by
  \begin{equation*}
    \breve{\mathcal B}^*_{\mathcal H_1(\Mat \Xi) \otimes
      \mathcal H_2(\Mat \Xi), 
      \mathcal K} (\Vek y)
    =
    \Mat T_2^*(\Mat \Xi) \;
    \breve{\mathcal B}^*_{\mathcal H_1 \otimes \mathcal H_2,
      \mathcal K} (\Vek y) \,
    \Mat T_1(\Mat \Xi).
  \end{equation*}
\end{Corollary}

\begin{Proof}
  Apply \thref{lem:adj-oper} two times to transform the adjoint
  firstly from
  $(\mathcal H_1 \otimes \mathcal H_2, \mathcal K)$ to
  $(\mathcal{H\!S}, \mathcal E)$ and, secondly, from
  $(\mathcal{H\!S}, \mathcal E)$ to
  $(\mathcal H_1(\Mat \Xi) \otimes \mathcal H_2(\Mat
  \Xi), \mathcal K)$.  \qed
\end{Proof}

The main benefit of \thref{cor:adj-oper} compared with
\thref{lem:adj-oper} is that the transformation can be done without
involving inverse matrices.  On the other side, one relies on an
efficient and direct implementation of
$\breve{\mathcal B}^*_{\mathcal H_1 \otimes \mathcal H_2,
  \mathcal K}$
for the unweighted spaces.  Remember that this adjoint may be
determined by \thref{lem:tens-free-adj-lift:bilin}.  The required
actions to compute the singular value threshold regarding the
reweighted spaces then have the following form.

\begin{Corollary}[Tensor-free adjoint bilinear lifting]
  \label{cor:tens-free-adj-bilin-lift}
  The left-hand and right-hand actions of the reweighted adjoint
  $\breve{\mathcal B}^*_{\mathcal H_1(\Mat \Xi) \otimes \mathcal
    H_2 (\Mat \Xi), \mathcal K}( \Vek y)$
  with $\Vek y \in \mathcal K$ are given by
  \begin{align*}
    \breve{\mathcal B}^*_{\mathcal H_1(\Mat \Xi) \otimes \mathcal
    H_2(\Mat \Xi), \mathcal K} (\Vek y) \, \Mat H_1(\Mat \Xi) \, \Vek u
    &=
      \Mat T_2^*(\Mat \Xi) \,
      [\mathcal B(\Vek e, \cdot)]^*_{\mathcal H_2, \mathcal K}(\Vek y)
    \\
    \shortintertext{and}
    [\breve{\mathcal B}^*_{\mathcal H_1(\Mat \Xi) \otimes \mathcal
    H_2(\Mat \Xi), \mathcal K}(\Vek y)]^* \, \Mat H_2(\Mat \Xi) \, \Vek v
    &=
      \Mat T_1^*(\Mat \Xi) \,
      [\mathcal B(\cdot, \Vek f)]^*_{\mathcal H_1, \mathcal K}(\Vek y)
  \end{align*}
  for $\Vek e \in \mathcal H_1(\Mat \Xi)$ and $\Vek f \in \mathcal H_2(\Mat \Xi)$.
\end{Corollary}

\begin{Proof}
  The assertion follows from \thref{cor:adj-oper} and
  Equation~\eqref{eq:weight-trans} by
  \begin{align*}
    \breve{\mathcal B}^*_{\mathcal H_1(\Mat \Xi) \otimes \mathcal
    H_2(\Mat \Xi), \mathcal K} (\Vek y) \, \Mat H_1(\Mat \Xi) \,
    \Vek u
    &=
      \Mat T_2^*(\Mat \Xi) \;
      \breve{\mathcal B}^*_{\mathcal H_1 \otimes \mathcal
      H_2, \mathcal K} (\Vek y) \, \Mat T_1(\Mat \Xi) \, \Mat H_1(\Mat
      \Xi) \, \Vek u 
    \\[\fskip]
    &= 
      \Mat T_2^*(\Mat \Xi) \;
      \breve{\mathcal B}^*_{\mathcal H_1 \otimes \mathcal
      H_2, \mathcal K} (\Vek y)  \, \Mat H_1 \, \Vek u.
  \end{align*}
  The second identity follows completely analogously.  \qed
\end{Proof}

Thanks to the above transformations, we can expand the tensor-free
primal-dual iteration for exact data in \thref{alg:tensor-free} by an
efficient reweighting step, which we perform every
$n_{\mathrm{rew}}$-th iteration.  If the set of promoted directions
$\Mat \Phi$ and $\Mat \Psi$ for $\mathcal H_1$ and $\mathcal H_2$ is
empty, we perform the unweighted algorithm with $\Mat \Xi = \Mat I$
and hence $\Mat H_1(\Mat \Xi) = \Mat H_1$ as well as
$\Mat H_2(\Mat \Xi) = \Mat H_2$.

In order to avoid a recursive reweighting, we always reweight the
original \PHilbert\ spaces.  Thus, if we want to promote the leading
singular vectors of $\Vek w^{(n)}$, we first have to compute the
singular value decomposition
$\Vek w^{(n)} = \sum_{k=0}^{S-1} \sigma_k' \, (\widetilde{\Vek u}_k'
\otimes \widetilde{\Vek v}_k')$ with respect to underlying, original
\PHilbert\ spaces $\mathcal H_1$ and $\mathcal H_2$.  Based on the
tensor-free characterization
$\Vek w^{(n)} = \sum_{k=0}^{R-1} \sigma_k^{(n)} \, (\widetilde{\Vek
  u}_k^{(n)} \otimes \widetilde{\Vek v}_k^{(n)})$ obtained form the
(weighted) singular value thresholding, this decomposition can be
computed by \thref{alg:subspace-iteration} or
\thref{alg:aug-lanc-bidiag}, where the required actions are simply
given by \eqref{eq:tensor-right-act} and \eqref{eq:tensor-left-act}.
Note that the involved inner products are here again the unweighted
versions.

After a reweighting step, the singular value thresholding must be
computed with respect to the new weight.  For this purpose, we adapt the
actions in \eqref{eq:sing-thres:right-act} and
\eqref{eq:sing-thres:left-act} by
\thref{cor:tens-free-adj-bilin-lift}.  In so doing, for
$\Vek w = \Vek w^{(n)} - \tau \, \breve{\mathcal B}^*_{\mathcal
  H_1(\Mat \Xi) \otimes \mathcal H_2(\Mat \Xi), \mathcal
  K}(\Vek y^{(n+1)})$, we obtain the right-hand action
\begin{equation}
  \label{eq:sing-thres:right-act:rew}
  \Vek w \Mat H_1 (\Mat \Xi) \, \Vek e 
  = 
  - \tau \, \Mat T_2^*(\Mat \Xi) \,
  [\mathcal B(\Vek e, \cdot)]^*_{\mathcal H_2, \mathcal K}
  \bigl(\Vek y^{(n+1)}\bigr)
  + \sum_{k=0}^{R-1} \sigma_k^{(n)} \, \iProdb{\Vek e}{\widetilde{\Vek
      u}_k^{(n)}}_{\mathcal H_1(\Mat \Xi)}
  \, \widetilde{\Vek v}_k^{(n)}
\end{equation}
and the left-hand action
\begin{equation}
  \label{eq:sing-thres:left-act:rew}
  \Vek w^* \Mat H_2(\Mat \Xi) \, \Vek f 
  = 
  - \tau \, \Mat T_1^*(\Mat \Xi) \,
  [\mathcal B( \cdot, \Vek f)]^*_{\mathcal H_1, \mathcal K}
  \bigl(\Vek y^{(n+1)}\bigr)
  + \sum_{k=0}^{R-1} \sigma_k^{(n)} \, \iProdb{\Vek f}{\widetilde{\Vek
      v}_k^{(n)}}_{\mathcal H_2(\Mat \Xi)}
  \, \widetilde{\Vek u}_k^{(n)}.
\end{equation}
The definition of the transformations $\Mat T_1(\Mat \Xi)$ and
$\Mat T_2(\Mat \Xi)$ is given in \eqref{eq:weight-trans}.  The
associated matrix $\Mat H_1(\Mat \Xi)$ in \eqref{eq:re-weight-H1}
leads to the inner product
\begin{equation*}
  \iProdn{\Vek u_1}{\Vek u_2}_{\mathcal H_1(\Mat \Xi)}
  = 
  \iProdn{\Vek u_1}{\Vek u_2}_{\mathcal H_1}
  -
  \sum_{k=0}^{S-1} \lambda_n \, 
  \iProdn{\Vek u_1}{\Vek \phi_k}_{\mathcal H_1}
  \iProdn{\Vek u_2}{\Vek \phi_k}_{\mathcal H_1},
\end{equation*}
where $S$ is the number of promoted directions.  For
$\iProdn{\cdot}{\cdot}_{\mathcal H_2(\Mat \Xi)}$, we obtain a similar
representation.

\begin{Algorithm}[reweighted tensor-free primal-dual for exact data]
  \label{alg:rew-tf-primal-dual:exact}
  \begin{enumerate}[(i)]
  \item Initiation: Fix the parameters $\tau, \sigma > 0$,
    $\theta \in [0,1]$, $\lambda \in [0,1)$, and
    $n_{\mathrm{rew}} > 0$.  Choose the start value
    $(\Vek w^{(0)}, \Vek y^{(0)}) = (\Vek 0 \otimes \Vek 0, \Vek 0)$
    in
    $(\mathcal H_1 \otimes \mathcal H_2) \times \mathcal
    K$,
    and set $\Vek w^{(-1)}$ to $\Vek w^{(0)}$.  Starting without
    weights, \ie\ $\Mat \Phi = []$ and $\Mat \Psi = []$.
  \item Iteration: For $n \ge 0$, update $\Vek w^{(n)}$ and $\Vek
    y^{(n)}$:
    \begin{enumerate}[(a)]
    \item Using the tensor-free computations in
      \thref{cor:tf-for-oper:bilin}, determine
      \begin{align*}
        \Vek y^{(n+1)}
        \coloneqq \Vek y^{(n)} + \sigma \, \bigl( (1 + \theta) \; 
        \breve{\mathcal B}(\Vek w^{(n)})
          -  \theta \; \breve{\mathcal B}(\Vek w^{(n-1)}) - \Vek g^\dagger \bigr).
      \end{align*}
    \item Compute a low-rank representation
      $\Vek w^{(n+1)} = \widetilde{\Mat V}^{(n+1)} \Mat \Sigma^{(n+1)}
      \widetilde{\Mat U}^{(n+1)}$
      of the singular value threshold
      \begin{equation*}
        \mathcal S_\tau ( \Vek w^{(n)} - \tau \, \breve{\mathcal
          B}^*_{\mathcal H_1(\Mat \Xi) \otimes \mathcal
          H_2(\Mat \Xi), \mathcal K}(\Vek y^{(n+1)}))
      \end{equation*}
      with \thref{alg:tf-threshold:lanczos}.  The required actions are
      given by \eqref{eq:sing-thres:right-act:rew} and
      \eqref{eq:sing-thres:left-act:rew}.
    \item Every $n_{\mathrm{rew}}$ iteration, re-weight the \PHilbert\
      spaces:
      \begin{itemize}
      \item  When $\Vek w^{(n+1)} = \Vek 0$, set
        $\Mat \Phi = []$ and $\Mat \Psi = []$. 
      \item Otherwise, use \thref{alg:aug-lanc-bidiag} to compute the
        (unweighted) singular value decomposition
        $\Vek w^{(n+1)} = \sum_{k=0}^{S-1} \sigma'_k \, (\Vek u'_k
        \otimes \Vek v'_k)$, i.e., with respect to
        $\mathcal H_1 \otimes \mathcal H_2$, where the needed actions
        are given in \eqref{eq:sing-thres:right-act:rew} and
        \eqref{eq:sing-thres:left-act:rew}.
      \item[] Set
        $\lambda_k \coloneqq \nicefrac{\lambda \sigma'_k}{\sigma'_0}$
        for $k=0,\ldots, S-1$, and
        $\Mat \Phi \coloneqq [\Vek u'_0, \dots, \Vek u'_{S-1}]$ as
        well as
        $\Mat \Psi \coloneqq [\Vek v'_0, \dots, \Vek v'_{S-1}]$.
      \end{itemize}
    \end{enumerate}
  \end{enumerate}
\end{Algorithm}

Analogously, one can apply the reweighting technique in
\thref{alg:rew-tf-primal-dual:exact}.ii.c to \thref{alg:tikh-reg} and
\ref{alg:inexact-data}.  All observations and heuristics in this
section also hold true or may be easily translated for the quadratic
setting.

\section{Masked \PFourier\ phase retrieval}
\label{sec:mask-phase-retr}

In this section, we apply the developed algorithm to the phase
retrieval problem.  Generally, the phase retrieval problem consists in
the recovery of an unknown signal from its \PFourier\ intensity.
Problems of this kind occur, for instance, in crystallography
\cite{Mil90, Hau91}, astronomy \cite{BS79, DF87}, and laser optics
\cite{SST04, SSD+06}.  To be more precise, in the following, we
consider the two-dimensional masked \PFourier\ phase retrieval problem
\cite{LCL+08,CSV13,CESV13,CLS15a,GKK17}, where the true signal
$\Vek u \in \C^{N_2 \times N_1}$ is firstly pointwise multiplied with
a set of known masks $\Vek d_\ell \in \C^{N_2 \times N_1}$ and
afterwards transformed by the two-dimensional $M_2 \times M_1$-point
\PFourier\ transform $\Fourier_{M_2 \times M_1}$ defined by
\begin{equation*}
  \bigl( \Fourier_{M_2 \times M_1}[\Vek v] \bigr)[m_2, m_1]
  \coloneqq \sum_{n_2 = 0}^{N_2 - 1} \sum_{n_1 = 0}^{N_1 - 1} \Vek
  v[n_2, n_1] \, \e^{-2\uppi\I(\nicefrac{n_2 m_2}{M_2} + \nicefrac{n_1 m_1}{N_1})}.
\end{equation*}
Denoting by $\odot$ the \pers{Hadamard} (or pointwise) product, 
the masked \PFourier\ phase retrieval problem can be
stated as follows.

\begin{Problem}[Masked \PFourier\ phase retrieval]
  \label{prob:mask-phase-retr}
  Recover the unknown complex-valued image
  $\Vek u \in \C^{N_2 \times N_1}$ from the masked \PFourier\
  intensities
  $\absn{\Fourier_{M_2 \times M_1}[\Vek d_\ell \odot \Vek u]}$ with
  $d = 0, \dots, L-1$.
\end{Problem}

In general, phase retrieval problems are ill-posed due to the loss of
the phase information in the frequency domain.  In one dimension, the
problem usually possesses an enormous set of non-trivial solutions,
which heavily differ from the true signal, see for instance
\cite{BP15} and references therein.  In the two-dimensional setting
considered by us, the situation changes dramatically since here almost
every signal can be uniquely recovered up to a global phase and up to
reflection and conjugation, see \cite{HM82,Hay82,BBE17}.

Before considering some numerical simulations, we study the quadratic
nature of the masked \PFourier\ phase retrieval problem, where we
employ the complex notation as mentioned in the introduction of
\autoref{sec:bilin-inverse-probl}.  For this purpose, we interpret
both, the domain $\C^{N_2 \times N_1}$ and the image
$\C^{L \times M_2 \times M_1}$ of the measurement operator in
\thref{prob:mask-phase-retr}, as real \PHilbert\ spaces.  To simplify
the notation, we vectorize the unknown image
$\Vek u \in \C^{N_2 \times N_1}$ and the given \PFourier\ intensities
$\absn{\Fourier_{M_2 \times M_1}[\Vek d_\ell \odot \Vek u]} \in
\C^{M_2 \times M_1}$ for a fixed mask columnwise.  Henceforth, the
vectorized variables are labeled with $\vec \cdot$.  On the
\PFourier\ side, we additionally attach the measurements for different
masks to each other.  Thus, the domain of the measurement operator in
\thref{prob:mask-phase-retr} becomes $\mathcal H = \C^{N_2 N_1}$ and
the image $\mathcal K = \C^{L M_2 M_1}$.  At the moment, the endowed
real inner product is not specified in detail.

In order to derive an explicit representation of the (vectorized)
forward operator, we write
the two-dimensional $(M_2 \times M_1)$-point \PFourier\ transform as
\begin{equation*}
  \Fourier_{M_2 \times M_1}[\Vek u] 
  = 
  \Mat F_{M_2} \Vek u \Mat F_{M_1}^\T
  =
  (\Mat F_{M_1} \otimes \Mat F_{M_2}) \, \vec{\Vek u},
\end{equation*}
with the \PFourier\ matrices
\begin{equation*}
  \Mat F_{M_2} 
  \coloneqq
  \Bigl(\e^{-2\uppi\I \, \frac{ n_2 m_2}{M_2}}\Bigr)_{m_2 = 0, n_2 =
    0}^{M_2-1,N_2-1}
  \qquad\text{and}\qquad
  \Mat F_{M_1} 
  \coloneqq
  \Bigl(\e^{-2\uppi\I \, \frac{ n_1 m_1}{M_1}}\Bigr)_{m_1 = 0, n_1 =
    0}^{M_1-1,N_1-1},
\end{equation*}
where $\otimes$ denotes the \PKronecker\ product of two matrices.  For
the vectorized version of the \PFourier\ transform
$\Fourier_{M_2 \times M_1}$, we henceforth use the notation
$\Mat F_{M_2 \times M_1} \coloneqq \Mat F_{M_2} \otimes \Mat F_{M_1}$.

In the same manner, we write the pointwise multiplication
$\Vek d_\ell \odot \Vek u$ as matrix vector multiplication
$\diag(\Vek d_\ell) \, \vec{\Vek u}$, where $\diag(\Vek d_\ell)$
denotes the matrix with diagonal $\Vek d_\ell$.  Combining the
interference with the given masks into one operator, we define the
matrix
\begin{equation*}
  \Mat D_L 
  \coloneqq
  \begin{pmatrix}
    \diag(\Vek d_1) \\
    \vdots \\
    \diag(\Vek d_L) \\
  \end{pmatrix}.
\end{equation*}
The action $\Mat D_L \vec{\Vek u}$ thus attaches the single masked
signals to each other.

Composing the two operations, and squaring the measurements, we notice
that \thref{prob:mask-phase-retr} is
equivalent to
\begin{equation}
  \label{eq:fourier-phase-retr}
  \text{recover} 
  \quad 
  \vec{\Vek u} \in \mathcal H
  \quad\text{from}\quad
  \absn{(\Mat I_L \otimes \Mat F_{M_2 \times M_1} ) \, \Mat D_L
    \vec{\Vek u})}^2 = \Vek g^\dagger,
  \tag{$\mathfrak{F}$}
\end{equation}
where $\Mat I_L \in \C^{L \times L}$ denotes the identity matrix, and
$\Vek g^\dagger \in \R^{L M_2 M_1}$ the vectorized exact squared, masked
\PFourier\ intensities of the looked for signal.  The associate
bilinear operator $\mathcal B_{\mathfrak F}$ of the quadratic forward
operator $\mathcal Q_{\mathfrak F}$ in \eqref{eq:fourier-phase-retr}
is now given by
\begin{equation}
  \label{eq:bilin-op-fourier}
  \begin{aligned}
    \mathcal B_{\mathfrak F}(\vec{\Vek u}, \vec{\Vek v}) 
    &\coloneqq
    \bigl[(\Mat I_L \otimes
    \Mat F_{M_2 \times M_1}) \, \Mat D_L \vec{\Vek v}\bigr]
    \odot
    \bigl[\overline{(\Mat I_L \otimes \Mat F_{M_2 \times M_1}) \, \Mat
      D_L \vec{\Vek u}}\bigr] 
    \\[\fskip]
    &=  \diag 
    \bigl(\bigl[(\Mat I_L \otimes
    \Mat F_{M_2 \times M_1}) \, \Mat D_L \vec{\Vek v}\bigr]
    \,
    \bigl[(\Mat I_L \otimes \Mat F_{M_2 \times M_1}) \, \Mat
    D_L \vec{\Vek u}\bigr]^* \bigr)
    \\[\fskip]
    &= \diag \bigl(
    (\Mat I_L \otimes \Mat F_{M_2 \times M_1}) \, \Mat D_L \vec{\Vek v} 
    \vec{\Vek u}^* \Mat D_L^* \, (\Mat I_L \otimes \Mat F_{M_2 \times
      M_1}^* ) \bigr),
  \end{aligned}
\end{equation}
where the function $\diag$ extracts the diagonal of a matrix.  Since
the last right-hand side is linear in $\vec{\Vek v} \vec{\Vek u}^*$ or
linear in $\vec{\Vek u} \vec{\Vek u}^*$ for
$\vec{\Vek v} = \vec{\Vek u}$, the quadratic lifting
$\breve{\mathcal Q}_{\mathfrak F} \colon \mathcal H
\otimes_{\mathrm{sym}} \mathcal H \rightarrow \mathcal K$ has to be
\begin{equation}
  \label{eq:quad-op-fourier}
  \breve{\mathcal Q}_{\mathfrak F}( \Vek w)
  = \diag \bigl(
  (\Mat I_L \otimes \Mat F_{M_2 \times M_1}) \,
  \Mat D_L \Vek w \Mat D_L^* \,
  (\Mat I_L \otimes \Mat F_{M_2 \times M_1}^*) \bigr).
  \addmathskip
\end{equation}
The last missing ingredient in order to apply our proximal algorithms
is the action of the adjoint $\breve{\mathcal Q}_{\mathfrak F}^*(\vec{\Vek
y})$ for a fixed $\vec{\Vek y} \in \mathcal K$.

\begin{Lemma}[Tensor-free adjoint lifting]
  \label{lem:tf-adj-lift:phase-retr}
  If the \PHilbert\ spaces $\mathcal H$ and $\mathcal K$ are endowed
  with the real \PEuclid{ian} inner product, \ie\ $\Mat H = \Mat I$
  and $\Mat K = \Mat I$, then the action of the
  adjoint operator $\breve{\mathcal Q}_{\mathfrak F}^*(\vec{\Vek y})$
  with fixed $\vec{\Vek y} \in \mathcal K$ is given by
  \begin{equation*}
    \breve{\mathcal Q}_{\mathfrak F}^*(\vec{\Vek y}) \,
    \vec{\Vek e} 
    =
    \Mat D_L^* \, (\Mat I_L \otimes \Mat F_{M_2 \times M_1}^*) \,
    \diag(\Re[\vec{\Vek y}]) \,
    (\Mat I_L \otimes \Mat F_{M_2 \times M_1}) \, \Mat D_L \vec{\Vek e}
  \end{equation*}
  for $\vec{\Vek e} \in \mathcal H$.
\end{Lemma}

\begin{Proof}
  We compute the action of the adjoint operator with the aid of
  \thref{lem:tens-free-adj-lift:quad}.  For this purpose, we first determine
  the left adjoint of $\mathcal B_{\mathfrak F}$ by considering
  \begin{align*}
    \iProdn{\mathcal B_{\mathfrak F}(\vec{\Vek f}, \vec{\Vek
    e})}{\vec{\Vek y}}
    &= 
      \iProdn{(\Mat I_L \otimes \Mat F_{M_2 \times M_1}) \, \Mat D_L \vec{\Vek e} 
    \vec{\Vek f}^* \Mat D_L^* \, (\Mat I_L \otimes \Mat F_{M_2 \times
      M_1}^* )}{ \diag(\vec{\Vek y})}
    \\[\fskip]
    &=
      \Re  \bigl[\tr\bigl(
      \diag(\vec{\Vek y})^* (\Mat I_L \otimes \Mat F_{M_2 \times M_1})
      \, \Mat D_L \vec{\Vek e}  
      \vec{\Vek f}^* \Mat D_L^* \, (\Mat I_L \otimes \Mat F_{M_2
      \times M_1}^* )
      \bigr) \bigr]
    \\[\fskip]
    &=
       \Re  \bigl[ \tr \bigl(
      \vec{\Vek f}^* \Mat D_L^* \, (\Mat I_L \otimes \Mat F_{M_2
      \times M_1}^* )
      \diag(\vec{\Vek y})^* (\Mat I_L \otimes \Mat F_{M_2 \times M_1})
      \, \Mat D_L \vec{\Vek e}  
      \bigr) \bigr]
    \\[\fskip]
    &= 
      \iProdn{\vec{\Vek f}}{ \Mat D_L^* \, (\Mat I_L \otimes \Mat F_{M_2
      \times M_1}^* )
      \diag(\vec{\Vek y})^* (\Mat I_L \otimes \Mat F_{M_2 \times M_1})
      \, \Mat D_L \vec{\Vek e}}
  \end{align*}
  for all $\vec{\Vek f} \in \mathcal H$ and fixed $\vec{\Vek e} \in
  \mathcal H$.  For the right adjoint, we analogously obtain
  \begin{equation*}
    \iProdn{\mathcal B_{\mathfrak F}(\vec{\Vek e}, \vec{\Vek
    f})}{\vec{\Vek y}}
    = 
    \iProdn{\vec{\Vek f}}{ \Mat D_L^* \, (\Mat I_L \otimes \Mat F_{M_2
    \times M_1}^* )
    \diag(\vec{\Vek y}) (\Mat I_L \otimes \Mat F_{M_2 \times M_1})
    \, \Mat D_L \vec{\Vek e}},
  \end{equation*}
  where diagonal in the middle is not conjugated.  Summation of the
  left and right adjoint as in \thref{lem:tens-free-adj-lift:quad}
  yields the assertion.  \qed
\end{Proof}

\begin{Remark}
  \label{rem:tf-adj-lift:phase-retr}
  Using \thref{lem:adj-oper}, we can now transform the actions of the
  \PEuclid{ian} adjoint to our actual \PHilbert\ spaces.  More
  precisely, we have
  \begin{equation*}
    \breve{\mathcal Q}_{\mathfrak F}^*(\vec{\Vek y}) \, \Mat H
    \vec{\Vek e} 
    =
    \Mat H^{-1}
    \Mat D_L^* \, (\Mat I_L \otimes \Mat F_{M_2 \times M_1}^*) \,
    \diag(\Re[\vec{\Mat K\Vek y}]) \,
    (\Mat I_L \otimes \Mat F_{M_2 \times M_1}) \, \Mat D_L \vec{\Vek
      e}.
    \tag*{\qed}
  \end{equation*}
\end{Remark}

One of the central reasons to choose the masked \PFourier\ phase
retrieval problem as application of the proposed algorithms and
heuristics is that the phase retrieval
problem~\eqref{eq:fourier-phase-retr}, although severely ill posed,
behaves nicely under convex relaxation.  More precisely, one can show
that under certain conditions the solution of the convex relaxed
problem
\begin{equation*}
  \minimize \quad \pNormn{\Vek w}^+_{\mathcal H \otimes_\uppi \mathcal
    H}
  \quad\text{subject to}\quad
  \breve{\mathcal Q}_{\mathfrak F}(\Vek w) = \Vek g^\dagger
\addmathskip
\end{equation*}
is unique and coincides with the true lifted solution
$\vec{\Vek u} \vec{\Vek u}^*$ with high probability, see \cite{CSV13,
  CESV13,CLS15a,GKK17}.  Therefore, we expect that the proposed
tensor-free proximal algorithms converge to a unique rank-one
solution.

\subsection{Effects of bidiagonalization and reweighting}
\label{sec:effects-bidi-re}

In the first numerical example, we want to study the effect of the
applied augmented \PLanczos\ bidiagonalization and of the reweighting
heuristic to the computation time and the convergence behavior.  For
all simulations in this paper, the proposed methods have been
implemented in MATLAB\textsuperscript{\textregistered} (R2017a,
64-bit) and are performed using an
Intel\textsuperscript{\textcopyright} Core\texttrademark\ i7-4790 CPU
(4$\times$3.60~GHz) and 32~GiB main memory.  The employed true two-dimensional
signal consists in a synthetic image referring to transmission
electron microscopy experiments with nanoparticles.  The test image
$\Vek u$ is rather small and is composed of $16 \times 16$ pixels.
The corresponding tensor $\Vek w \coloneqq \Vek u \Vek u^*$ is already
of dimension $256$.

\begin{figure}
  \centering
  \subfloat[Recovered signal.]  {\includegraphics{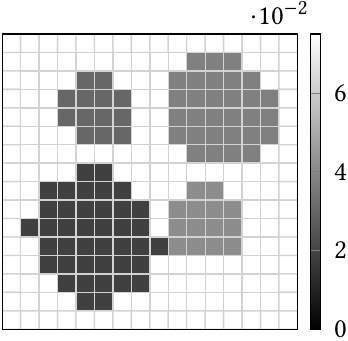}}
  \qquad
  \subfloat[Absolute error.]
  {\includegraphics{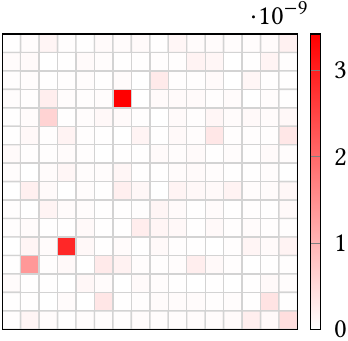}}
  \caption{Masked phase retrieval based on
    Algorithm~\ref{alg:exact-data} with 1\,000 iterations and without
    any modification.  The eight masks have been chosen with respect
    to a \PRademacher\ distribution.  Each pixel is covered by at
    least one mask.}
  \label{fig:no-bidi-rw}
\end{figure}

Based on the true image, we compute synthetic and noise-free data by
applying the masked \PFourier\ transform.  The entries of the eight
employed masks have been randomly generated with respect to
independent \PRademacher\ random variables.  More precisely, the
entries of the masks are distributed with respect to the model
\begin{equation}
  \label{eq:rademacher-mask}
  \Vek d_\ell [n_2, n_1]
  \sim
  \begin{cases}
    \sqrt2 & \text{with probability $\nicefrac14$,} \\
    0 & \text{with probability $\nicefrac12$,} \\
    - \sqrt2 & \text{with probability $\nicefrac14$.} \\
  \end{cases}
\end{equation}
Masks of this kind have been studied in \cite{CLS15a,GKK17} in order
to de-randomize the generic phase retrieval problem considered in
\cite{CSV13,CESV13}.  In our experiment, we employ the
$32 \times 32$-point \PFourier\ transform such that the complete
autocorrelation of the masked signals is encoded in the given
\PFourier\ intensities.  A first reconstruction of the true signal
based on \thref{alg:exact-data} is shown in \autoref{fig:no-bidi-rw},
where we compute the singular value threshold with the aid of a full
singular value decomposition of the tensor $\Vek w^{(n)}$.  For the
involved \PHilbert\ spaces $\mathcal H$ and $\mathcal K$, we simply
choose the corresponding \PEuclid{ian} spaces.  In other words, we
employ the \PHilbert--\PSchmidt\ inner product for the (vectorized)
matrices in $\mathcal H = \C^{N_2 N_1}$.

Although the reconstruction is quite accurate, the main drawbacks of a
direct application of \thref{alg:exact-data} are the computation time
and memory requirements.  For 1\,000 iterations, we need about
6.38~minutes to recover the true signal.  Next, we employ the
tensor-free variant of the primal-dual iteration in
\thref{alg:tensor-free}, where we apply the augmented \PLanczos\
bidiagonalization to determine the singular value thresholding
(\thref{alg:tf-threshold:lanczos}).  Using this modification, we
merely need about 58 seconds to perform the reconstruction.  Since we
can control the accuracy of the partial singular value decomposition,
the performed iteration essentially coincides with the original
iteration.

The influence of the augmented \PLanczos\ bidiagonalization on the
computation time is presented in \autoref{tab:time-saving}.  The
parameter $k$ here describes the maximal size of the bidiagonal matrix
$\Mat B_k$.  Further, $\ell$ denotes the number of fixed \PRitz\ pairs
in the augmented restarting procedure.  Since the approximation
property of the incomplete \PLanczos\ method becomes worse for small
$k$, we require more restarts in order to observe an accurate partial
singular value decomposition.  For higher-dimensional input images,
the time saving aspect becomes much more important.  

\begin{table}
  \centering
  \small
  \begin{tabular}{|l|ccccc|c|}
    \cline{7-7}
    \multicolumn{6}{c|}{} & \tiny reweighting\\
    \hline
    \PLanczos\ vectors $k$ & $\infty$ & 100 & 50 & 20 & 10 & 10 \\
    iteration vectors $\ell$ & --- & 50 & 25 & 10 & 5 & 5 \\
    \hline
    time in seconds & 383.23 & 139.51 & 67.94  & 60.03 & 57.62 & 35.74 \\
    average restarts & --- & 0.00 & 0.06 & 2.89 & 10.01 & 4.89 \\
    \lasthline
  \end{tabular}
  \caption{Required computation time for the reconstruction of a $16
    \times 16$ image by Algorithm~\ref{alg:tensor-free} with
    augmented \PLanczos\ process and performing
    1\,000 iterations.}
  \label{tab:time-saving}
\end{table}

Considering the evolution of the non-zero singular values of
$\Vek w^{(n)}$ during the iteration, see \autoref{tab:evo-sval}, we
observe that the projective norm heuristic here enforces a very low
rank.  After 500 iterations, we nearly obtain a rank-one tensor such
that the additional reconstruction error caused by the rank-one
projection of $\Vek w^{(n)}$ to extract the recovered image becomes
negligible.  After 2\,000 iterations, the tensor $\Vek w^{(n)}$ has
converged to a rank-one tensor.  

\begin{table}
  \centering
  \small
  \begin{tabular}{|r|*{3}{r@{\,$\cdot$\,}l}|}
    \firsthline
    \multicolumn{1}{|c|}{$n$} 
    & \multicolumn{2}{c}{$\sigma_0$}
    & \multicolumn{2}{c}{$\sigma_1$}
    & \multicolumn{2}{c|}{$\sigma_2$} 
    \\
    \hline 
    $100$ 
    & $9.89$ & $10^{-1}$ 
    & $2.00$ & $10^{-3}$ 
    & $6.02$ & $10^{-4}$
    \\
    $500$
    & $1.00$ & $10^{\pm 0}$
    & $7.06$ & $10^{-7}$
    & $4.59$ & $10^{-7}$
    \\
    $1\,000$
    & $1.00$ & $10^{\pm 0}$
    & $1.13$ & $10^{-10}$
    & \multicolumn{2}{c|}{}
    \\
    $2\,000$
    & $1.00$ & $10^{\pm 0}$
    & \multicolumn{2}{c}{}
    & \multicolumn{2}{c|}{}
    \\
    \lasthline
  \end{tabular}
  \caption{Evolution of the non-zero singular values using
    Algorithm~\ref{alg:tensor-free} with  augmented
    \PLanczos\  process.}
  \label{tab:evo-sval}
\end{table}

In order to promote the rank-one solutions of the masked \PFourier\
phase retrieval problem even further, we next exploit the reweighting
approach proposed in \autoref{sec:reduc-rank-adapt}.  For our current
simulation, this means that we replace the inner product matrices
$\Mat H = \Mat I$ by the modified matrices $\Mat H(\Mat \Xi)$ defined
in \eqref{eq:re-weight-H1}, where we build up the basis
$\{\Vek \phi_{k}\}$ from the singular value decomposition of
$\Vek w^{(n)}$.  In the \PEuclid{ian} setting considered in this
simulation, the transformed directions $\widetilde{\Vek \phi}_k$
simply coincide with $\Vek \phi_{k}$. The reweighting is here applied
every ten iterations with relative weight $\lambda \coloneqq \nicefrac12$.

\begin{figure}
  \centering
  \subfloat[Recovered signal]
  {\includegraphics{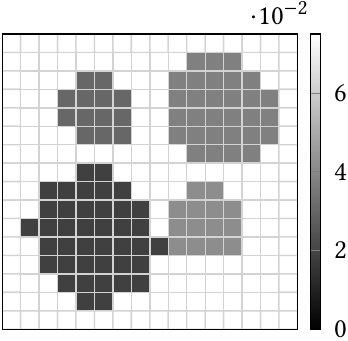}}
  \qquad
  \subfloat[Absolute error]
  {\includegraphics{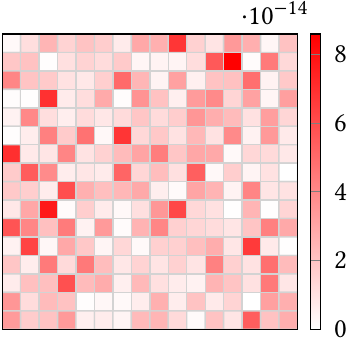}}
  \caption{Masked phase retrieval based on
    Algorithm~\ref{alg:rew-tf-primal-dual:exact} incorporating the
    augmented \PLanczos\ process and the reweighting heuristic.  The
    algorithm is terminated after 1\,000 iterations, the reweighting
    is repeated every 10 iterations with $\lambda \coloneqq \nicefrac12$.}
  \label{fig:with-bidi-rw}
\end{figure}

\begin{figure}
  \centering
  \subfloat[Evolution of the tensor rank.]
  {\includegraphics{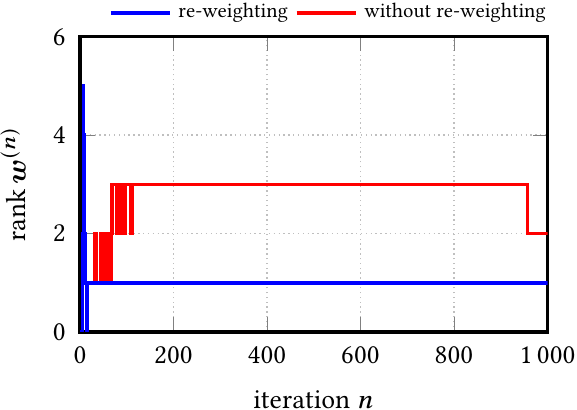}}
  \qquad
  \subfloat[Evolution of the data fidelity term.]
  {\includegraphics{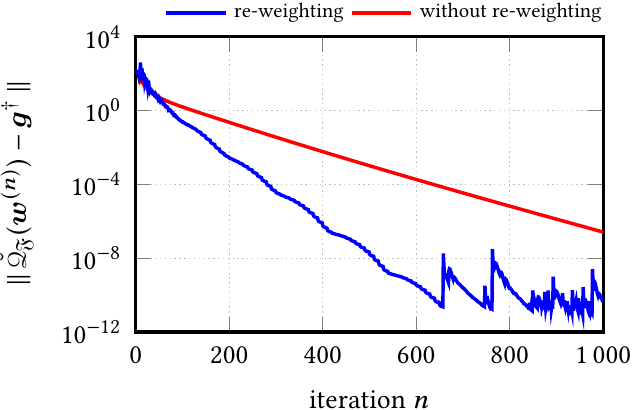}}
  \caption{Evolution of the rank and data fidelity term during the
    masked phase retrieval problem based on
    Algorithm~\ref{alg:exact-data} and
    Algorithm~\ref{alg:rew-tf-primal-dual:exact}. }
  \label{fig:rank-data-fid}
\end{figure}

The results of the tensor-free masked \PFourier\ phase retrieval with
augmented \PLanczos\ process and \PHilbert\ space reweighting
(\thref{alg:rew-tf-primal-dual:exact}) are shown in
\autoref{fig:with-bidi-rw}.  Although the reconstruction looks
comparable, we want to point out that the absolute errors are several
magnitudes smaller.  If we compare the evolution of the ranks,
see~\autoref{fig:rank-data-fid}, we can see that the proposed
reweighting heuristic reduces the rank quite efficiently.  More
precisely, most of the iterations have rank one.  Due to this
reduction, the reweighting has also a positive effect on the
computation time and the average number of restarts of the \PLanczos\
process, see \autoref{tab:time-saving}.  Moreover, we may notice that
the data fidelity term
$\pNormn{\breve{\mathcal Q_{\mathfrak F}}(\Vek w^{(n)}) - \Vek
  g^\dagger}$ decreases with a higher rate.  Here the convergences
stops after about 650 iterations due to numerical reasons.

\subsection{Incorporating smoothness properties}
\label{sec:incorp-smoothn-prop}

One of the central difference between our tensor-free reweighting algorithm
and Phase\-Lift \cite{CSV13,CESV13} consists in the modeling of the
domain and image space of the phase retrieval problem.  Where
PhaseLift is based on the standard \PEuclid{ian} setting, we rely on
arbitrary discrete \PHilbert\ spaces $\mathcal H$ and $\mathcal K$.
Especially, in two-dimensional phase retrieval, we may thus exploit
relationships between neighbored pixels like finite differences.
More precisely, we here study the influence of the a-priori smoothness
property formulated in terms of the two-dimensional discretized
\PSobolev\ norm.

In order to discretize the (weighted) \PSobolev\ space $W^{1,2}$, we
employ the forward differences
\begin{align*}
  \uppartial_1 \Vek u[n_2,n_1] &\coloneqq \Vek u [n_2, n_1 +1] - \Vek u[n_2,n_1]
  \qquad\text{for}\qquad
  \left\{
  \begin{aligned}
    n_1&=0, \dots, N_1-2\\n_2 &= 0, \dots, N_2-1
  \end{aligned}
  \right\}
  \\
  \shortintertext{and}
  \uppartial_2 \Vek u[n_2,n_1] &\coloneqq \Vek u [n_2 +1, n_1] - \Vek u[n_2,n_1]
  \qquad\text{for}\qquad
  \left\{
  \begin{aligned}
    n_1&=0, \dots, N_1-1\\n_2 &= 0, \dots, N_2-2
  \end{aligned}
  \right\}
\end{align*}
to approximate the first partial derivatives.  The associate linear
mappings for the vectorized image $\vec{\Vek u}$ are here denoted by
$\Mat D_1$ and $\Mat D_2$.  Thus the weighted \PSobolev\ space $\mathcal H =
W^{1,2}_{\Vek \mu}$ with weights $\Vek \mu \coloneqq (\mu_{\Mat I},
\mu_{\Mat D_1}, \mu_{\Mat D_2})^*$ corresponds to the matrix
\begin{equation*}
  \Mat H \coloneqq \mu_{\Mat I} \, \Mat I + \mu_{\Mat D_1} \, \Mat D_1^*
  \Mat D_1 + \mu_{\Mat D_2} \Mat D_2^* \Mat D_2.
\end{equation*}
In comparison, the discretized $L^2$-norm corresponds to the standard
\PEuclid{ian} setting associated with the identity matrix.

\begin{figure}
  \centering
  \subfloat[True signal.]
  {\includegraphics{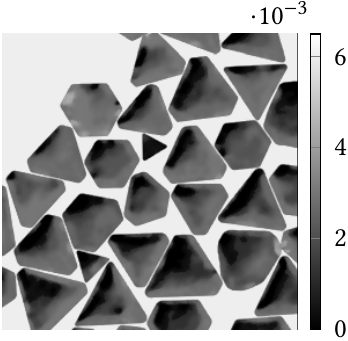}}
  \qquad
  \subfloat[Number of covering \PRademacher\  masks.]
  {\includegraphics{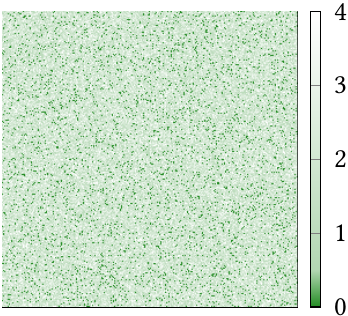}}
  \caption{The \PFourier\ data for the second experiment
    ($256 \times 256$ pixels) have been created on the basis of the TEM
    micrograph of gold nanoparticles \cite[Figure~6C]{LSP+16}.  The
    employed masks have again been chosen with respect to the
    \PRademacher\ distribution.  In this instance, about a sixteenth
    of all pixels are blocked in each of the four masks.}
  \label{fig:smooth:true-sol}
\end{figure}

The masked \PFourier\ intensities of the second example has been
created on the basis of a transmission electron microscopy (TEM)
reconstruction in \cite{LSP+16}.  The image has a dimension of
$265 \times 265$ pixels such that the related tensor possesses
$2^{32}$ complex-valued entries and requires 64~GiB memory (double
precision complex numbers).  Further, we apply four random masks of
\PRademacher-type \eqref{eq:rademacher-mask}.  Since the masks are
generated entirely random, about a one-sixteenth of the pixel are
blocked by all masks.  The test image together with the number of
masks covering a certain pixel are shown in
\autoref{fig:smooth:true-sol}.

\begin{figure}
  \centering
  \subfloat[Recovered signal.]
  {\includegraphics{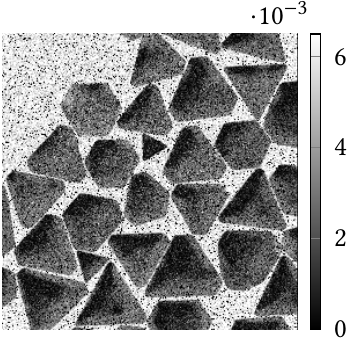}}
  \qquad
  \subfloat[Absolute error.]
  {\includegraphics{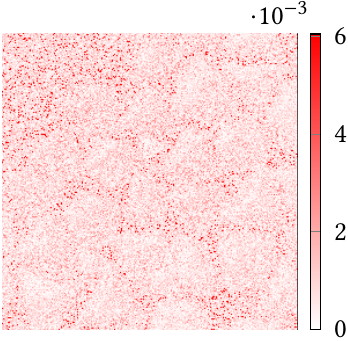}}
  \caption{Masked phase retrieval based on
    Algorithm~\ref{alg:rew-tf-primal-dual:exact}.  The pre-image space
    $\C^{N_2 N_1}$ is equipped with the discrete $L^2$ inner product.
    The reconstruction is terminated after 100 iterations.  In order
    to compare the retrieval with the true signal, the pixels are
    presented in the same range as the true image, resulting in the
    truncation of higher intensities.}
  \label{fig:smooth:eucl}
\end{figure}

\begin{figure}
  \centering
  \subfloat[Recovered signal.]
  {\includegraphics{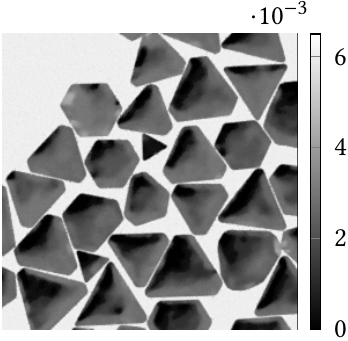}}
  \qquad
  \subfloat[Absolute error.]
  {\includegraphics{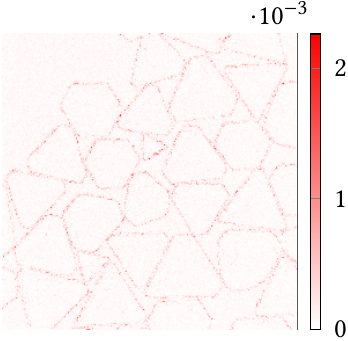}}
  \caption{Masked phase retrieval based on
    Algorithm~\ref{alg:rew-tf-primal-dual:exact}.  The pre-image space
    $\C^{N_2 N_1}$ is equipped with the discrete \PSobolev\ inner
    product based on the weight $\Vek \mu = (0.25, 1.00, 1.00)^*$.
    The reconstruction is terminated after 100 iterations.}
  \label{fig:smooth:sobo}
\end{figure}

To solve the corresponding masked phase retrieval problem, we apply
\thref{alg:rew-tf-primal-dual:exact}, where we reweight the \PHilbert\
spaces every 10 steps with the relative weight
$\lambda \coloneqq \nicefrac12$.  Since the algorithm tends to
higher-rank tensors in the starting phase, we only compute a partial
singular value decomposition with at most five leading singular values
using ten \PLanczos\ vectors $\Vek e_n$ and $\Vek f_n$.  Hence, we
perform \thref{alg:tensor-free} in an inexact manner.  After a few
iterations the rank of $\Vek w^{(n)}$ decreases such that the method
becomes again exact, and that the convergence is ensured.

\begin{figure}
  \centering
  \subfloat[Evolution of the tensor rank.]
  {\includegraphics{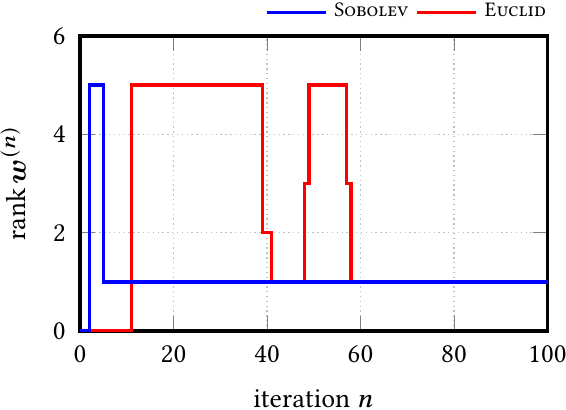}}
  \qquad
  \subfloat[Evolution of the data fidelity term.]
  {\includegraphics{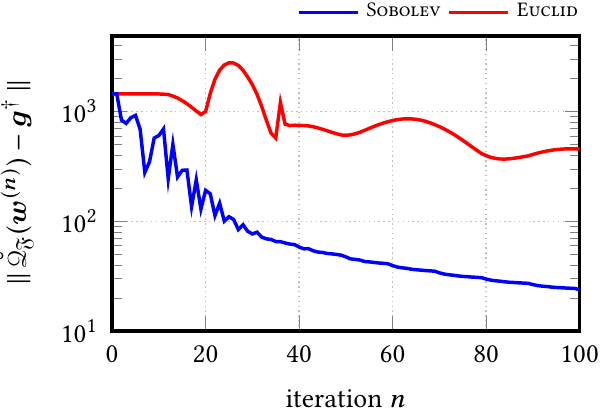}}
  \caption{Evolution of the rank and data fidelity term during the
    masked phase retrieval based on
    Algorithm~\ref{alg:rew-tf-primal-dual:exact}.  Both terms are compared for the
    discrete $L^2$ norm (\PEuclid) and the discrete \PSobolev\ norm
    with weight $\Vek \mu = (0.25, 1.00, 1.00)^*$.}
  \label{fig:smooth:rank-data-fid}
\end{figure}

The reconstructions for the \PEuclid{ian} and \PSobolev\ setting are
presented in \autoref{fig:smooth:eucl} and \ref{fig:smooth:sobo},
respectively.  Due to the small number of four masks, the convergence
of the algorithm using the discretized $L^2$-norm is very problematic.
Although the method converges for the chosen parameters, the
convergence rate is very low.  Moreover, pixels that are not covered
by any masks cannot be recovered and cause reconstruction defects
characterized by black holes.  Using instead the discretized
\PSobolev\ norm with weight $\Vek \mu \coloneqq (\nicefrac14,1,1)$ and
the same parameters, we obtain a much faster convergence and rank
reduction, see \autoref{fig:smooth:rank-data-fid}.  Further, the
required number of dual updates in order to produce a non-zero primal
update is reduced.  A small drawback is that the \PSobolev\ norm tends
to smooth out the edges of the particles in the reconstruction.  One
the other side, the a-priori smoothness condition allows us to recover
pixels not covered by the given data.

\subsection{Phase retrieval for large-scale images}
\label{sec:phase-retrieval-high}

\begin{figure}
  \centering
  \subfloat[Recovered signal.]
  {\includegraphics{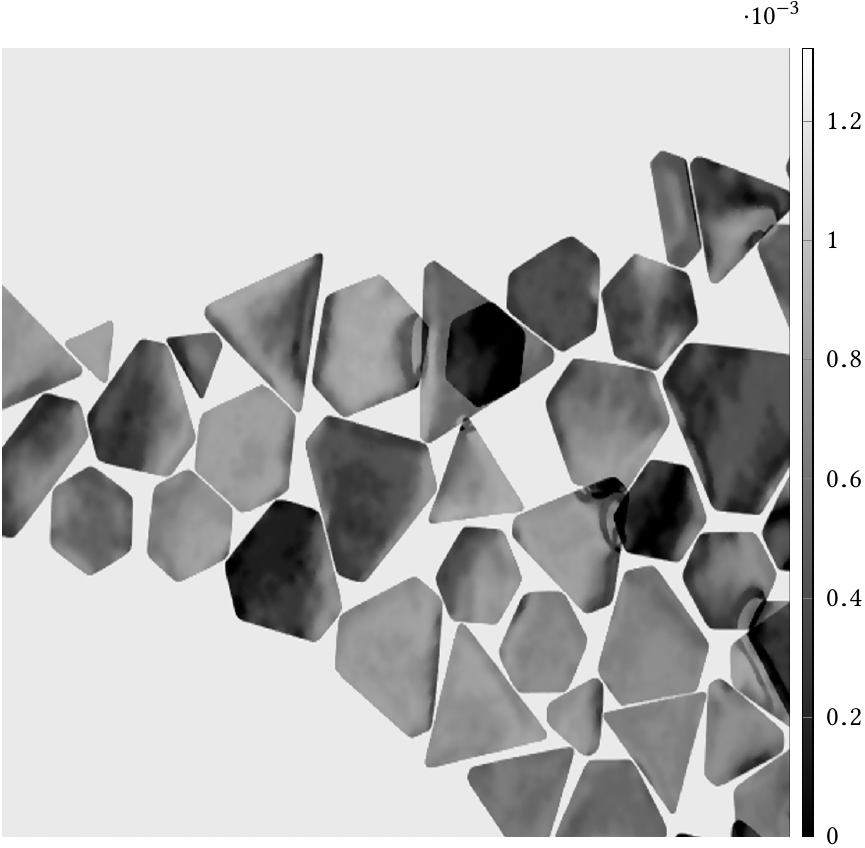}}
  \qquad
  \subfloat[Absolute error.]
  {\includegraphics{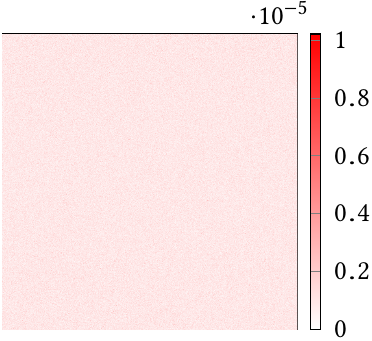}}
  \caption{The \PFourier\ data for the third experiment
    ($1\,024 \times 1\,024$ pixels) have been created on the basis of
    the TEM micrograph \cite[Figure~6B]{LSP+16}.  The eight masks have
    been generated regarding a \PGauss{ian} distribution.
    Algorithm~\ref{alg:rew-tf-primal-dual:exact} has been terminated
    after 1\,000 iterations.}
  \label{fig:high-dim:rec-sol}
\end{figure}

\begin{figure}
  \centering
  \subfloat[Evolution of the tensor rank.]
  {\includegraphics{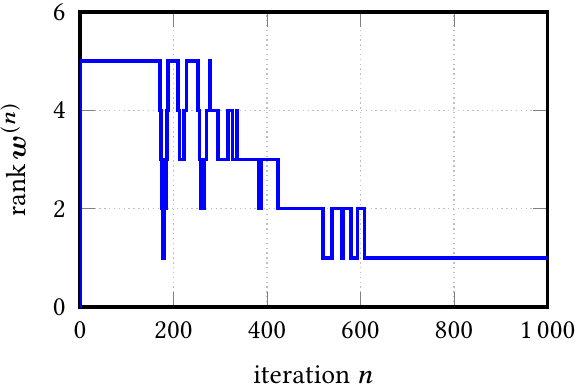}}
  \qquad
  \subfloat[Evolution of the data fidelity term.]
  {\includegraphics{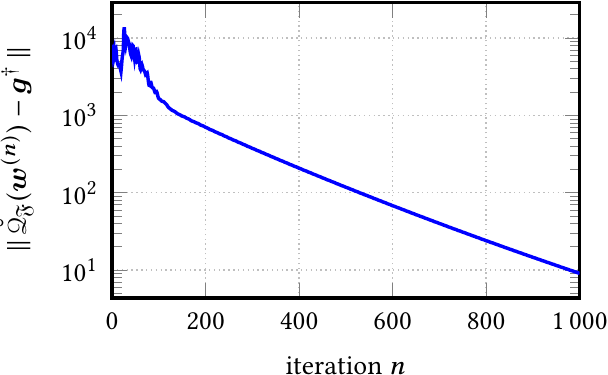}}
  \caption{Evolution of the rank and data fidelity term during the
    masked phase retrieval problem based on
    Algorithm~\ref{alg:rew-tf-primal-dual:exact} using eight
    \PGauss{ian} masks and terminating after 1\,000 iterations.}
  \label{fig:high-dim:rank-data-fid}
\end{figure}

Using the proposed reweighting heuristic to reduce the rank of the
iteration $\Vek w^{(n)}$, we are able to perform
\thref{alg:rew-tf-primal-dual:exact} for much larger test images.  In
this numerical experiment, we consider an $1\,024 \times 1\,024$ pixel
image, whose \PFourier\ data are again based on a TEM micrograph of
gold nanoparticles \cite{LSP+16}.  The lifted image here already
requires 16~TiB memory in order to hold the $2^{40}$ complex-valued
entries with double precision.  Differently from the previous
examples, we here apply eight \PGauss{ian} masks following the
standard normal distribution
\begin{equation*}
  \Vek d_\ell[n_2, n_1] \sim \mathcal N(0,1).
  \addmathskip
\end{equation*}
The recovered signal for the \PEuclid{ian} inner product is
shown in \autoref{fig:high-dim:rec-sol} together with the evolution of
the rank and data fidelity in \autoref{fig:high-dim:rank-data-fid}.
Analogously to the above experiments, the \PHilbert\ space $\mathcal
H$ is reweighted every ten iterations with a relative weight $\lambda
\coloneqq \nicefrac12$.

\subsection{Corruption by noise}
\label{sec:corruption-noise}

In the last numerical example, we study the influence of noise to the
proposed tensor-free primal-dual algorithm.  For simplicity, we only
study the behavior of the proposed method with respect to white or
\PGauss{ian} noise of the form
$\Vek g^\epsilon \coloneqq \Vek g^\dagger + \Vek \zeta$ where
$\Vek \zeta$ is a normal distributed random vector.  For the noise
level $\epsilon \coloneqq \pNormn{\Vek g^\dagger - \Vek g^\epsilon}$,
we consider different percentages of the norm
$\pNormn{\Vek g^\dagger}$.  Similarly to the first numerical examples,
we again apply four \PRademacher-type masks of the form
\eqref{eq:rademacher-mask}.  The synthetic data
$\Vek g^\dagger$ for the $256 \times 256$ test image are again based
on a TEM reconstruction of gold nanoparticles \cite[Figure~S1B]{LSP+16a} and the
$512 \times 512$-point \PFourier\ transform.  The domain $\mathcal H$
is again equipped with the discretized \PSobolev\ norm with weight
$\Vek \mu \coloneqq (\nicefrac14, 1, 1)^*$.  Due to the noise, we have
adapted \thref{alg:rew-tf-primal-dual:exact} to the \PTikhonov\
regularization in \thref{alg:tikh-reg}, which simply means that we
have to multiply $\Vek y^{(n+1)}$ in step (ii.a) with
$\nicefrac{1}{\sigma+1}$ additionally and to scale the threshold in
step (ii.b) to $\mathcal S_{\tau \alpha}$.  Since $\alpha$ affects the
influence of the projective norm heuristic, this parameter has to be
chosen relatively large.  In this brief test run with noisy
measurements, we chose $\alpha \coloneqq 10^3$ independent of the
noise level.  Surely, the results can be improved by more sophisticated
parameter choice rules.

\begin{figure}
  \centering
  \begin{tabular}{l@{\qquad}l@{\qquad}l}
    \subfloat[Reconstruction (5\%).]
    {\includegraphics{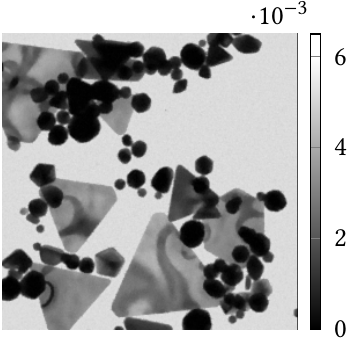}\;\:}
    &
      \subfloat[Reconstruction (10\%).]
      {\includegraphics{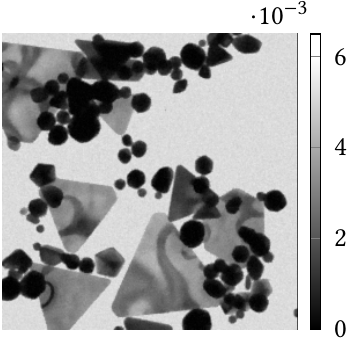}\;\:}
    &
      \subfloat[Reconstruction (25\%).]
      {\includegraphics{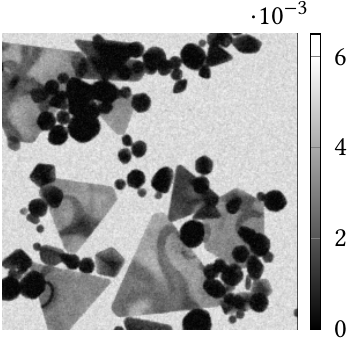}\;\:}
    \\
    \subfloat[Absolute error (5\%).]
    {\includegraphics{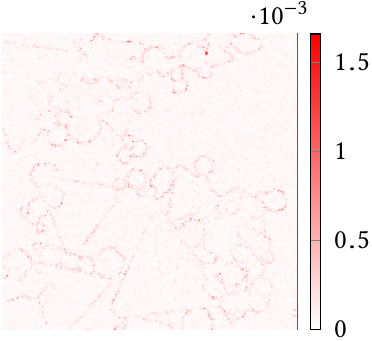}}
    &
      \subfloat[Absolute error (10\%).]
      {\includegraphics{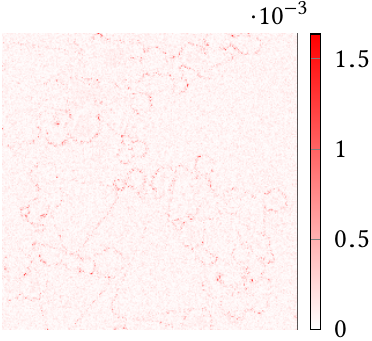}}
    &
      \subfloat[Absolute error (25\%).]
      {\includegraphics{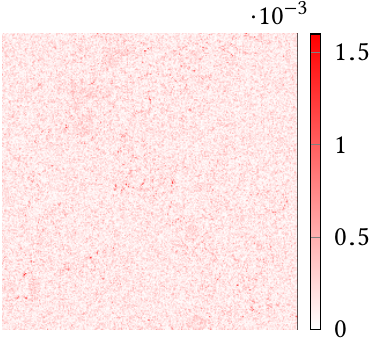}}
  \end{tabular}
  \caption{The influence of noise to masked phase retrieval for a
    $256 \times 256$ test image based on the TEM reconstruction of
    gold nanoparticles \cite[Figure~S1B]{LSP+16a}.  The adaption of
    Algorithm~\ref{alg:rew-tf-primal-dual:exact} with respect to the
    \PTikhonov\ regularization in Algorithm~\ref{alg:tikh-reg} is
    terminated after 100 iterations.}
  \label{fig:noise:rec-sol}
\end{figure}

\begin{figure}
  \centering
  \subfloat[Evolution of the tensor rank.]
  {\includegraphics{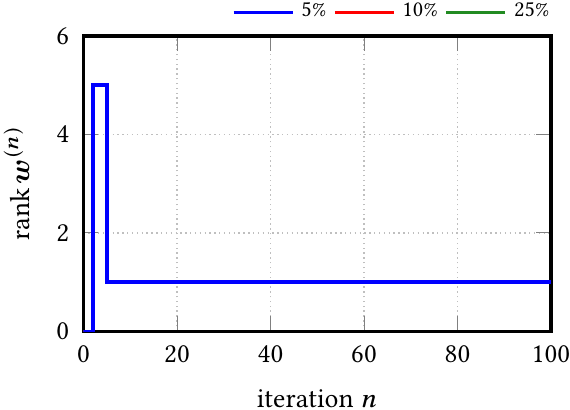}}
  \qquad
  \subfloat[Evolution of the data fidelity term.]
  {\includegraphics{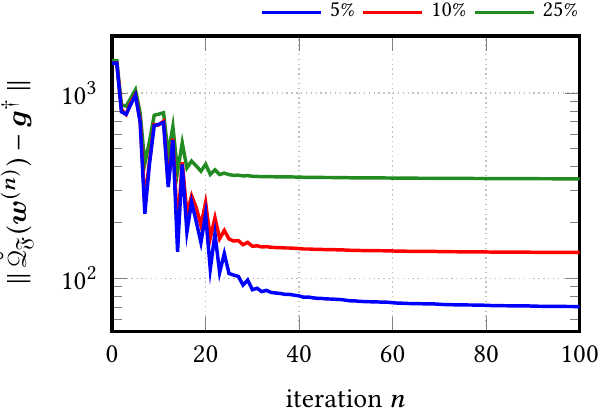}}
  \caption{Evolution of the rank and data fidelity term during masked
    phase retrieval for noisy measurements.  The noise level is given
    with respect to the norm $\pNormn{\Vek g^\dagger}$ of the
    noise-free measurements.}
  \label{fig:noise:rank-data-fid}
\end{figure}

The recovered signals and the evolution of the rank and data fidelity
are shown in \autoref{fig:noise:rec-sol} and
\ref{fig:noise:rank-data-fid}, respectively.  Due to the noise, we
cannot ensure that the recovered tensor has rank one and corresponds
to a meaningful approximation of the true signal.  If we endow the
domain $\mathcal H$ with the \PEuclid{ian} inner product in analogy to
classical PhaseLift, then we observed that the rank of the iteration
$\Vek w^{(n)}$ increases uncontrollable, and the proposed algorithm
diverges because of the limited data provided by only four masks.
Using instead the weighted \PSobolev\ norm and the \PHilbert\ space
reweighting, the rank becomes one after a short starting phase, where
the maximal rank is restricted by five.  Because of the same starting
value and the same regularization parameter, the primal-dual iteration
initially performs nearly identical for the three considered noise
levels such that the rank evolutions coincide. Further, for all three
cases, we here obtain reasonable reconstructions.  Moreover, the
pixels not covered by any mask are filled up, and the influence of the
noise to the reconstruction is smoothed out.  Further numerical
experiments suggest that \thref{alg:rew-tf-primal-dual:exact} in
combination with a-priori smoothness properties recovers the unknown
signal in a stable manner.

\section{Conclusion}
\label{sec:conclusion}

In this paper, we developed a novel proximal algorithm to solve
bilinear and quadratic inverse problems.  The basic idea was to
exploit the universal property to lift the considered problem to
linear inverse problems on the tensor product.  In order to deal with
the rank-one constraint, we applied a nuclear or projective norm
heuristic, which is known to produce low-rank solutions.  The
relaxation of the lifted problem yields a constrained minimization
problem, which has been solved by applying the first-order primal-dual
algorithm proposed by \pers{Chambolle} and \pers{Pock}.  Since the
choice of the underlying algorithm is somehow arbitrary, there are
several further options to develop minimization methods for the lifted
and relaxed problem.  For instance, one may apply the alternating
direction method of multipliers (ADMM) \cite{BPC+11},
forward-backward splitting \cite{LM79,CW05}, or FISTA \cite{BT09}.

The flexibility to adapt the domain of the bilinear or
quadratic operator allows us to incorporate smoothness assumptions or
neighborhood relations.  As demonstrated for the masked \PFourier\
phase retrieval problem, this freedom enables us to recover pixels that
are blocked by the applied masks such that they do not contribute to
the given measurements.  Further, the smoothing properties of the
discretized \PSobolev\ norm greatly improve the numerical observed
convergence rates.  Moreover, one can exploit this flexibility to
reweight the pre-image spaces in order to promote low-rank
solutions.  In the moment, we rely on a projective norm based on the
corresponding \PHilbert\ norms.  Here, the question arises if one can
employ nuclear ``norms'' based on semi-norms like the total variation or
total generalized variation.

For the masked \PFourier\ phase retrieval problem, we have shown that
the developed algorithms seem to be stable under noise.  More
precisely, we have studied the influence of white noise, which can be
treated by choosing the squared \PEuclid{ian} or discretized
$L^2$-norm as data fidelity term of the \PTikhonov\ functional in
\eqref{eq:tikh-func-bilin} and \eqref{eq:tikh-func-quad}.  In order to
incorporate more realistic noise models like \PPoisson\ noise into
phase retrieval, one can, for instance, replace the data fidelity by
the \pers{Kullback}--\pers{Leibler} divergence.  In so doing, one only
has to update the proximal mapping to compute the dual variable.

\bigskip\noindent
\parbox{\linewidth}{\footnotesize {\sffamily\bfseries
    Acknowledgments.} \quad We gratefully acknowledges the funding of
  this work by the Austrian Science Fund (FWF) within the project
  P28858.  The Institute of Mathematics and Scientific Computing off
  the University of Graz, with which the authors are affiliated, is a
  member of NAWI Graz (https://www.nawigraz.at/).}

{\footnotesize 
\newcommand{\etalchar}[1]{$^{#1}$}
}

\end{document}